\documentclass[12pt]{amsart}
\usepackage[utf8]{inputenc}
\usepackage[margin=1in]{geometry}
\usepackage[all]{xy}
\usepackage{tikz}
\usepackage{tikz-cd}
\usepackage{enumitem}
\usepackage{amssymb}
\usepackage{url}
\usepackage[colorlinks,allcolors=blue]{hyperref}
\newtheorem{theorem}[subsection]{Theorem}
\newtheorem{lemma}[subsection]{Lemma}
\newtheorem{proposition}[subsection]{Proposition}
\newtheorem{corollary}[subsection]{Corollary}

\theoremstyle{definition}
\newtheorem{definition}[subsection]{Definition}
\newtheorem{construction}[subsection]{Construction}

\theoremstyle{remark}
\newtheorem{example}[subsection]{Example}
\newtheorem{remark}[subsection]{Remark}

\newtheorem{warning}[subsection]{Warning}

\makeatletter
\let\c@figure\c@subsection

\makeatother

\newcommand{\Top}{\mathsf{Top}_*}

\newcommand{\Spt}{\mathsf{Spt}}
\newcommand{\TQ}{\mathsf{TQ}}
\newcommand{\Sp}{\Sigma^{\infty}}
\newcommand{\U}{\Omega^{\infty}}
\newcommand{\Tot}{\operatorname{\mathsf{Tot}}}
\newcommand{\holim}{\operatorname{holim}}

\newcommand{\SymSeq}{\mathsf{SymSeq}}

\newcommand{\Oper}{\mathsf{Oper}}

\newcommand{\Map}{\operatorname{Map}}
\newcommand{\cro}{\operatorname{cr}}

\newcommand{\hofib}{\operatorname{hofib}}

\newcommand{\colim}{\operatorname{colim}}

\newcommand{\id}{\mathrm{id}}
\newcommand{\Id}{\mathbb{I}}
\newcommand{\Cobar}{\operatorname{Cobar}}

\newcommand{\pt}{*}

\newcommand{\Cal}[1]{\mathcal{#1}}
\newcommand{\Ncirc}[1]{\mathbf{N}^{\hat{\circ}#1}}

\newcommand{\boxcirc}{\mathring{\square}}

\newcommand{\Nl}{\mathbf{N}_{\mathsf{lev}}}

\newcommand{\SD}{\Sigma{\cdot}\mathbf{\Delta}}
\newcommand{\bdel}{\mathbf{\Delta}}

\newcommand{\uS}{\underline{S}}
\newcommand{\Par}{\mathsf{Par}}
\newcommand{\ccirc}{\check{\circ}}

\newcommand{\Prim}{\operatorname{\mathsf{Prim}}}

\newcommand{\uSz}{\uS^{0}}                       
\newcommand{\BC}{\mathsf{B}}                     
\newcommand{\Omw}{\Omega^{\mathsf{w}}}           
\newcommand{\TreeF}{\mathcal{T}}                 
\newcommand{\wt}{\mathsf{w}}                     
\newcommand{\Kop}{\mathcal{K}}                   
\newcommand{\Lop}{\mathcal{L}}                   
\newcommand{\Aop}{\mathcal{A}}                   
\newcommand{\Thet}{\Theta}                       
\newcommand{\coEnd}{\operatorname{coEnd}}
\newcommand{\fat}[1]{\lVert #1 \rVert}           
\newcommand{\Dres}{\Delta^{\mathsf{res}}}
\newcommand{\alive}[2]{#1\langle #2\rangle}      
\newcommand{\basept}{\ast}
\newcommand{\anc}{\operatorname{anc}}            
\newcommand{\Dg}[1]{D(#1)}                       

\title[Derivatives of identity in spaces]{The partition poset complex and the Goodwillie derivatives of the identity in spaces}
\author{Duncan A. Clark}
\email{daclark311@gmail.com}

\begin{document}
\date{\today}
\maketitle

\begin{abstract}
    We produce a canonical highly homotopy-coherent operad structure on the derivatives of the identity functor in spaces via a pairing of cosimplicial objects, providing a new description of an operad structure on such objects first described by Ching. We prove that the two structures agree: both are restrictions of a single algebra over an operad of windowed cut systems on weighted trees, whose component spaces are contractible, along maps of structuring operads which are all levelwise weak equivalences. In addition, we show the derived primitives of a commutative coalgebra in spectra form an algebra over this operad.
\end{abstract}

\section{Introduction}
The Goodwillie derivatives of the identity on the category of based spaces $\Top$ (denoted $\partial_* \Id)$ play a central role in the homotopy theory of based spaces. As the coefficients of the Taylor tower, the sequence controls a natural filtration of a space $X$ which begins with stabilization $QX$ and converges to $X$ whenever $X$ is $1$-connected\footnote{In fact, this tower will converge for nilpotent spaces as well.}. An operad structure on $\partial_* \Id$ was first constructed by Ching in \cite{Ching-thesis}. Arone-Ching \cite{AC-1} have further shown that the Goodwillie derivatives of a reduced homotopy functor $F\colon \Top \to \Top$ naturally form a $\partial_*\Id$-bimodule, which can be used to rebuild the Taylor tower of $F$ \cite{AC-coalg}, \cite{AC-2}.

The author has shown in \cite{Clark-1} that the Goodwillie derivatives of the identity on the category of algebras over a (reduced) operad of spectra $\Cal{O}$ forms a natural ``highly homotopy-coherent'' operad, and moreover, that this operad structure is equivalent to that on $\Cal{O}$. In this document, we show that similar techniques may be utilized to provide an alternate construction for an operad structure on the Goodwillie derivatives of the identity in based spaces, specifically we show the following.

\begin{theorem} \label{thm:main-1}
The symmetric sequence $\partial_*\Id$ is a highly homotopy coherent operad (i.e., $A_{\infty}$-operad in the language of \cite{Clark-1}).
\end{theorem}

Our main interest in operads is in what structure is described on their algebras. Ching further shows that the homology of $\partial_*\Id$ gives a ``Lie operad'' \cite{Ching-thesis}, which is \textit{Koszul dual} to the commutative cooperad \cite{Ginz-Kap}. The operad $\partial_*\Id$ of spectra is similarly Koszul dual to the commutative cooperad in spectra (as in \cite{FG}) and so is often referred to as the \textit{spectral Lie operad}.

In addition, we prove the following theorem (Theorem \ref{thm:main-2}). As the commutative cooperad in $\Spt$ is also an operad, Theorem \ref{thm:main-2} may be thought of as an amplified version of the classic result from commutative algebra which states that the primitives of a Hopf algebra naturally form a Lie algebra (see, e.g. \cite{Abe}).

\begin{theorem} \label{thm:main-2}
    The derived primitives (Definition \ref{def:derived-prim}) of a commutative coalgebra\footnote{For us, coalgebra means \textit{coalgebra with divided powers} (see, e.g. \cite{FG}).} in spectra admit a natural action by $\partial_*\Id$.
\end{theorem}

The above theorem is also not necessarily new. For a commutative coalgebra $Y$, the derived primitives $\Prim(Y)$ is equivalent to $\TQ(Y^{\vee})^{\vee}$ (see Example \ref{ex:SW-TQ}), where here $Y^{\vee}=\Map(Y,S)$ is the Spanier-Whitehead dual, and the latter is known to admit an action by $\partial_* \Id$ (e.g. as in \cite{Ching-thesis}, \cite{BR}, \cite{Heuts-vn}).  However, our constructions provide an easy to use, algebraic alternative to the pioneering operad structure on $\partial_*\Id$ found in \cite{Ching-thesis}. Operad structures on the Goodwillie derivatives of the identity in settings beyond $\Top$ have also been constructed, by $\infty$-categorical methods; see \cite{Ching-day-conv}, \cite{BH-charLie}.

Since the constructions of this document and those of \cite{Ching-thesis} produce operad structures on the same underlying symmetric sequence by quite different means, it is natural to ask whether the two structures agree. Ching's structure is strict and arises by dualizing a cooperad structure on a topological bar construction, while ours is homotopy-coherent and arises from a $\boxcirc$-monoid structure on a cosimplicial symmetric sequence. Our third main theorem answers this affirmatively. We state it informally here; the precise statement is Theorem \ref{thm:main}, which is proved in Appendix \ref{app:comparison} (see \ref{sec:ack} for how this proof is produced).

\begin{theorem} \label{thm:main-3}
    The highly homotopy coherent operad structure on $\partial_* \Id$ of Theorem \ref{thm:main-1} agrees, as a highly homotopy coherent operad, with the operad structure on $\partial_* \Id$ constructed by Ching in \cite{Ching-thesis}.
\end{theorem}

\subsection{Outline of the argument} \label{sec:outline-comparison}
We make use of the models for $\partial_n \Id$ described by Johnson \cite{Joh} and Arone-Mahowald \cite{Ar-Mah} as the Spanier-Whitehead dual of the $n$-th \textit{partition poset complex} $\Par(n)$ (Definition \ref{def:par-poset}), \[ \partial_n \Id \simeq \Map(\Par(n),S).\]

We then show that $\Map(\Par(n),S)$ can be modeled as the totalization of the cobar resolution $C(\uS)[n]$ on the commutative cooperad of spectra. Essentially, $C(\uS)[n]$ has $p$-simplices given by the $p$-fold composition $\uS^{\circ p}[n]$, where $\uS$ is the reduced symmetric sequence with $\uS[n]=S$ (with trivial $\Sigma_n$ action) for $n\geq 1$. We show that the cosimplicial symmetric sequence $C(\uS)$ admits a natural $\boxcirc$-monoid structure (see \cite[Definition 5.6]{Clark-1}) induced, essentially, by the canonical isomorphisms $\uS^{\circ p} \circ \uS^{\circ q} \cong \uS^{\circ p+q}$. This $\boxcirc$-monoid structure then induces the desired ``highly homotopy coherent'' operad structure on $\partial_* \Id$ upon passing to totalization.

We similarly show that the derived primitives may be calculated as the totalization of a certain cosimplicial spectrum $C(Y)$. In particular, $C(Y)$ is precisely the dual notion of derived indecomposables which underlies the construction of \textit{topological Quillen homology} ($\TQ$) \cite{Bast}, \cite{Bast-Man}, \cite{Kuh}, \cite{Har-TQ} of structured ring spectra described as an algebra over an operad in spectra. We show that $C(Y)$ is a $\boxcirc$-module over the cobar resolution $C(\uS)$ and that this module structure induces an action of $\partial_*\Id$ on the derived primitives upon passing to totalization.

\subsection{Comparison with Ching's construction} \label{sec:outline-ching}
Ching's original construction \cite{Ching-thesis} proceeds by exhibiting the partition poset complexes as a topological bar construction $\BC(n)$ on the commutative operad of spaces, equipping the bar construction with a cooperad structure by cutting weighted trees, and dualizing arity-wise. Our main technical lemma (Lemma \ref{lem:tech}) is reminiscent of the ``tree ungrafting'' arguments found in \cite{Ching-thesis}, and our notion of $\boxcirc$-monoid is similar to Ching's notion of ``pre-cooperad'' in \cite{Ching-barcobar}. Despite these similarities, the defining formulas of the two operad structures do not match: Ching's cuts are not tied to any fixed height, whereas the box-product formalism forces cuts at height $\tfrac{1}{2}$ (Remark \ref{rem:naive-frontback-fails}). Neither structure is a formal consequence of the other.

The comparison is carried out in Appendix \ref{app:comparison}, where Theorem \ref{thm:main-3} is stated precisely and proved as Theorem \ref{thm:main}; an outline is given in \S\ref{sec:app-outline}. In brief: an explicit $\Sigma_n$-equivariant homeomorphism $|P(n)| \cong \BC(n)$ (Theorem \ref{thm:levelling}) places both structures on a single object, on which each is given by \textit{cutting} weighted trees along flags of partitions and rescaling the resulting pieces affinely, the two differing only in where the cuts are made. Both prescriptions are then exhibited as cells of a single $\Nl$-operad $\Kop$ of \textit{windowed cut systems}, whose component spaces are contractible, and the two structures are identified as restrictions of a single $\Kop$-algebra along levelwise weak equivalences of structuring operads. 

In the $\infty$-categorical literature, the derivatives of the identity and their operad structure are treated via Koszul duality and Day convolution \cite{Ching-day-conv}, \cite{Heuts-vn}, \cite{BR}; most relevantly, Blans-Heuts \cite{BH-charLie} characterize the spectral Lie operad abstractly and prove that the functor $\partial_*$ induces this operad structure, noting that such a proof had not previously appeared. These results do not subsume Theorem \ref{thm:main-3}: the comparison of Ching's original structure \cite{Ching-thesis} with the $\infty$-categorical ones remains open. Indeed, \cite[Example 3.3]{Ching-day-conv} records that the agreement of the structure of \cite{Ching-thesis} with the Day convolution structure is not obvious, and \cite{BH-charLie} does not address the original construction. The present document proves the point-set statement directly, and no result of it relies on the $\infty$-categorical literature. See also \cite{BB-chainrule} for a highly structured form of the Arone-Ching chain rule \cite{AC-1} (see also \cite{Klein-Rog}, \cite{Yeak}) over $\partial_*\Id$, confirming a conjecture of Lurie.

\subsection{Applications and related results}
Any space $X$ gives rise to a commutative coalgebra $\Sp X$ and so the map \[X\mapsto \Prim C(\Sp X)\] provides a functor from $\Top$ to algebras over $\partial_*\Id$. Work of Heuts \cite{Heuts-vn} shows that this functor induces an equivalence (on spaces which converge to their Goodwillie tower) after localization with respect to $v_n$-periodic homotopy equivalences (see, e.g. \cite{BR}); in the same work the chromatic localizations $M_n^f$ of the $\infty$-category of based spaces are identified with the $\infty$-category of spectral Lie algebras (i.e., algebras over $\partial_*\Id$) in $T(n)$-local spectra.

The key idea is that for $X\in \Top$ there is a coaugmentation (see also \cite{BR-1}) \[ \Phi_n(X)\to C(\Sp_{T(n)}X) \] where $\Phi_n$ denotes the Bousfield-Kuhn functor  (see, e.g. \cite{Bou-BK-functor}, \cite{Kuh-BK-functor}, \cite{Kuh-tel}, \cite{Heuts-vn}) and $C(\Sp_{T(n)} X)$ is the cosimplicial diagram which computes the derived primitives of the commutative coalgebra $\Sp_{T(n)}X =L_{T(n)}\Sp X$ in $T(n)$-local spectra. For instance, when $n=0$ this is precisely the result of Quillen \cite{Quillen} which provides algebraic models for rational homotopy theory as either Lie algebras or commutative coalgebras in chain complexes.

\subsection{Conventions and notation} \label{sec:conventions}
We make use of the pointed, compactly generated Hausdorff spaces as our model for $\Top$, and use the theory of EKMM $S$-modules developed in \cite{EKMM} as our model for a closed symmetric monoidal category of spectra. The main technical benefit for us is that every spectrum will be fibrant. We write $\Spt$ for the category of $S$-modules and refer to such objects as \textit{spectra}. We will use $\Map(-,-)$ to denote the mapping spectrum in $\Spt$.

We write $\Delta^{\mathsf{res}}$ for the subcategory of the usual simplex category $\Delta$ obtained by omitting degeneracy maps. We will write $\Delta^{\bullet}$ for the cosimplicial space of topological simplices. The category $\Spt$ is topologically enriched and we write $\Tot$ for the \textit{restricted totalization} $\mathrm{Tot}^{\mathsf{res}}=\Map_{\Delta^{\mathsf{res}}} (\Delta^{\bullet}, -)$. A useful observation is that since all spectra are fibrant, given any $Y\in \Spt^{\Delta}$ there are equivalences \[
    \holim_{\Delta} Y\xrightarrow{\sim} \holim_{\Delta^{\mathsf{res}}} Y \simeq \Tot Y\]
without separately requiring that $Y$ be Reedy fibrant (see, e.g., \cite[\textsection 8]{Ching-Harper-Koszul}).

We denote by $\SymSeq$ the category of symmetric sequences in spectra. A symmetric sequence $A$ consists of objects $A[n] \in \Spt$ with (right) action by $\Sigma_n$ for each $n\geq 0$. We will often assume that our symmetric sequences are reduced, that is, $A[0]=\pt$. The category $\SymSeq$ admits a monoidal product, the composition product $\circ$ of Definition \ref{comp-prod}, whose unit we denote $I$.

\subsection{Organization of the paper}
Section 2 provides the necessary background on functor calculus and the Goodwillie derivatives of the identity functor on $\Top$. In Section 3 we develop the necessary language of cooperads and their coalgebras. Section 4 is dedicated to proving Theorems \ref{thm:main-1} and \ref{thm:main-2}.

Appendix \ref{app:comparison} states Theorem \ref{thm:main-3} precisely (Theorem \ref{thm:main}) and proves it; it carries its own outline in \S\ref{sec:app-outline}. Appendix \ref{sec:anchor-low} exhibits the comparison explicitly in arity 3.

\subsection{Acknowledgements} \label{sec:ack}
The author would like to thank John Harper, Michael Ching, and Connor Malin for useful discussions and feedback. Additionally, the author wishes to disclose that the proof of Theorem \ref{thm:main-3} (Appendix \ref{app:comparison}) was developed with substantial assistance from Anthropic's Claude models Fable 5/5.1, and validated numerically\footnote{See \url{https://github.com/dclark202/math/tree/main/derivatives-id-spaces}} in low arities. The author was partially supported by National Science Foundation grant DMS-1547357 and Simons Foundation: Collaboration Grants for Mathematicians \#638247.

\section{Functor calculus and the derivatives of the identity}

Functor calculus was introduced by Goodwillie in \cite{Goo-Calc2}, \cite{Goo-Calc3} as a method for analyzing homotopy preserving functors between $\Top$ and $\Spt$. A central construction in functor calculus is that of the \textit{Taylor tower} of $n$-excisive approximations associated to a functor $F\colon \Top \to \Top$ as follows
\begin{equation}\label{Taylor-tower}
    \xymatrix{
    & & D_nF\ar[d] & &&\\
    F\ar[r] & \cdots \ar[r] & P_nF \ar[r] & P_{n-1}F \ar[r] & \cdots \ar[r] & P_0F.}
\end{equation}

The functors $P_nF$ are called the \textit{$n$-th excisive approximation} to $F$ and are initial in the homotopy category of $n$-excisive functors receiving a map from $F$. All of our approximations are based at the zero object $\pt \in \Top$. 

For $n\geq 1$, we set the $n$-th \textit{homogeneous layer}, denoted $D_nF$, to be \[D_nF :=\hofib(P_nF\to P_{n-1}F).\] 
The following proposition summarizes the salient properties of $D_n F$ (see \cite{Goo-Calc3}). 

\begin{proposition}
\label{prop:functor-calc-recap-top}
Let $F\colon \Top\to \Top$ be a homotopy functor and $n\geq 1$. Then: \begin{itemize}
\item $D_nF$ is \textit{$n$-homogeneous}\footnote{Recall a functor $G$ is $n$-homogeneous if $G$ is $n$-excisive and $P_{n-1}G\simeq \pt$.}.

\item $D_nF$ naturally factors through $\Spt$ as $ D_nF\simeq \U \circ \mathbb{D}_nF \circ \Sp$ such that $\mathbb{D}_n F$ is $n$-homogeneous.

\item $D_nF$ is characterized by a spectrum with (right) $\Sigma_n$-action, $\partial_n F$ called the \textit{$n$-th derivative} of $F$, and there is an equivalence\footnote{For $X$ which have the homotopy type of finite CW-complexes; or on an arbitrary space $X$ if $F$ commutes with filtered colimits (i.e., $F$ is \textit{finitary}). } \[
    D_n F(X) \simeq \U (\partial_nF \wedge_{\Sigma_n} (\Sp X)^{\wedge n}).\]
    
\item The spectrum $\partial_nF$ may be calculated via cross effects \cite[\textsection 3]{Goo-Calc3} as \[\partial_n F \simeq \cro_n \mathbb{D}_n F (S,\dots,S)\] with $\Sigma_n$-action given by permuting the inputs $S$.
\end{itemize}
\end{proposition}

\subsection{The Taylor tower of the identity functor} 
The Taylor tower of the identity is a central object of homotopy theory. From the definitions in \cite{Goo-Calc3}, it is not hard to show \begin{equation} \label{P_1Id}
    P_1\Id(X)\simeq D_1\Id(X)\simeq \U \Sp (X)
    \end{equation}
the stabilization of a space $X$. The higher Blakers-Massey theorems \cite[2.5, 2.6]{Goo-Calc2} show that $\Id$ is $1$-analytic and therefore the Taylor tower of the identity in $\Top$ offers an interpolation between a simply connected space $X\simeq \holim_n P_n\Id(X)$ and its stabilization $\U \Sp X$. 

The equivalences from (\ref{P_1Id}) show that the first derivative of $\Id$ is just the sphere spectrum $S$. Johnson \cite{Joh} and later Arone-Mahowald \cite{Ar-Mah} give a description of the higher homogeneous layers and derivatives of $\Id$ in terms of the\textit{ partition poset complex} (Definition \ref{def:par-poset}). Specifically, they show 
\begin{equation} \label{D_n I}
    D_n \Id (X)\simeq \U\Big( \Map \big(\mathsf{Par}(n), \Sp X^{\wedge n}\big)_{h \Sigma_n}\Big)
\end{equation}
and similarly \begin{equation} \label{der_n I}\partial_n \Id\simeq \Map( \mathsf{Par}(n), S)\end{equation}
for all $n\geq 1$. In particular, the $n$-th derivative of $\Id$ is just the Spanier-Whitehead dual to $\Par(n)$.

\subsection{The partition poset complex}
For $n\geq 1$ we denote by $\mathbf{n}$ the set $\{1,\dots,n\}$
. A \textit{partition} $\lambda$ of $\mathbf{n}$ is a decomposition $\mathbf{n}=\coprod_{i\in I} T_i$ into nonempty subsets (here $I$ is required to be a nonempty set). Given partition $\lambda=\{T_i\}_{i\in I}$ and $\lambda'=\{T'_j\}_{j\in J}$ of $\mathbf{n}$ we say that $\lambda \leq \lambda'$ if there is a surjection $f\colon J\to I$ such that $T_i=\coprod_{j\in f^{-1}(i)} T'_j$ for all $i\in I$. 

Note that the set of partitions of $\mathbf{n}$ has a minimal element $\mathsf{min}$ consisting of only the trivial partition $\{1,\dots,n\}$, and a maximal element $\mathsf{max}$ consisting of the partition of $\mathbf{n}$ into singletons, i.e. $\{\{1\},\dots, \{n\}\}$. The set of partitions of $\mathbf{n}$ then forms a poset with respect to $\leq$, and so may be interpreted as a category. The partition poset complex as defined below is (a quotient of) the nerve of this category.

\begin{definition} \label{def:par-poset} Define the \textit{$n$-th partition poset complex} $\Par(n)$ to be the geometric realization of the pointed simplicial set $P(n)$ defined as follows. The $k$-simplices of $P(n)$ are given by sequences \[ \lambda_0 \leq \lambda_1 \leq \cdots \leq \lambda_{k-1} \leq \lambda_k \] of partitions of $\mathbf{n}$ such that any chain that does not satisfy $\lambda_0=\mathsf{min}$ and $\lambda_k=\mathsf{max}$ is identified with the basepoint. 

Face maps $d_i\colon P(n)_k\to P(n)_{k-1}$ are given by removing the $i$-th entry $\lambda_i$ and degeneracy maps $s_j\colon P(n)_k\to P(n)_{k+1}$ are given by repeating the $j$-th entry $\lambda_j$. Note that the image of $d_0$ (resp. $d_k$) is only the basepoint if $\lambda_1\neq \mathsf{min}$ (resp. $\lambda_{k-1}\neq \mathsf{max}$).
\end{definition}

More generally, for a finite set $T$ we define $\Par(T)$ analogously, e.g. by setting $|T|=n$ and specifying a bijection $T\cong \mathbf{n}$. 

\begin{remark}
    Note that $\Par(n)$ inherits a natural action of $\Sigma_n$ by permuting the elements of $\mathbf{n}$. A useful observation is that non-basepoint elements $\alpha \in P(n)_k$ are in bijective correspondence with the collection of rooted trees with $n$ labelled leaves and $k$ levels up to isomorphism.
\end{remark}

\begin{example} \label{ex:par-low}
    In low arities, $\Par(1)\cong S^0$ and $\Par(2)\cong S^1$ with trivial $\Sigma_2$ action, while $\Par(3)$ may be identified with the $2$-sphere with a disc glued in at the equator, $\Sigma_3$ acting by permuting the three $2$-discs (top hemisphere, bottom hemisphere and equator); see also Appendix \ref{sec:anchor-low}.

\end{example}

We briefly describe one route of arriving at the model \eqref{der_n I}, using the approach of Arone-Kankaanrinta \cite{AK} to analyze $\Id$ by the Snaith splitting. 

\subsection{Analysis of the Taylor tower of the identity by higher stabilization}
Associated to the stabilization adjunction $(\Sp, \U)$ between $\Top$ and $\Spt$, for any space $X$, there is a coaugmented cosimplicial diagram $X\to Q^{\bullet+1}X$. Here, $Q^{\bullet+1}X$ is the cobar resolution \[Q^{\bullet+1}X:=\Cobar(\U, \Sp\U, \Sp X)\] and the coaugmentation is provided by the unit map $X\to \U \Sp X$. (We reserve the symbol $C(-)$ for the cosimplicial objects of $\Spt$ introduced in \S\ref{com-coop} and Definition \ref{Q-coalg}.) The resolution $Q^{\bullet+1}X$ is functorial in $X$ and provides a cosimplicial functor \begin{equation}
\label{eq:loop-susp-cosimplicial-res}
    Q^{\bullet+1} = \left( \xymatrix{
    \U \Sp  \ar@<.5 ex>[r] \ar@<-.5ex>[r] &
    (\U \Sp)^2 \ar@<1 ex>[r] \ar[r] \ar@<-1 ex>[r] &
    (\U \Sp)^3 \cdots
    }\right)
\end{equation} whose coface maps are induced by inserting the unit map $\Id \to Q:= \U \Sp$ and codegeneracy maps are induced by the counit map $\Sp\U=: \mathsf{K}\to \Id_{\Spt}$. 

Blomquist-Harper \cite{BH} utilize the cubical higher Blakers-Massey theorems of \cite{Goo-Calc2} to recover a classical result of Bousfield-Kan \cite{BK} (see also \cite{Carl}): a $1$-connected $X\in \Top$ is equivalent to its completion with respect to stabilization, $X\xrightarrow{\sim} \holim_{\Delta}Q^{\bullet+1}X\simeq X_{\U \Sp}^{\wedge}$. Their connectivity estimates give, by the argument of \cite[Proposition 4.7]{Clark-1}, equivalences $P_n \Id \xrightarrow{\sim} \holim_{\Delta^{\leq k-1}} P_n (Q^{\bullet+1})$ for $k\geq n\geq 1$ (see also \cite[\textsection 16]{AC-1}), and likewise on $D_n$ and $\partial_n$. The upshot for us is that $\partial_n (Q^{k+1})$ is readily computable via the Snaith splitting, as follows.

\subsection{The Snaith splitting}
Let $\uS$ denote the symmetric sequence in $\Spt$ such that $\uS[n]=S$ with trivial $\Sigma_n$ action. The \textit{Snaith splitting} (see e.g., \cite{Snaith}, \cite{Coh-May-Tay}) provides equivalences\footnote{Here and throughout, $X$ is implicitly assumed to be cofibrant (e.g., of the homotopy type of a CW complex) so that the displayed orbit constructions have the correct homotopy type.}
\begin{equation} 
\label{Snaith}
    \Sp \U \Sp X \simeq \bigvee_{k\geq 1} \Sp X^{\wedge k}_{h\Sigma_k} \simeq \bigvee_{k\geq 1} S\wedge_{\Sigma_k} (\Sp X)^{\wedge k}
\end{equation}
Said differently, the Taylor tower for $\mathsf{K}=\Sp\U$ splits as a coproduct of its homogeneous layers when evaluated on a suspension spectrum and that $\partial_* \mathsf{K}\simeq \uS$. 

A result of Arone-Kankaanrinta \cite{AK} uses the above splittings to recover the model for $n$-th homogeneous layers and $n$-th derivatives of the identity in spaces in (\ref{D_n I}) and (\ref{der_n I}), respectively. The crux of their argument is that iterating the Snaith splitting provides equivalences \[\partial_n (Q^{k+1})\simeq \partial_n (\mathsf{K}^{k})\simeq \uS ^{\circ k}[n].\] 
Here, $\uS^{\circ k}$ denotes the $k$-fold composition of the symmetric sequence $\uS$.

\begin{remark}A key observation is that $\uS^{\circ k}[n]$ is just a wedge of copies of the sphere spectrum $S$ indexed by the $k$-simplices of $P(n)$. This symmetric sequence $\uS$ admits both an operad and \textit{cooperad} structure, and the derivatives of the identity may be further recognized as the totalization of the cobar complex with respect to this cooperad structure. 
\end{remark}

\section{Cooperads and their coalgebras}

The aim of this section is to describe what we mean by cooperads and their coalgebras. We recall that an operad $\Cal{O}$ is just a monoid for the composition product $\circ$ of symmetric sequences; the dual notion of cooperad is a more delicate matter.

\subsection{Reinterpreting the composition product}
We will utilize the terminology of $\Nl$-objects developed previously in \cite{Clark-1} and encourage the curious reader to find more detail in \cite[Section 6]{Clark-1}. We first recall that given $n,k>0$ that $\Ncirc{k}[n]$ denotes the collection of sequences of orbits of nonzero integers of the form \begin{equation} \label{orbit-pro} \underline{p}=\left( n^1, (n^2_1,\cdots, n^2_{n^1})_{\Sigma_{n^1}}, \cdots, (n^k_1, \cdots, n^{k}_{n^{k-1}})_{\Sigma_{n^{k-1}}} \right)\end{equation}
where for $1<\ell \leq k$, $n^{\ell}$ is inductively defined as \[n^{\ell} := n^{\ell}_1+\cdots+n^{\ell}_{n^{\ell-1}}, \quad n^{\ell}_j > 0\] and $n^k=n$. We often drop the orbit notation in \eqref{orbit-pro} when clear from context. We call an element $\underline{p}\in \Ncirc{k}[n]$ a \textit{profile}, noting that such profiles naturally occur as the indexing set of $k$-fold iterates of the composition product of symmetric sequences. 

Given a non-basepoint element $\alpha\in P(n)_k$, we let $|\alpha|$ denote the corresponding profile in $\Ncirc{k}[n]$ and $\alpha_{j,i}$ be such that \[|\alpha|=\left(\alpha_{1,1}, (\alpha_{2,1},\cdots,\alpha_{2,\alpha_1}), \cdots, (\alpha_{k,1},\cdots,\alpha_{k,\alpha_{k-1}} ) \right).\]
We note that $\alpha_j$ is the number of equivalence classes in the partition $\lambda_j$, and $\alpha_{j,1},\dots, \alpha_{j,\alpha_{j-1}}$ are the sizes of such equivalence classes appearing in $\lambda_{j-1}$ for $j=1,\dots,k$. Note that $\alpha$ is \textit{not} determined by $|\alpha|$, and that the entries $\alpha_{j,i}$ are only defined up to the orbit relations in \eqref{orbit-pro}.

\begin{definition} \label{comp-prod} Let $A_1,\dots,A_k$ be reduced symmetric sequences. We define their \textit{composition product} as follows.
\[(A_1{\circ} A_2 {\circ} \cdots {\circ} A_k)[n]=\bigvee_{\alpha\in P(n)_k} (A_1 \otimes \cdots \otimes A_k)[\alpha] \]
Here, we use the notation \[
    (A_1 \otimes \cdots \otimes A_k)[\alpha]=
        A_1[\alpha_1] \wedge \bigwedge_{i=1}^{\alpha_1} A_2[\alpha_{2,i}] \wedge \cdots \wedge \bigwedge_{i=1}^{\alpha_{k-1}} A_k[\alpha_{k,i}]
        \]  as $\Sigma_{\alpha}$-objects\footnote{We denote by $\Sigma_{\alpha}\leq \Sigma_n$ the subgroup of permutations $\sigma$ which fix the chain $\alpha$.} if $\alpha$ is not the basepoint, and set $(A_1\otimes \cdots \otimes A_k)$ of the basepoint to be the terminal spectrum $*$. The $\Sigma_n$ action on $(A_1\circ \cdots \circ A_k)[n]$ is induced by the permutation action of $\Sigma_n$ on $P(n)_k$ together with the $\Sigma_{\alpha}$-structure of the factors (and likewise for the dual composition product below).

Similarly, their \textit{dual composition product} is defined as \[ (A_1\ccirc A_2 \ccirc \cdots \ccirc A_k)[n]=\prod_{\alpha \in P(n)_k} (A_1 \otimes \cdots \otimes A_k)[\alpha]. \]
\end{definition}

We write $A^{\otimes k}$ for the $k$-fold product $A\otimes \cdots \otimes A$ and note that $A^{\otimes 0}:= I$. Note that this symbol $\otimes$ is different from the tensor of symmetric sequences as in \cite{Rezk-thesis}, \cite{Har-1}. Let $P(n)^{\circ}_k$ denote the set of non-basepoint $k$-simplices of $P(n)$. Since $\Spt$ is \textit{stable}, finite coproducts and products are equivalent and hence the natural comparison \begin{equation} \label{comp} A_1{\circ} \cdots {\circ} A_k\xrightarrow{\;\sim\;} A_1\ccirc \cdots \ccirc A_k\end{equation} is a weak equivalence of symmetric sequences. 

\begin{remark} 
The dual composition product is rarely strictly associative, and therefore we cannot say that a cooperad is a comonoid with respect to $\ccirc$. The issue is that the smash product of spectra will rarely commute with limits; however for $F$ a small diagram of spectra and $X \in \Spt$, the induced maps  $(\lim F ) \wedge X\to \lim ( F\wedge X)$ make $\ccirc$ \textit{oplax monoidal} (see, e.g. \cite{Day-St}, \cite{Ching-oplax}). 

We eschew the full development of (op)lax (co)monoids and only state what we need for a symmetric sequence $\Cal{Q}$ to be a \textit{cooperad}. Note that what we are calling a cooperad is more precisely a \textit{coaugmented cooperad with divided powers}; similarly, our notion of coalgebra is that of a \textit{coalgebra with divided powers} \cite{FG}.
\end{remark} 

\subsection{Cooperads} 
Informally, a cooperad is a reduced symmetric sequence $\Cal{Q}$ that admits \textit{cocomposition} maps of the form \[ \delta \colon \Cal{Q}[k]\to \Cal{Q}[n] \wedge \Cal{Q}[k_1] \wedge \cdots \wedge \Cal{Q}[k_n] \] for all $n,k=k_1+\dots+ k_n\geq 1$, along with a \textit{counit} $\epsilon \colon \Cal{Q}[1]\to S$, which satisfy certain associativity, unitality and equivariance conditions. 

Equivalently, we may write the above cocomposition maps as a collection \[\Cal{Q}[d_1\alpha] \to \Cal{Q}^{\otimes 2}[\alpha] \quad \quad (\alpha\in P(n)_2).\] Further our cooperads are coaugmented in that the counit $\Cal{Q}[1]\to S$ is required to admit a retract $\eta\colon S\to \Cal{Q}[1]$. 

\begin{definition} A symmetric sequence $\Cal{Q}$ which admits comultiplication maps is a \textit{cooperad} if there are well-defined cosimplicial objects in $\Spt^{\Sigma_n}$ \[ C(\Cal{Q})[n] = {\prod}_{\alpha\in P(n)_{\bullet}}\Cal{Q}^{\otimes \bullet} [\alpha]\]
for all $n\geq 1$. Coface and codegeneracy maps are induced by the face and degeneracy maps from $P(n)_{\bullet}$ as follows. The coface maps $d^i\colon C(\Cal{Q})[n]^{k-1}\to C(\Cal{Q})[n]^{k}$ are induced by $\Cal{Q}^{\otimes k-1}[d_i \alpha]\to \Cal{Q}^{\otimes k}[\alpha]$; that is, by inserting the comultiplication maps from $\Cal{Q}$ at the $i$-th position for $0< i<k$, and by applying the coaugmentation $S\to \Cal{Q}[1]$ in the first (resp. last) slots for $i=0$ (resp. $i=k$)--see Example \ref{coop-ex}. 

Similarly, for $0\leq j \leq k$, \[s^j\colon C(\Cal{Q})[n]^{k+1}\to C(\Cal{Q})[n]^{k}\] is induced by the analogous codegeneracy map $s^j \colon \Cal{Q}^{\otimes k+1} [s_j \alpha] \to \Cal{Q}^{\otimes k} [\alpha].$
\end{definition}

\begin{example} \label{coop-ex}
Recall that $\Cal{Q}^{\otimes 0}=I$, so $C(\Cal{Q})[n]$ begins \begin{equation}
        C(\Cal{Q})[n]\cong \left( \xymatrix{
            I[n] \ar@<.5ex>[r] \ar@<-.5ex>[r] &
            \Cal{Q}[n] \ar@<1ex>[r] \ar[r] \ar@<-1ex>[r] &
            \prod_{\alpha\in P(n)_2^{\circ}}  \Cal{Q}[s]\wedge \Cal{Q}[t_1]\wedge\cdots\wedge \Cal{Q}[t_{s}] \; \cdots
        }\right)
    \end{equation}
    where we write $(s,(t_1,\dots,t_s))=|\alpha|$. The two maps $I[n]\to \Cal{Q}[n]$ are induced by the coaugmentation $\eta\colon S\to \Cal{Q}[1]$ if $n=1$ and are the initial maps $\pt\to \Cal{Q}[n]$ if $n\geq 2$; $d^i\colon C(\Cal{Q})[n]^1\to C(\Cal{Q})[n]^2$ is induced by the maps \begin{align*}
        d^0 \colon & \Cal{Q}[n]\cong S\wedge \Cal{Q}[n] \xrightarrow{\eta\wedge \id} \Cal{Q}[1]\wedge \Cal{Q}[n]
        \\
        d^1 \colon & \Cal{Q}[n] \xrightarrow{\delta} \prod_{\alpha\in P(n)^{\circ}_2} \Cal{Q}[s] \wedge \Cal{Q}[t_1]\wedge\cdots \wedge \Cal{Q}[t_s]
        \\
        d^2 \colon & \Cal{Q}[n]\cong \Cal{Q}[n] \wedge S^{\wedge n} \xrightarrow{\id\wedge \eta^{\wedge n}} \Cal{Q}[n] \wedge \Cal{Q}[1]^{\wedge n}
    \end{align*}
    and the trivial map $\Cal{Q}[n]\to \pt$ on all other factors. Here $\delta$ collects the comultiplication maps of $\Cal{Q}$. For any $n\geq 1$ the codegeneracies $s^j$ are induced by the counit $\Cal{Q}[1]\to S$ and $\Cal{Q}[k]\to \pt$ ($k\geq 2$).
\end{example}

\subsection{The commutative cooperad of spectra} \label{com-coop}
The symmetric sequence $\uS$ introduced before admits a natural cooperad structure with comultiplication $\delta$ induced by the natural isomorphisms $S\to S\wedge S\wedge \cdots \wedge S$.

In particular, since $\uS^{\otimes k}[\alpha]=S$ for any non-basepoint $\alpha\in P(n)_k$, it follows that 
\[ C(\uS)[n]^k= {\prod}_{P(n)^{\circ}_{k}} S\]
and the coface (resp. codegeneracy) maps are just induced by the face (resp. degeneracy) maps of $P(n)$. In what follows we will often abbreviate $\prod_{P(n)^{\circ}_k}$ by $\prod_{P(n)_k}$, the factor indexed by the basepoint being terminal.

\begin{proposition} \label{prop:model-equiv}
    There is an equivalence $\partial_* \Id\simeq \Tot C(\uS)$.
\end{proposition}

\begin{proof}
    Using the model from \eqref{der_n I}, this follows from the equivalences \begin{align*}
        \partial_n \Id \simeq \Map \left(|P(n)_{\bullet}|, S\right) 
        \simeq \Tot {\prod}_{P(n)_{\bullet}} \Map(S^0, S) 
         \simeq \Tot C(\uS)[n]
    \end{align*} 
    for all $n\geq 1$.
\end{proof}

\subsection{Coalgebras over a cooperad}
Let $\Cal{Q}$ be a cooperad. Informally, $\Cal{Q}$-coalgebra structure on a spectrum $Y$ consists of comultiplication maps \[ Y\to (\Cal{Q}[n] \wedge Y^{\wedge n} )_{\Sigma_n} \quad \quad (n\geq 1)\] which are required to further satisfy associativity, unitality and equivariance conditions. Note these comultiplications induce a map \[ Y\mapsto \Cal{Q}\ccirc (Y)= \prod_{n\geq 1}\left( \mathcal{Q}[n]\wedge Y^{\wedge n}\right)_{\Sigma_n}\] which we essentially require to be ``oplax comonoidal'' \cite{Day-St}, and that $\Cal{Q}$-coalgebras are the coalgebras for this comonoid. 

First we require some technical observations. Let $\widetilde{P}(n)_k$ denote the quotient set of $P(n)_k$ by the $\Sigma_n$ action on $\mathbf{n}$. Note that an element $\widetilde{\alpha} \in \widetilde{P}(n)_k$ is precisely the collection of all $\alpha\in P(n)_k$ with the same profile $|\alpha|\in \Ncirc{k}[n]$. For any $n,k$, there is then a natural comparison \begin{equation} \label{Pnk-orbits}
    \left( \prod_{\alpha\in P(n)_k} \Cal{Q}^{\otimes k}[\alpha]\wedge Y^{\wedge n} \right)_{\Sigma_n}
    \simeq
    \prod_{ \widetilde{\alpha} \in \widetilde{P}(n)_k } \Cal{Q}^{\otimes k}[\alpha] \wedge_{\Sigma_{\alpha}} Y^{ \wedge n}
\end{equation}
given by the zig-zag through the corresponding wedge \[\Big(\bigvee_{\alpha\in P(n)_k} \Cal{Q}^{\otimes k}[\alpha]\wedge Y^{\wedge n}\Big)_{\Sigma_n} \cong \bigvee_{\widetilde{\alpha}\in \widetilde{P}(n)_k} \Cal{Q}^{\otimes k}[\alpha]\wedge_{\Sigma_{\alpha}} Y^{\wedge n};\] here orbits commute with wedges, though \textit{not} with products, so that the two sides of \eqref{Pnk-orbits} are not literally isomorphic; and for each fixed $n,k$ the wedge-to-product comparisons are stable equivalences (for suitably cofibrant $Y$) since $P(n)_k$ is finite. In the sequel we take the right-hand side of \eqref{Pnk-orbits} as our working model; in particular the coface maps below are constructed on this model.
On the right hand side of \eqref{Pnk-orbits}, $\alpha$ denotes a chosen preimage of $\widetilde{\alpha}$ under the quotient map $P(n)_k\to \widetilde{P}(n)_k$; note that for any such choice of preimage the groups $\Sigma_{\alpha}$ are isomorphic, and $\Cal{Q}^{\otimes k}[\alpha]$ does not depend on this choice.

\begin{definition} \label{Q-coalg}
    Let $\Cal{Q}$ be a cooperad. A $\Cal{Q}$-coalgebra is a spectrum $Y$ that admits a well-defined cosimplicial object $C(Y)$ as follows. For $k\geq 0$, \[ C(Y)^k=\prod_{n\geq 1} \left( \prod_{\alpha \in P(n)_k} \Cal{Q}^{\otimes k}[\alpha] \wedge Y^{\wedge n}\right)_{\Sigma_n }.\]
    (As above, we work with the model afforded by the right-hand side of \eqref{Pnk-orbits}.)

Let $k\geq 0$, coface maps $C(Y)^{k}\to C(Y)^{k+1}$ are induced by the face maps from $P(n)$ for each $n\geq 1$ along with the diagonal maps on $Y$ as follows. \begin{itemize}
    \item $d^0$ is induced by the trivial map $Y^{\wedge n}_{\Sigma_n} \to \pt$ at factors $\alpha\in P(n)_{k+1}$ such that $d_0\alpha\in P(n)_{k}$ is the basepoint, and otherwise by the coaugmentation maps $\Cal{Q}^{\otimes k}[d_0\alpha]\to \Cal{Q}^{\otimes k+1}[\alpha]$ of $C(\Cal{Q})$ smashed with the identity on $Y^{\wedge n}$.
    \item For $i=1,\dots,k$, $d^i$ is induced by $\Cal{Q}$-cooperad structure maps \[ \Cal{Q}^{\otimes k}[d_i\alpha]\to \Cal{Q}^{\otimes k+1}[\alpha]\] for $\alpha\in P(n)_{k+1}$.
    \item $d^{k+1}$ is induced as follows. For $m\geq n$, if $\alpha\in P(m)_{k+1}$ is given by the chain $\lambda_0\leq \cdots \leq \lambda_{k+1}$ with $k$-th entry $\lambda_k=\{T_1,\dots,T_n\}$, then let $\alpha'\in P(n)_k$ be the chain $\lambda_0\leq \cdots \leq \lambda_k$ regarded as a chain of partitions of $\mathbf{n}$; i.e., the result of quotienting the set $\{1,\dots,m\}$ by the relation $a\sim b$ if $a,b\in T_i$ (so that $\lambda_k$ becomes $\mathsf{max}$). Set $t_j=|T_j|$. There is an induced map \begin{align*} 
    \Cal{Q}^{\otimes k}[\alpha']\wedge Y^{\wedge n} 
    &
    \to \Cal{Q}^{\otimes k}[\alpha']\wedge (\Cal{Q}[t_1]\wedge_{\Sigma_{t_1}} Y^{\wedge t_1}) \wedge \cdots \wedge (\Cal{Q}[t_n]\wedge_{\Sigma_{t_n}} Y^{\wedge t_n}) 
    \\
    & \cong \Cal{Q}^{\otimes k}[\alpha']\wedge (\Cal{Q}[t_1]\wedge\cdots \wedge \Cal{Q}[t_n]) \wedge_{\Sigma_{t_1}\times\cdots\times \Sigma_{t_n}} Y^{\wedge t_1+\cdots + t_n} \\
    & 
    \cong \Cal{Q}^{\otimes k+1}[\alpha]\wedge_{\Sigma_{t_1}\times\cdots\times\Sigma_{t_n}} Y^{\wedge m}
    \end{align*} for any $\alpha\in P(m)_{k+1}$ and $m\geq n$. The induced map on $\Sigma_{\alpha'}$ orbits (which acts on $Y^{\wedge m}$ by block permutation on the factors $Y^{\wedge t_i}$) produces a map of the form \begin{equation} \label{alpha'alpha}  \Cal{Q}^{\otimes k} [\alpha']\wedge_{\Sigma_{\alpha'}} Y^{\wedge n} \to (\Cal{Q}^{\otimes k+1}[\alpha]\wedge_{\Sigma_{t_1}\times\cdots\times\Sigma_{t_n}} Y^{\wedge m})_{\Sigma_{\alpha'}}\cong\Cal{Q}^{\otimes k+1}[\alpha] \wedge_{\Sigma_{\alpha}} Y^{\wedge m}.\end{equation} 
    Together with the identifications from \eqref{Pnk-orbits}, $d^{k+1}$  is obtained as the map\begin{equation*}
        \prod_{n\geq 1} \,\prod_{\widetilde{\beta} \in \widetilde{P}(n)_k } \Cal{Q}^{\otimes k}[\beta]\wedge_{\Sigma_{\beta}} Y^{\wedge n} 
        \to
        \prod_{m\geq 1} \,\prod_{\widetilde{\alpha} \in \widetilde{P}(m)_{k+1} } \Cal{Q}^{\otimes k+1}[\alpha]\wedge_{\Sigma_{\alpha}} Y^{\wedge m}
    \end{equation*} induced at a particular $\widetilde{\alpha}\in \widetilde{P}(m)_{k+1}$ by \eqref{alpha'alpha} after first projecting onto the factor $\widetilde{\beta}=\widetilde{\alpha'}$. 
\end{itemize}
Codegeneracy maps are induced by $\Cal{Q}^{\otimes k+1}[s_j\alpha]\to \Cal{Q}^{\otimes k}[\alpha]$ for $\alpha\in P(n)_k$ as in $C(\Cal{Q})$.
\end{definition}

\begin{remark}
    Note that $C(Y)$ is essentially the cobar resolution on $Y$ with respect to the comultiplication map $Y\to \Cal{Q}\ccirc (Y)$.
\end{remark}
\subsection{Derived primitives} 
Let $Y$ be a $\Cal{Q}$-coalgebra. The \textit{primitives} of $Y$, denoted $\bar{Y}$, is given by the (coreflexive) equalizer in $\Spt$ \begin{equation} \label{eq:primitives}
    \bar{Y} = \lim \left( \xymatrix{ Y\ar@<.5ex>[r] \ar@<-.5ex>[r] & \prod_{n\geq 1} (\Cal{Q}[n]\wedge Y^{\wedge n})_{\Sigma_n}} \right). 
\end{equation}
Here, the top map is induced by the coaugmented structure \[Y\cong S\wedge Y\xrightarrow{\eta\wedge \id} \Cal{Q}[1]\wedge Y\]
on the factor $n=1$ and by the trivial map on the factors $n\geq 2$, and the bottom map is induced by the comultiplication maps $Y\to (\Cal{Q}[n]\wedge Y^{\wedge n})_{\Sigma_n}$ for $n\geq 1$. There is a common retract of both maps given by applying $\epsilon\wedge \id\colon \Cal{Q}[1]\wedge Y\to S\wedge Y\cong Y$ and $\Cal{Q}[n]\to \pt$ for $n\geq 2$. 

This gives the precise analog of primitives of a coalgebra from commutative algebra, but fails to be homotopy invariant in general. The derived primitives as defined below gives the homotopy-theoretic analog to this construction, and arises as the dual notion of \textit{topological Quillen homology} ($\TQ$) for algebras over operads in spectra. 

\begin{definition} \label{def:derived-prim} For $Y$ a $\Cal{Q}$-coalgebra we define the \textit{derived primitives} $\Prim Y$ to be the totalization $\Tot C(Y)$. 
\end{definition}

Note that $C(Y)|_{\Delta^{\leq 1}}$ is precisely the equalizer diagram defining $\bar{Y}$ from \eqref{eq:primitives}.

\begin{example} \label{ex:SW-TQ}
    Any $X\in \Top$ gives rise to an $\uS$-coalgebra $\Sp X$ whose comultiplication maps are induced by the diagonals on $X$, i.e., \[\Sp X\to \Sp (X\wedge \cdots \wedge X)\cong \Sp X\wedge \cdots \wedge \Sp X.\]
    Moreover, for any $n\geq 1$ these diagonals naturally commute with the $\Sigma_n$ action on $(\Sp X)^{\wedge n}$ and thus descend to a map \[\kappa_n\colon \Sp X\to (\Sp X)^{\wedge n}_{\Sigma_n} \cong (\uS[n]\wedge (\Sp X)^{\wedge n} )_{\Sigma_n} \]
    coassociative in the sense prescribed by Definition \ref{Q-coalg}.

    In particular, for $Y = \Sp X$ the dual $Y^{\vee}=\Map(Y,S)$ is a commutative ring spectrum, i.e. an algebra over the commutative operad $\mathsf{Com}$, by dualizing the diagonal of $X$; the dual of a general $\uS$-coalgebra is a commutative algebra with divided powers, its structure maps being defined on the fixed points $\big((Y^{\vee})^{\wedge n}\big)^{\Sigma_n}$. If $Y$ is finite (specifically, we need $Y$ to be equivalent to its double Spanier-Whitehead dual), there is further an equivalence  $\TQ(Y^{\vee})^{\vee}\simeq \Prim (Y)$ (see, e.g., \cite{BR}, \cite{Heuts-vn}).
    \end{example}

\section{Operad structure on the derivatives of the identity}
As in \cite{Clark-1}, we make use of the box product $\square$ of cosimplicial objects as introduced by Batanin \cite{Bat-box}.

\subsection{The box product} \label{sec:box-prod} Given a monoidal category $(\mathsf{C},\otimes, \mathbf{1})$ and cosimplicial objects $X,Y\in \mathsf{C}^{\Delta}$, we define the cosimplicial object $X\square Y$ as follows   \begin{equation} \label{box-prod}
    (X\square Y)^{n} := \colim \left(\xymatrix{ \displaystyle{\coprod_{r+s=n-1} X^r \otimes Y^s} \ar@<.75ex>[r] \ar@<-.5ex>[r] & \displaystyle{\coprod_{p+q=n} X^p \otimes Y^q}}\right)
\end{equation} 
The maps in (\ref{box-prod}) are induced by $\id \otimes d^0$ and $d^{r+1}\otimes \id$, see e.g. \cite[Definition 4.13]{Ching-Harper-Koszul} for the cosimplicial structure maps of $X\square Y$.  

Given $c\in \mathsf{C}$ we let $\underline{c}$ denote the constant cosimplicial object on $c$.

\begin{remark} A useful fact is that if $\mathsf{C}$ is closed symmetric monoidal then $(\mathsf{C}^{\Delta}, \square, \underline{\mathbf{1}})$ is a monoidal category. The box product has been used before to construct $A_{\infty}$-pairings on totalizations of $\square$-monoids in $\mathsf{C}^{\Delta}$: McClure-Smith \cite{MS-box} show that if $\mathcal{X}$ is a $\square$-monoid in $\Top^{\Delta}$, then $\Tot \mathcal{X}$ is an $A_{\infty}$-monoid (see also \cite{AC-2}, \cite{Ching-Harper-Koszul}).\end{remark}

We write $\boxcirc$ for the box product with respect to the composition product of symmetric sequences. In our previous work \cite{Clark-1}, we produce an $A_{\infty}$ operad structure on the derivatives of the identity in the category of algebras over a reduced operad of spectra by an underlying pairing of certain cosimplicial symmetric sequences with respect to $\boxcirc$. This is precisely what we will implement in the following section. 

Note that $\boxcirc$ is \textit{not} a strictly monoidal product for $\SymSeq^{\Delta}$. In \cite{Clark-1} we introduced the notion of a \textit{$\boxcirc$-monoid} (\cite[Definition 5.6]{Clark-1}) which is sufficient for a cosimplicial sequence $\Cal{X}$ to admit an $A_{\infty}$-operad structure on its totalization $\Tot \Cal{X}$.

\subsection{Operad structure on the derivatives of the identity in spaces}

We are now prepared to show that $C(\uS)$ admits $\boxcirc$-monoidal structure. We need a technical lemma first; note the following is similar to the ``tree ungrafting'' argument from \cite{Ching-thesis}. 

\begin{lemma} \label{lem:tech}
    For $k,p,q\geq 0$, there is a $\Sigma_k$-equivariant ``decomposition'' map \[\Psi_{k,(p,q)} \colon P(k)_{p+q}\to \prod_{\alpha\in P(k)_2 } \left( P(n)_p\times  P(k_1)_q\times \cdots \times  P(k_n)_q \right)\]
    where $n$ and $k_1,\dots,k_n$ are obtained by setting $| \alpha |=(n,(k_1,\dots,k_n))$.
\end{lemma}

\begin{proof}
    Let $n\geq 1$ and $T_1,\dots,T_n$ be a partition of $\mathbf{k}$. Let $\beta_j\in P(T_j)_q$ be given by $\mu_0^j\leq \cdots \leq \mu^j_q$ for $j=1,\dots,n$, $\gamma\in P(n)_p$ given by $\lambda_0 \leq \cdots \leq \lambda_p$, and let $\lambda'_j$ denote the partition of $\mathbf{k}$ obtained by replacing each block $B\in \lambda_j$ by $\coprod_{s\in B} T_s$. There is then an element $\gamma\circ (\beta_1,\cdots,\beta_n) \in P(k)_{p+q}$ given by \[ \lambda_0'\leq \cdots \leq \lambda_{p-1}' \leq \lambda_p' \cong \coprod_{i=1}^{n} \mu^i_0 \leq \coprod_{i=1}^n \mu^i_1 \leq \cdots \leq \coprod_{i=1}^n \mu^i_q.\]
    
    Given $\alpha\in P(k)_2$, let $T_1,\cdots,T_n$ be the corresponding partition of $\mathbf{k}$ determined by $\alpha$. Given $\gamma'\in P(k)_{p+q}$, the image $\Psi_{k,(p,q)}(\gamma')$ at the $\alpha$-factor of the product is defined to be the string $(\gamma,\beta_1,\cdots,\beta_n)$ if there is a decomposition $\gamma'=\gamma\circ (\beta_1,\cdots,\beta_n)$ where $\gamma\in P(n)_p$ and $\beta_j\in P(T_j)_q$ for $j=1,\dots,n$, and by the basepoint otherwise. 

    We note for later use that a nonbasepoint $\gamma'\in P(k)_{p+q}$, written as $\lambda_0\leq \cdots \leq \lambda_{p+q}$, admits such a decomposition with respect to $\alpha$ if and only if $\lambda_p=\{T_1,\dots,T_n\}$; in this case the decomposition is unique, with $\gamma$ obtained from $\lambda_0\leq \cdots \leq \lambda_p$ by collapsing each $T_i$ to a point and $\beta_j$ obtained by restricting $\lambda_p\leq \cdots \leq \lambda_{p+q}$ to $T_j$. In particular, each nonbasepoint $\gamma'\in P(k)_{p+q}$ decomposes with respect to exactly one $\alpha\in P(k)_2$, namely $\alpha=\lambda_p$.

\end{proof}

\begin{proposition} \label{prop:m_p,q}
    For $p,q\geq 0$ there are maps $m_{p,q} \colon C(\uS)^p \circ C(\uS)^q \longrightarrow C(\uS)^{p+q}$.
\end{proposition}

\begin{proof}
    Let $p,q\geq 0$. For $k\geq 1$, $m_{p,q}$ at level $k$ is given by the following composite
     \begin{align*}
        ( C(\uS)^p \circ C(\uS)^q)[k] & 
        = \bigvee_{\alpha\in P(k)_2} \left( \left( \prod_{P(n)_p} S\right) \wedge \bigwedge_{i=1}^{n} \left(\prod_{P(k_i)_q} S\right) \right)\\ 
        & \to \bigvee_{\alpha \in P(k)_2} \left( \prod_{P(n)_p} S  \wedge \prod_{P(k_1)_q \times \cdots \times P(k_n)_q}  S\wedge\cdots \wedge S \right)\\
        & \to \prod_{\alpha\in P(k)_2}\, \prod_{P(n)_p \times P(k_1)_q \times \cdots \times P(k_n)_q} S\wedge S\wedge \cdots \wedge S\\
        & \xrightarrow{\Psi_{k,(p,q)}^*} \prod_{P(k)_{p+q}} S\wedge S \wedge \cdots \wedge S\cong C(\uS)^{p+q}[k]
    \end{align*}
    where $n,k_1,\dots,k_n$ are such that $ |\alpha|=(n,(k_1,\dots,k_n))$. Here, the map $\Psi_{k,(p,q)}^*$ is well-defined by the uniqueness statement in Lemma \ref{lem:tech}: the coordinate of the image at $\gamma'\in P(k)_{p+q}$ is obtained by first projecting onto the factor indexed by the unique $\alpha\in P(k)_2$ with respect to which $\gamma'$ decomposes, followed by projection onto the coordinate indexed by the resulting decomposition $(\gamma,\beta_1,\dots,\beta_n)$. 
\end{proof}

\begin{proposition}
    \label{prop:box-monoid}
    The cosimplicial symmetric sequence $C(\uS)$ admits a natural $\boxcirc$-monoid structure, i.e., there are maps $m\colon C(\uS)\boxcirc C(\uS)\to C(\uS)$ and $u \colon \underline{I}\to C(\uS)$ which satisfy associativity and unitality. 
\end{proposition}

\begin{proof}
    For $p,q\geq 0$ the maps described in Proposition \ref{prop:m_p,q} provide a square \begin{equation} \label{com-square} \xymatrix{ 
        C(\uS)^p \circ C(\uS)^{q+1} \ar[r]^-{m_{p,q+1}}
        &
        C(\uS)^{p+q+1}
        \\
        C(\uS)^p\circ C(\uS)^{q} \ar[u]^-{\id \circ d^0} \ar[r]^-{d^{p+1} \circ \id}
        &
        C(\uS)^{p+1} \circ C(\uS)^{q} \ar[u]_-{m_{p+1,q}}
    }\end{equation}
    commutativity of which follows from the commutativity of the associated squares \[ \xymatrix{
    \prod_{P(k)_2} (P(n)_{p}\times \prod_{i=1}^n P(k_i)_{q+1} )\ar[d]^-{\id \times \prod_{i=1}^n d_0}
    &&
    P(k)_{p+q+1} \ar[ll]_-{\Psi_{k,(p,q+1)}} \ar[d]^-{\Psi_{k,(p+1,q)}}
    \\
    \prod_{P(k)_2} \left(P(n)_{p} \times \prod_{i=1}^n P(k_i)_{q} \right)
    & &
    \prod_{P(k)_2} (P(n)_{p+1}\times \prod_{i=1}^n P(k_i)_{q}) \ar[ll]_-{d_{p+1} \times \prod_{i=1}^n \id}
    }\] where the products range over all $\alpha\in P(k)_2$, we write $(n,(k_1,\dots,k_n))$ for a given $|\alpha|$, and $\Psi_{k,-}$ is as in Lemma \ref{lem:tech}. The multiplication $m$ is induced at level $\ell$ by considering the diagram built from squares of the form \eqref{com-square} with $p+q+1=\ell$. The same squares of pointed sets, with an arbitrary face $d_i$ or degeneracy $s_j$ in place of $d_0$ and $d_{p+1}$, show that $m$ commutes with the cosimplicial structure maps, so that $m$ is a map of cosimplicial symmetric sequences. 

    The unit map is induced by the coaugmentation $I\to C(\uS)$ given by the identity on $I$. Associativity and unitality of $m$ and $u$ follow from the same argument as in \cite[Proposition 5.8]{Clark-1} with $C(\Cal{O})$\footnote{Though, perhaps to reflect the notation in this article, the cosimplicial object $C(\Cal{O})$ in \cite{Clark-1} should be written as $C(B(\Cal{O}))$, as the bar construction $B(\Cal{O})$ is a cooperad \cite{Ching-thesis}, \textit{not} $\Cal{O}$.} replaced by $C(\uS)$.
\end{proof}

\begin{proof}[Proof of Theorem \ref{thm:main-1}]
    From Proposition \ref{prop:model-equiv}, we know $\partial_* \Id \simeq \Tot C(\uS)$. The theorem then follows from \cite[Proposition 5.9]{Clark-1}.
\end{proof}

\subsection{Derived primitives are spectral Lie algebras}
We now show explicitly that the derived primitives of a $\uS$-coalgebra $Y$ admit a left action by $\partial_*\Id$. In particular, this defines the functor discussed in the introduction, \[ \Xi \colon \Top \to \mathsf{Alg}_{\partial_*\Id}, \quad \quad X\mapsto \Tot C(\Sp X)\]

\begin{proof}[Proof of Theorem \ref{thm:main-2}]The idea is to show that $C(Y)$ is a left $\boxcirc$-module over $C(\uS)$ concentrated in symmetric sequence level $0$. That is, the following associativity \eqref{eq:assoc} and unitality \eqref{eq:unit} diagrams commute \begin{equation} \label{eq:assoc}
    \xymatrix{
        C(\uS)\boxcirc C(\uS)\boxcirc C(Y) \ar[r]^{\theta} \ar[d]^-{\cong}
        &
        C(\uS)\boxcirc (C(\uS)\boxcirc C(Y)) \ar[r]^-{\id\boxcirc \mu}
        &
        C(\uS) \boxcirc C(Y) \ar[d]^-{\mu}
        \\
        (C(\uS)\boxcirc C(\uS))\boxcirc C(Y) \ar[r]^{m\boxcirc \id}
        &
        C(\uS)\boxcirc C(Y) \ar[r]^-{\mu}
        & 
        C(Y)
    }
    \end{equation}
    and 
    \begin{equation} \label{eq:unit} \xymatrix{
        C(\uS) \boxcirc C(Y) \ar[r]^{\mu}
        &
        C(Y)
        \\
        \underline{I}\boxcirc C(Y) \ar[u]^-{u \boxcirc \id} \ar[ur]^-{\cong}
    }
    \end{equation}
Once we have shown this structure, it follows that $\Tot C(Y)$ is a $\partial_*\Id$-algebra in the sense of \cite[Definition 6.36]{Clark-1}. We first produce maps $\mu_{p,q}\colon C(\uS)^p \circ C(Y)^q \to C(Y)^{p+q}$. 

Let $p,q\geq 0$ be given, $\mu_{p,q}$ is then the composite
\begingroup\allowdisplaybreaks
\begin{align*}
         C(&\uS)^p \circ C(Y)^q  = \bigvee_{n\geq 1} C(\uS)^p[n] \wedge_{\Sigma_n} (C(Y)^q)^{\wedge n} \\
         & \to \prod_{n\geq 1} \left( \prod_{P(n)_p} S \right) \wedge_{\Sigma_n} \left( \bigwedge_{i=1}^n \prod_{k_i\geq 1}\left(  \prod_{P(k_i)_q} Y^{\wedge k_i} \right)_{\Sigma_{k_i}} \, \right)
         \\
         & \to \prod_{n\geq 1} \left( \prod_{P(n)_p} S \right) \wedge_{\Sigma_n} \left( \prod_{k\geq 1} \prod_{k=k_1+\dots+k_n} \bigwedge_{i=1}^n \left( \prod_{P(k_i)_q} Y^{\wedge k_i} \right)_{\Sigma_{k_i}} \, \right)
         \\
         & \to \prod_{k\geq 1} \prod_{n\geq 1}  \left(\prod_{P(n)_p} S\right) \wedge_{\Sigma_n} \left( \prod_{k=k_1+\cdots+k_n} \left(\prod_{P(k_1)_q\times \cdots \times P(k_n)_q} \bigwedge_{i=1}^n Y^{\wedge k_i}\right)_{\Sigma_{k_1}\times \cdots \times \Sigma_{k_n}}\right)
         \\
         & \to \prod_{k\geq 1}  \left(  \prod_{\alpha\in P(k)_2} \prod_{P(n)_p\times P(k_1)_q\times \cdots \times P(k_n)_q}  Y^{\wedge k}\right)_{\Sigma_k} \\
         & \xrightarrow{\Psi_{k,(p,q)}^*} \prod_{k \geq 1} \left( \prod_{P(k)_{p+q}} Y^{\wedge k} \right)_{\Sigma_k}
    \end{align*}\endgroup
    Here, the fourth map is induced on orbits by the evident assembly of indices, which is equivariant with respect to the subgroup $\Sigma_n\ltimes(\Sigma_{k_1}\times \cdots \times \Sigma_{k_n})\leq \Sigma_k$ of block permutations, and $\Psi^*_{k,(p,q)}$ is induced on $\Sigma_k$-orbits by the $\Sigma_k$-equivariant projections determined by Lemma \ref{lem:tech}, as in the proof of Proposition \ref{prop:m_p,q}. (We have displayed $C(Y)^q$ in the form of the left-hand side of \eqref{Pnk-orbits}; on the working model, the right-hand side, the same maps are defined verbatim, using $\Sigma_{\gamma'}\cong \Sigma_{\gamma}\ltimes\prod_i \Sigma_{\beta_i}$ for the decomposition $\gamma' = \gamma\circ(\beta_i)$ exactly as in \eqref{alpha'alpha}, and the fourth map is the composite of the isomorphism $(X_{\alpha_0})_{\Sigma_n\ltimes\prod_i\Sigma_{k_i}}\cong (\bigvee_{\alpha} X_{\alpha})_{\Sigma_k}$ with the wedge-to-product comparison.)

    The map $\mu \colon C(\uS)\boxcirc C(Y)\to C(Y)$ is induced at cosimplicial level $0$ by \[ C(\uS)^0\circ C(Y)^0\cong S\wedge Y\xrightarrow{\cong} Y\] and for $n\geq 0$ at level $n+1$ by the commuting squares \[\xymatrix{
        C(\uS)^{p}\circ C(Y)^{q+1} \ar[r]^-{\mu_{p,q+1}} 
        &
        C(Y)^{p+q+1}
        \\
        C(\uS)^{p} \circ C(Y)^{q} \ar[u]^-{\id \circ d^0} \ar[r]^-{d^{p+1} \circ \id}
        &
        C(\uS)^{p+1}\circ C(Y)^{q} \ar[u]_-{\mu_{p+1,q}}
    }\]
    Here $p,q\geq 0$ are such that $p+q=n$. For each $k\geq 1$, commutativity of these squares follows from the commutativity of the same squares of pointed sets appearing in the proof of Proposition \ref{prop:box-monoid}: both composites are trivial on the factor indexed by a chain $\gamma'=(\lambda_0\leq \cdots \leq \lambda_{p+q+1})\in P(k)_{p+q+1}$ unless $\lambda_p=\lambda_{p+1}$, in which case both are induced by projection onto the factor indexed by the common decomposition of $\gamma'$ from Lemma \ref{lem:tech}.

    We emphasize that the comultiplication maps on $Y$ do not participate in the above descent data; rather, they enter through the compatibility of the maps $\mu_{p,q}$ with the coface maps of $C(Y)$. In particular, compatibility with the final cofaces amounts to the identities \[ d^{p+q+1} \mu_{p,q}=\mu_{p,q+1} (\id \circ d^{q+1})\] which hold since the comultiplication maps of $Y$, applied blockwise at the finest partition level as in Definition \ref{Q-coalg}, commute with the $C(\uS)^p$-coordinates and with the reindexing maps $\Psi_{k,(p,q)}^*$; the remaining cofaces and the codegeneracies are handled exactly as for $m$ in Proposition \ref{prop:box-monoid}.
    
    Associativity and unitality of $\mu$ again follows from a similar straightforward modification of the  argument from \cite[Proposition 5.8]{Clark-1} (see also Corollary 8.4 and Remark 8.5 in  \cite{Clark-1}).
\end{proof}

\begin{remark}
    No connectivity hypotheses on $Y$ are required above: the maps $\mu_{p,q}$ are built only from the canonical comparison $\bigvee \to \prod$ and from projections, never their inverses. Connectivity enters only when comparing $\Tot C(Y)$ with other models for the derived primitives (compare \cite[Corollary 8.4]{Clark-1}, where $0$-connectedness is assumed).
\end{remark}

\begin{remark}
    A straightforward modification to the proofs of Theorems \ref{thm:main-1} and \ref{thm:main-2} shows that $\Tot C(\Cal{Q})$ is an $A_{\infty}$-operad and, if $Y$ is a $\Cal{Q}$-coalgebra, that $C(Y)$ is an algebra over this operad. In particular, our constructions provide a point-set model for $\Tot C(\Cal{Q})$ as a (homotopy coherent) operad which is \textit{Koszul dual} to $\Cal{Q}$, and $Y\mapsto \Tot C(Y)$ provides a comparison from $\Cal{Q}$-coalgebras to algebras over this operad (see \cite{FG}). 
\end{remark}

\appendix

\section{Proof of Theorem \texorpdfstring{\ref{thm:main-3}}{1.3}}\label{app:comparison}

This appendix proves Theorem \ref{thm:main-3}, as Theorem \ref{thm:main} in \S\ref{sec:main}; the intervening sections construct the objects named in it.

\subsection{Outline of the proof}\label{sec:app-outline}
Sections \ref{sec:cotensors}--\ref{sec:theta} set up the weighted-tree spaces $\BC(T)$ underlying Ching's bar construction on the commutative operad, in the coordinates of this document. We prove $\BC(n)\cong\Par(n) = |P(n)|$ by an explicit $\Sigma_n$-equivariant \emph{levelling homeomorphism} $\Thet$ (Theorem \ref{thm:levelling}), and write $\Omw(\uS) := \Map(\BC(-),S)$ for the arity-wise dual, which is Ching's model for $\partial_*\Id$.

Sections \ref{sec:grafting}--\ref{sec:bridge} express the cutting maps and the strict operad structure of \cite{Ching-thesis} in these coordinates, and connect the two models by a zig-zag of canonical maps. Both structures then appear as cutting weighted trees along flags and rescaling the pieces affinely. They differ only in where the cuts fall. Ching cuts each branch at the parent vertex of its block, while the structure of Theorem \ref{thm:main-1}, transported through $\Thet$, cuts at fixed heights.

Sections \ref{sec:Nlev}--\ref{sec:main} define, for each flag, spaces of \emph{framed} and \emph{windowed cut systems} containing both prescriptions among their cells. The windowed systems form an $\Nl$-operad $\Kop$ whose composition is affine substitution, everywhere defined and strictly associative, with component spaces contracting onto Ching's cuts. The tautological $\Kop$-algebra $\tau$ on $\Omw(\uS)$ restricts to Ching's structure along $o\colon \Oper\to\Kop$, and to the uniform-cut $\Lop$-algebra $\zeta$ along $\sigma\colon\Lop\to\Kop$. The latter is strictly intertwined with $\lambda^{\mathrm{lev}}$ through $j\colon\Lop\to\Aop$ and $\beta$.

We use the conventions of this document, with Warning \ref{warn:reversal} as the dictionary to those of \cite{Ching-thesis}. The results of \cite{Ching-thesis} are taken as given. Throughout, ``level slot'' refers to the grading of the $\Nl$-formalism: composition of an $\ell$-cell and an $\ell'$-cell happens at one of the $\ell$ levels and refines that entire level, producing an $(\ell+\ell'-1)$-cell.

\subsection{Motivation: loop spaces, one level up}\label{sec:loop-analogy}
Every construction below is the corresponding construction for based loops, one level up. \emph{Moore loops}, of variable length and concatenated end to end, form a strictly associative monoid; loops on a \emph{fixed interval} concatenate by squeezing both into $[0,1]$ at the cost of a reparametrization point, canonically $\tfrac12$, and form only an $A_\infty$-monoid. Ching's cocomposition is the Moore convention, and the structure of Theorem \ref{thm:main-1} the fixed-interval one, the box-product colimit forcing every cut to height $\tfrac12$ (Construction \ref{con:halfcut}). What follows transports the classical comparison from concatenation to cutting:
\begin{center}
\renewcommand{\arraystretch}{1.3}
\begin{tabular}{|p{0.40\textwidth}|p{0.50\textwidth}|}
\hline
\textbf{Loops} & \textbf{Weighted trees} \\
\hline
concatenation of loops & cocomposition by cutting weighted trees (\S\ref{sec:cutmaps}) \\
\hline
Moore loops (variable intervals) & Ching's tree cuts $\nu^{\mathrm{tr}}$ (Theorem \ref{thm:omw-operad}) \\
\hline
fixed interval, reparametrize at $\tfrac12$ & the half-cut $\psi_0$ and $\lambda^{\mathrm{lev}}$ (Construction \ref{con:halfcut}) \\
\hline
reparametrization choices & cut labellings: the levels operad $\Lop$ (Definition \ref{def:L}) \\
\hline
carrying interval lengths as data & windows as arguments: the operad $\Kop$ (Definition \ref{def:K}) \\
\hline
contractible choice spaces & contractible cell spaces (Proposition \ref{prop:K-operad}) \\
\hline
Moore loops $\simeq \Omega X$ as $A_\infty$-monoids & Theorem \ref{thm:main} \\
\hline
\end{tabular}
\end{center}
Two features have no analogue for loops. A cut at a given height is only partially defined, valid where the prescribed blocks are alive at that height and extended by the basepoint elsewhere, so continuity across the degeneration walls of $\BC(n)$ has to be arranged; this is what the axioms (W3), (W5), (W6) of Definition \ref{def:WCS} do. And the comparison is made on the structuring operads rather than on the objects.

\subsection{Preliminaries}\label{sec:cotensors} The category $\Spt$ is enriched, tensored and cotensored over $\Top$; for $X\in \Top$ and $Y\in \Spt$ we write $Y^{X}\in \Spt$ for the cotensor and, following the (mild) abuse in Proposition \ref{prop:model-equiv} below, we write $\Map(X, Y) := Y^{X}$, reserving $\Map(-,-)$ between two spectra for the mapping spectrum. The cotensor satisfies the standard adjunction isomorphisms, natural in all variables:
\begin{equation}\label{eq:cotensor}
    Y^{S^{0}} \cong Y, \quad\quad Y^{X\wedge X'} \cong (Y^{X'})^{X}, \quad\quad Y^{\colim_i X_i} \cong \lim_i Y^{X_i},
\end{equation}
the last for any small diagram of spaces. In particular, for a finite pointed set $F$ (discrete space) we have $S^{F}\cong \prod_{F\smallsetminus \basept} S$, the product of copies of the sphere spectrum indexed by the non-basepoint elements. For a finite set $T$ we write $P(T)$ for the pointed simplicial set of chains of partitions of $T$, as in Definition \ref{def:par-poset}.

For a pointed simplicial set $Z_{\bullet}$ (regarded as a simplicial discrete space) we write $|Z_{\bullet}|$ for its geometric realization and
\[ \fat{Z_{\bullet}} := \Big( \coprod_{p\geq 0} Z_p \wedge (\Delta^{p})_{+} \Big) \Big/ \big( (d_i z, u) \sim (z, \delta^i u) \big) \]
for its \emph{fat realization}, using only the face relations. The quotient $\fat{Z_{\bullet}}\to |Z_{\bullet}|$, which imposes the degeneracies, is a based homotopy equivalence (see \cite[Appendix A]{Segal-cat-coh} for the unbased statement; $\fat{\basept}$ is contractible). Write $\Tot_{\Delta} := \Map_{\Delta}(\Delta^{\bullet},-)$ for the unrestricted totalization.

We record a compactness principle, used below whenever a formula into a smash product degenerates on one factor: if $X, X'$ are compact pointed spaces and $x_i \to \basept$ in $X$, then $x_i \wedge x'_i \to \basept$ in $X\wedge X'$ for \emph{any} $(x'_i)$ in $X'$, by the tube lemma. Every space to which we apply it is a quotient of a finite disjoint union of compact polytopes, hence compact.

We coordinatize the topological $k$-simplex by \emph{heights}:
\[ \Delta^k \;\cong\; \nabla^k := \{ u=(u_1,\dots,u_k) \in [0,1]^k : u_1\leq u_2\leq \cdots \leq u_k \},\]
via $u_j = t_0+\cdots+t_{j-1}$ for barycentric coordinates $(t_0,\dots,t_k)$. In these coordinates the cofaces and codegeneracies are
\begin{align}
    \delta^0(u_1,\dots,u_{k-1}) &= (0, u_1,\dots,u_{k-1}), \notag \\
    \delta^i(u_1,\dots,u_{k-1}) &= (u_1,\dots,u_i,u_i,\dots,u_{k-1}) \quad (0<i<k), \label{eq:heights-cofaces}\\
    \delta^k(u_1,\dots,u_{k-1}) &= (u_1,\dots,u_{k-1}, 1), \notag
\end{align}
and $\sigma^{j}$ deletes $u_{j+1}$: $\delta^0$ prepends height $0$, $\delta^k$ appends height $1$, and interior cofaces duplicate a height.

\begin{lemma}\label{lem:tot-is-map}
Let $Z_{\bullet}$ be a pointed simplicial set with finite levels, and let $\Map(Z_{\bullet}, S)$ denote the cosimplicial spectrum $p\mapsto S^{Z_p} \cong \prod_{Z_p\smallsetminus \basept} S$, with cofaces $(d_i)^{*}$ and codegeneracies $(s_j)^{*}$. There are isomorphisms of spectra, natural in $Z_{\bullet}$,
\[ \Tot \Map(Z_{\bullet}, S) \;\cong\; \Map\big( \fat{Z_{\bullet}}, S \big), \quad\quad
   \Tot_{\Delta} \Map(Z_{\bullet}, S) \;\cong\; \Map\big( |Z_{\bullet}|, S \big),\]
under which the natural map $\Tot_{\Delta}\to \Tot$, restricting $\Delta$-natural families to $\Dres$-natural families, is identified with restriction along $\fat{Z_{\bullet}}\to |Z_{\bullet}|$. In particular $\Tot_{\Delta} \Map(Z_{\bullet}, S)\to \Tot \Map(Z_{\bullet},S)$ is a weak equivalence.
\end{lemma}

\begin{proof}
By \eqref{eq:cotensor}, $S^{(-)}$ sends colimits of spaces to limits of spectra, so the limit over $\Dres$ defining $\Tot\Map(Z_{\bullet},S)$ is $S$ cotensored by the colimit defining $\fat{Z_{\bullet}}$, and likewise over $\Delta$ with $|Z_{\bullet}|$; both restriction maps are induced by the inclusion of $\Dres$-indexed into $\Delta$-indexed data. The last claim holds since $\fat{Z_\bullet}\to|Z_\bullet|$ is a based homotopy equivalence and $S^{(-)}$ preserves such.
\end{proof}

Combining Proposition \ref{prop:model-equiv} with Lemma \ref{lem:tot-is-map} applied to $Z_{\bullet} = P(n)_{\bullet}$:
\begin{equation}\label{eq:tot-fat-model}
    \Tot C(\uS)[n] \;\cong\; \Map\big(\fat{P(n)}, S\big) \;\simeq\; \Map\big(\Par(n), S\big)\;\simeq\; \partial_n\Id .
\end{equation}

\subsection{Trees as laminar families}\label{sec:laminar}

\begin{definition}\label{def:laminar}
Let $T$ be a finite nonempty set. A \emph{tree on $T$} is a collection $\TreeF$ of nonempty subsets of $T$ such that:
\begin{enumerate}[label=(\roman*)]
    \item $T\in \TreeF$ and $\{x\}\in \TreeF$ for every $x\in T$;
    \item any two members of $\TreeF$ are either nested or disjoint.
\end{enumerate}
Members of $\TreeF$ are called \emph{vertices}; $T$ is the \emph{root vertex}, singletons are \emph{leaves}, and non-singleton vertices are \emph{internal}. We write $V(\TreeF)$ for the set of internal vertices. For $B\in \TreeF$ with $B\neq T$, the \emph{parent} $p(B)$ is the smallest member of $\TreeF$ properly containing $B$; the \emph{children} of $B\in\TreeF$ are the maximal members of $\TreeF$ properly contained in $B$. For a leaf $\{x\}$ we abbreviate $v(x) := p(\{x\})$.
\end{definition}

\begin{lemma}\label{lem:laminar-facts}
The children of a non-singleton $B\in\TreeF$ are pairwise disjoint with union $B$, and there are at least two of them. In particular a tree has no ``unary vertices'', and $|V(\TreeF)|\leq |T|-1$, with equality iff every internal vertex is binary.
\end{lemma}

Both claims are elementary from maximality and laminarity, the count by induction over the root's children. We picture $\TreeF$ as a rooted tree drawn with the root at the bottom and the leaves at the top: a \emph{root edge} hangs below the root vertex, each non-root $B$ is joined to $p(B)$ by an edge $e_B$, and \emph{leaf edges} terminate in the elements of $T$. Bijections $T\cong T'$ act on trees in the evident way.

\begin{warning}[Orientation dictionary]\label{warn:reversal}
Our chains run \emph{coarse to fine}: $\lambda_0=\mathsf{min}$ is the one-block partition and $\lambda_k=\mathsf{max}$ the discrete one. In \cite{Ching-thesis} the opposite convention is used, chains running from $\hat{0}$ discrete to $\hat{1}$ one-block. The dictionary is the levelwise order-reversal
\[ (\lambda_0\leq \cdots \leq \lambda_k) \longleftrightarrow (\lambda_k\geq \cdots \geq \lambda_0), \]
a simplicial isomorphism onto the reversed simplicial set, under which our $d_i$ is Ching's $d_{k-i}$ and our $d_0$ his top face. Our chains thus read a tree from the root upward. The dictionary is applied silently whenever \cite{Ching-thesis} is invoked.
\end{warning}

\subsection{Weightings by heights}\label{sec:weightings}

\begin{definition}\label{def:weighting}
A \emph{weighting} of a tree $\TreeF$ on $T$ is a function $h\colon V(\TreeF)\to [0,1]$ that is order-preserving: $B\subseteq B'$ implies $h(B')\leq h(B)$ for internal vertices $B, B'$. We write $\wt(\TreeF)$ for the space of weightings, topologized as a subspace of $[0,1]^{V(\TreeF)}$. We picture $h(B)$ as the height of the vertex $B$, the root of the tree sitting at height $0$ and all leaves at height $1$; thus the root edge has length $h(T)$, the edge $e_B$ has length $h(B)-h(p(B))$, and the leaf edge above $x\in T$ has length $1 - h(v(x))$. (For $|T|=1$ there are no internal vertices and $\wt(\TreeF)=\{*\}$; the unique edge is both root and leaf edge and has length $1$.)
\end{definition}

\begin{lemma}\label{lem:weightings-disc}
$\wt(\TreeF)$ is a compact convex polytope of dimension $|V(\TreeF)|$ with non\-empty interior in $[0,1]^{V(\TreeF)}$; in particular $\wt(\TreeF)\cong D^{|V(\TreeF)|}$.
\end{lemma}

Indeed $\wt(\TreeF)$ is cut out of the cube by the inequalities $h(p(B))\leq h(B)$, and contains the strictly order-preserving weightings with values in $(0,1)$, an open subset; a compact convex subset of Euclidean space with nonempty interior is a disc.

\begin{remark}\label{rem:ching-weightings}
In \cite[Definitions 3.6--3.10]{Ching-thesis} weightings are nonnegative edge lengths with all root-to-leaf distances $1$; passing to cumulative heights from the root identifies that description with Definition \ref{def:weighting}, and Lemma \ref{lem:weightings-disc} recovers \cite[Lemma 3.9]{Ching-thesis}. Our trees have no unary vertices; those appear in \cite{Ching-thesis} through degenerate simplices of the simplicial bar construction, and are accounted for below by $\Thet$ being defined on the thin realization $|P(n)|$, where degeneracies are already collapsed.
\end{remark}

\subsection{The bar construction as a space of weighted trees}\label{sec:Bdef}

\begin{definition}\label{def:B}
For a finite nonempty set $T$, let
\[ \BC(T) := \Big( \coprod_{\TreeF \text{ tree on } T} \wt(\TreeF) \Big) \Big/ \sim, \]
pointed by the class $\basept$ described in (b), where $\sim$ is the equivalence relation generated by:
\begin{enumerate}[label=(\alph*)]
    \item (\emph{collapse}) for an internal non-root vertex $B$ with $h(B)=h(p(B))$:
    \[(\TreeF, h) \sim (\TreeF\smallsetminus\{B\},\, h|_{V(\TreeF)\smallsetminus\{B\}});\]
    \item (\emph{degeneracy}) every $(\TreeF, h)$ with $h(T)=0$ (zero root edge) or with $h(v(x))=1$ for some $x\in T$ (a zero leaf edge) is equivalent to the basepoint.
\end{enumerate}
We write $\BC(n) := \BC(\mathbf{n})$; bijections of $T$ act on $\BC(T)$, and in particular $\BC(n)$ is a pointed $\Sigma_n$-space. Note $\BC(T)\cong S^0$ when $|T|=1$.
\end{definition}

Removing $B$ preserves Definition \ref{def:laminar}, and (b) is stable under (a), since a collapse never changes $h(T)$ and a collapsed $v(x)$ leaves the new leaf-parent at the same height. We identify canonical representatives.

\begin{definition}\label{def:reduced}
A weighting $h$ of $\TreeF$ is \emph{reduced} if $h(B) > h(p(B))$ for every internal non-root vertex $B$, i.e., all internal edges have positive length.
\end{definition}

Notice that in a reduced pair a vertex of maximal height has only leaf children, since an internal child would be strictly higher. We use this twice below.

\begin{lemma}\label{lem:canonical-reps}
For $(\TreeF, h)$ set
\begin{align*}
\TreeF^{\mathrm{red}} &:= \{T\}\cup \{\text{leaves}\} \cup \{B\in V(\TreeF),\ B\neq T : h(B) > h(p(B))\}, \\
h^{\mathrm{red}} &:= h|_{V(\TreeF^{\mathrm{red}})}.
\end{align*}
Then $(\TreeF^{\mathrm{red}}, h^{\mathrm{red}})$ is reduced, and it is the unique reduced pair equivalent to $(\TreeF,h)$ under the collapse relation (a). Consequently two pairs are equivalent under (a) if and only if they have the same reduction, and the non-basepoint classes of $\BC(T)$ are in bijection with reduced pairs satisfying $h(T)>0$ and $h(v(x))<1$ for all $x\in T$.
\end{lemma}

Indeed, along any ancestor chain the heights form weakly decreasing constant runs, and $\TreeF^{\mathrm{red}}$ keeps exactly the root-most vertex of each; each move (a) preserves the reduction, and finitely many of them connect $(\TreeF,h)$ to it.

\begin{lemma}\label{lem:B-cpt-hd}
$\BC(T)$ is compact Hausdorff and the quotient map $q\colon \coprod_{\TreeF} \wt(\TreeF) \to \BC(T)$ is closed.
\end{lemma}

\begin{proof}
The source is compact Hausdorff, and the equivalence relation is closed in its square: the basepoint class is cut out by the closed conditions of Definition \ref{def:B}(b), and ``same reduced pair'' is a finite union of graphs of affine partial maps on closed polyhedral subsets (Lemma \ref{lem:canonical-reps}). A quotient of a compact Hausdorff space by a closed equivalence relation is compact Hausdorff with closed quotient map.
\end{proof}

\subsection{The alive partition}\label{sec:alive}

\begin{definition}\label{def:alive}
We adopt once and for all the conventions $h(p(T)) := 0$ and $h(\{x\}) := 1$ for leaves. Let $(\TreeF, h)$ be a weighted tree on $T$ and let $t\in (0,1]$. A vertex $B\in \TreeF$ is \emph{alive at $t$} if $h(B)\geq t$ and $h(p(B)) < t$; a leaf is thus alive at $t$ if and only if $h(v(x))<t$. The vertices alive at $t$ form a partition $\alive{(\TreeF,h)}{t}$ of $T$, consisting of the maximal $B$ with $h(B)\geq t$ together with the leaves $\{x\}$ with $h(v(x)) < t$.
\end{definition}

\begin{lemma}\label{lem:alive-facts}
Fix $(\TreeF, h)$ on $T$ with $h$ reduced. Then:
\begin{enumerate}[label=(\arabic*)]
\item $\alive{(\TreeF,h)}{t}$ is a partition of $T$ for each $t\in(0,1]$, refined by $\alive{(\TreeF,h)}{t'}$ for $t\leq t'$;
\item $\alive{(\TreeF,h)}{t} = \mathsf{min}$ for $0 < t\leq h(T)$, and $\alive{(\TreeF,h)}{t} = \mathsf{max}$ for $\max_x h(v(x)) < t \leq 1$;
\item $t\mapsto \alive{(\TreeF,h)}{t}$ is left-continuous with jumps exactly at the heights $\{h(B) : B\in V(\TreeF)\}$: for $t$ in the half-open interval $(h', h]$ between consecutive distinct vertex heights $h' < h$ (with $h'=0$ below the least), the alive partition is constant, and its blocks are the vertices $B$ with $h(B)\geq t > h(p(B))$;
\item $(\TreeF, h)$ is recovered from the function $t\mapsto \alive{(\TreeF,h)}{t}$: the internal vertices are the non-singleton blocks of the alive partitions, and $h(B) = \max\{t : B \text{ is a block of } \alive{(\TreeF,h)}{t}\}$.
\end{enumerate}
\end{lemma}

All four are elementary from the definitions; reducedness enters only through the exactness of the jump set in (3) and the nonemptiness of the range $h(p(B))<t\leq h(B)$ in (4).

\subsection{The levelling homeomorphism}\label{sec:theta}

Recall the height coordinates $\nabla^k\cong \Delta^k$ of \S\ref{sec:conventions}. Let $\alpha = (\lambda_0\leq \cdots \leq \lambda_k)\in P(n)^{\circ}_k$ be a non-basepoint chain and put
\[\TreeF(\alpha) := \{ B \subseteq \mathbf{n} : B \text{ a block of some } \lambda_j \},\]
a tree on $\mathbf{n}$ in the sense of Definition \ref{def:laminar}: blocks of comparable partitions are nested or disjoint, $\mathbf{n}$ is the block of $\lambda_0$, and the singletons are the blocks of $\lambda_k$.

Two combinatorial facts about $\TreeF(\alpha)$ will be used throughout. Both are immediate from disjointness of blocks within a partition and minimality of the parent.

\emph{Interval property.} If $B$ is a block of $\lambda_a$ and of $\lambda_b$ then it is a block of every $\lambda_j$ in between. We may therefore write $a(B)\leq b(B)$ for the first and last levels of a block $B$.

\emph{Separation of levels.} If $B\subsetneq B'$ in $\TreeF(\alpha)$ then $b(B')<a(B)$, and moreover
\begin{equation}\label{eq:parent-level}
a(B) = b(p(B))+1 \quad\quad (B\neq\mathbf{n}).
\end{equation}

\begin{construction}\label{con:theta}
Define, for $\alpha\in P(n)^{\circ}_k$ and $u\in \nabla^k$,
\[ \Thet(\alpha, u) := \big(\TreeF(\alpha),\, h_{\alpha,u}\big), \quad\quad h_{\alpha, u}(B) := u_{\,b(B)+1} \quad (B\in V(\TreeF(\alpha))),\]
using that $b(B)\leq k-1$ for non-singleton $B$ (level $k$ consists of singletons). Basepoint simplices are sent to the basepoint. Set $u_{k+1}:=1$ and $u_0 := 0$.
\end{construction}

Each block of $\alpha$ is placed at the height of its last level. The comparison of \S\ref{sec:j-transport} uses the following description of the alive partitions of $\Thet(\alpha,u)$ rather than Theorem \ref{thm:levelling} itself.

\begin{lemma}\label{lem:theta-dictionary}
Let $\alpha\in P(n)^{\circ}_k$ be nondegenerate and let $u\in\nabla^k$ be interior, so that $0<u_1<\cdots<u_k<1$. Then:
\begin{enumerate}[label=(\arabic*)]
    \item $h_{\alpha,u}$ is a reduced weighting of $\TreeF(\alpha)$, with $h_{\alpha,u}(\mathbf{n})>0$ and $h_{\alpha,u}(v(x))<1$ for every $x$;
    \item $\alive{\Thet(\alpha,u)}{t} = \lambda_j$ for every $t\in(u_j, u_{j+1}]$ and $0\leq j\leq k-1$;
    \item the distinct values of $h_{\alpha,u}$ are exactly $u_1<\cdots<u_k$.
\end{enumerate}
In particular $(\alpha,u)$ is recovered from $\Thet(\alpha,u)$.
\end{lemma}

\begin{proof}
(1) For internal $B\subsetneq B'$, separation of levels gives $b(B')<a(B)\leq b(B)$, so $h_{\alpha,u}(B')\leq h_{\alpha,u}(B)$ by monotonicity of $u$, and the inequality is strict when $B' = p(B)$ by \eqref{eq:parent-level} and strictness of $u$. Also $h_{\alpha,u}(\mathbf{n}) \geq u_1 > 0$ and $h_{\alpha,u}(v(x)) \leq u_k < 1$.

(2) By \eqref{eq:parent-level}, a non-singleton $B$ is a block of $\lambda_j$ if and only if $b(p(B)) < j \leq b(B)$, which by strictness of $u$ says $h(p(B)) < u_{j+1}\leq h(B)$, that is, that $B$ is alive at $u_{j+1}$. Singletons persist through level $k$, and a singleton $\{x\}$ is a block of $\lambda_j$ iff $h(v(x))\leq u_j < u_{j+1}$, so the same holds for them. This gives the case $t = u_{j+1}$, and the alive partition is constant on $(u_j,u_{j+1}]$ by Lemma \ref{lem:alive-facts}(3).

(3) Each $u_j$ is attained: since $\lambda_{j-1}\neq\lambda_j$, some block of $\lambda_{j-1}$ is absent from $\lambda_j$, and such a block has last level $j-1$, hence height $u_j$. It is non-singleton because singletons persist through level $k$. No other values occur, by the definition of $h_{\alpha,u}$.
\end{proof}

\begin{theorem}\label{thm:levelling}
The maps of Construction \ref{con:theta} descend to a $\Sigma_n$-equivariant homeomorphism
\[ \Thet = \Thet_n \colon \Par(n) = |P(n)| \xrightarrow{\;\cong\;} \BC(n).\]
\end{theorem}

That the two spaces are homeomorphic is not new (\cite{Fresse-Koszul}, \cite{Heu-Moer-trees}, Remark \ref{rem:theta-vs-ching}); what is used below is that the map of Construction \ref{con:theta} is a homeomorphism.

\begin{proof}
(I) \emph{$h_{\alpha,u}$ is a weighting for every $u\in\nabla^k$.} For $B\subsetneq B'$, separation of levels gives $b(B')<a(B)\leq b(B)$, so $h_{\alpha,u}(B')\leq h_{\alpha,u}(B)$ by monotonicity of $u$.

(II) \emph{$\Thet$ respects the relations.} The realization $|P(n)|$ is the quotient of $\coprod_k P(n)^{\circ}_k\times \nabla^k$ by the relations $(d_i\alpha, u)\sim (\alpha, \delta^i u)$ and $(s_j \alpha, u)\sim (\alpha, \sigma^j u)$. A relation whose left side is a basepoint simplex sends the right-hand points to the basepoint. We check these against Definition \ref{def:B}.

\emph{Degeneracies and inner faces.} These are bookkeeping with the indices $b(B)$: $s_j$ shifts $b_\alpha(B)$ by one exactly when $b_\alpha(B)\geq j$ while $\sigma^j$ deletes $u_{j+1}$, giving $\Thet(s_j\alpha,u)=\Thet(\alpha,\sigma^ju)$ on the nose; and for $0<i<k$, deleting the $i$-th entry removes exactly the blocks with $a(B)=b(B)=i$, no two of which are nested, each of which acquires a zero-length parent edge under $\delta^i$ (its height and its parent's are both $u_i$, using $a(B)=b(p(B))+1$), so that the collapse relation (a) identifies $\Thet(\alpha,\delta^iu)$ with $\Thet(d_i\alpha,u)$.

\emph{Outer faces.} For $i=0$: $\delta^0$ prepends $0$. If $\lambda_1 \neq \mathsf{min}$ then $d_0\alpha=\basept$; on the other side $b(\mathbf{n}) = 0$, so $h_{\alpha, \delta^0 u}(\mathbf{n}) = (\delta^0 u)_1 = 0$: the root edge has length zero and relation (b) sends the point to the basepoint. If $\lambda_1 = \mathsf{min}$ then every block of $\alpha$ has $b_\alpha(B)\geq 1$ (for $B = \mathbf{n}$ directly; for $B\neq \mathbf{n}$ since $b_\alpha(B)\geq a_\alpha(B) = b_\alpha(p(B))+1\geq 1$), the blocks of $d_0\alpha$ are those of $\alpha$ with $b_{d_0\alpha}(B) = b_\alpha(B)-1$, and $(\delta^0u)_{b_\alpha(B)+1} = u_{b_\alpha(B)}$, so all heights match. For $i=k$: $\delta^k$ appends $1$. If $\lambda_{k-1}\neq \mathsf{max}$, some singleton $\{x\}$ has $a(\{x\}) = k$, so its parent satisfies $b(v(x)) = k-1$ and $h_{\alpha,\delta^ku}(v(x)) = (\delta^ku)_k = 1$: a zero leaf edge, hence basepoint by (b), matching $d_k\alpha = \basept$. If $\lambda_{k-1}=\mathsf{max}$ then all non-singleton blocks have $b(B)\leq k-2$ and heights match verbatim.

Thus $\Thet$ descends to a continuous pointed map $|P(n)|\to \BC(n)$ (continuity: on each closed cell $\{\alpha\}\times\nabla^k$ the map to $\wt(\TreeF(\alpha))$ is the linear coordinate-selection $u\mapsto (u_{b(B)+1})_B$, and both quotients carry the quotient topology). Equivariance is clear, as $\sigma\in\Sigma_n$ permutes chains, blocks, and trees compatibly and does not touch heights.

(III) \emph{Bijectivity.} Every non-basepoint class of $|P(n)|$ has a unique representative $(\alpha, u)$ with $\alpha$ nondegenerate and $u$ interior, and every non-basepoint class of $\BC(n)$ a unique reduced representative with $h(\mathbf{n})>0$ and all $h(v(x))<1$ (Lemma \ref{lem:canonical-reps}). By Lemma \ref{lem:theta-dictionary}(1), $\Thet$ carries the first set into the second, and by Lemma \ref{lem:theta-dictionary}(2) and (3) it is injective on it. It remains to construct a preimage.

Let $(\TreeF,h)$ be reduced with $h(\mathbf{n})>0$ and all $h(v(x))<1$. Let $w_1<\cdots<w_k$ be the distinct values of $h$ on $V(\TreeF)$ and put $u := (w_1,\dots,w_k)$. This is an interior point of $\nabla^k$: we have $w_1 = h(\mathbf{n})>0$, and a vertex of maximal height has only leaf children (Definition \ref{def:reduced}), so $w_k = \max_x h(v(x))<1$. Setting
\[\lambda_j := \alive{(\TreeF,h)}{w_{j+1}} \;\;(0\leq j\leq k-1), \quad\quad \lambda_k := \mathsf{max}\]
gives a nondegenerate chain starting at $\mathsf{min}$, since each $w_j$ is the height of a vertex alive at $w_j$ but not at $w_{j+1}$. By Lemma \ref{lem:alive-facts}(4) we get $\TreeF(\alpha(\TreeF,h)) = \TreeF$ and $h_{\alpha(\TreeF,h),u} = h$, because a non-singleton $B$ is a block of $\lambda_j$ exactly when $h(p(B))<w_{j+1}\leq h(B)$.

(IV) $|P(n)|$ is compact, since $P(n)$ has finitely many nondegenerate simplices, and $\BC(n)$ is Hausdorff by Lemma \ref{lem:B-cpt-hd}. A continuous bijection from a compact space to a Hausdorff space is a homeomorphism.
\end{proof}

\begin{remark}\label{rem:theta-vs-ching}
Under the dictionary of Warning \ref{warn:reversal}, Theorem \ref{thm:levelling} is the composite of \cite[Proposition 4.13]{Ching-thesis} and \cite[Lemma 8.6]{Ching-thesis} for $P = \uSz$: the proof of the former weights the tree of a chain by cumulative barycentric coordinates, which in height coordinates is Construction \ref{con:theta}. We have included a proof because \S\ref{sec:j-transport} needs Lemma \ref{lem:theta-dictionary}.
\end{remark}

\subsection{Grafting and degrafting}\label{sec:grafting}

Fix $k\geq 1$ and a partition $\alpha = \{T_1,\dots,T_n\}$ of $\mathbf{k}$, allowing $n=1$ and the partition into singletons; we identify $\alpha$ with the $2$-simplex $\mathsf{min}\leq \alpha\leq \mathsf{max}$ of $P(k)$. Let $q = q_\alpha\colon \mathbf{k}\to \mathbf{n}$ denote the block-collapse map, $q(x) = i$ for $x\in T_i$.

\begin{definition}\label{def:graft}
Given a tree $\Cal{G}$ on $\mathbf{n}$ and trees $\Cal{U}_i$ on $T_i$ ($1\leq i\leq n$), their \emph{graft} is
\[ \Cal{G}\circ(\Cal{U}_1,\dots,\Cal{U}_n) := \{ q^{-1}(D) : D\in \Cal{G}\} \cup \Cal{U}_1\cup \cdots \cup \Cal{U}_n,\]
a collection of subsets of $\mathbf{k}$.
\end{definition}

\begin{lemma}\label{lem:graft-degraft}
$\Cal{G}\circ(\Cal{U}_i)$ is a tree on $\mathbf{k}$ containing every $T_i$. Conversely, a tree $\Cal{V}$ on $\mathbf{k}$ is of the form $\Cal{G}\circ (\Cal{U}_i)$ if and only if $T_i\in \Cal{V}$ for every $i$, and in that case the \emph{degrafting} is unique:
\[ \Cal{U}_i = \{B\in\Cal{V} : B\subseteq T_i\}, \quad\quad \Cal{G} = \{ q(B) : B\in \Cal{V},\; B\not\subsetneq T_i \text{ for all } i \}.\]
\end{lemma}

\begin{proof}
This is the binary grafting and degrafting of \cite[Definition 4.14, Lemma 4.15]{Ching-thesis}, iterated over the blocks of $\alpha$; in the laminar description the verification is a direct check of Definition \ref{def:laminar} (a member $B\in\Cal{V}$ with $B\not\subseteq T_i$ for all $i$ is a union of blocks by nestedness, so $B = q^{-1}(q(B))$), which we omit.
\end{proof}

\begin{remark}\label{rem:graft-vs-ching}
Compare the uniqueness clause of Lemma \ref{lem:tech}: a levelled chain $\gamma'$ decomposes with respect to $\alpha$ if and only if $\lambda_p = \alpha$, which under $\Thet$ says that all the $T_i$ occur as vertices \emph{at a common level}. Lemma \ref{lem:graft-degraft} drops that requirement, and the difference between the two requirements accounts, in the sense established in \S\S\ref{sec:flags}--\ref{sec:main}, for the difference between the two operad structures.
\end{remark}

\subsection{The cut maps}\label{sec:cutmaps}

\begin{definition}\label{def:cut}
Fix $\alpha = \{T_1,\dots,T_n\}$ as above. Define, for each tree $\Cal{V}$ on $\mathbf{k}$ and each weighting $h\in \wt(\Cal{V})$,
\[ c^{\mathrm{tr}}_{\alpha}(\Cal{V}, h) \in \BC(n)\wedge \BC(T_1)\wedge \cdots \wedge \BC(T_n)\]
as follows. If some non-singleton $T_i\notin \Cal{V}$, set $c^{\mathrm{tr}}_\alpha(\Cal{V},h) := \basept$. Otherwise let $(\Cal{G}, (\Cal{U}_i))$ be the degrafting of Lemma \ref{lem:graft-degraft} and set
\[ c^{\mathrm{tr}}_\alpha(\Cal{V},h) := \big( (\Cal{G}, h_{\mathrm{bot}}) ,\, (\Cal{U}_1, h_1), \dots, (\Cal{U}_n, h_n) \big), \]
where, writing $\pi_i := h(p(T_i))$ (with the conventions $p(\mathbf{k})$ of height $0$, and, for a singleton $T_i = \{x\}$, $(\Cal{U}_i, h_i)$ the non-basepoint of $\BC(T_i)\cong S^0$):
\begin{align}
    h_{\mathrm{bot}}(D) &:= h\big(q^{-1}(D)\big) && (D\in V(\Cal{G})), \label{eq:cut-bot}\\
    h_i(B) &:= \frac{h(B) - \pi_i}{1 - \pi_i} && (B \in V(\Cal{U}_i),\ \text{$T_i$ non-singleton, } \pi_i < 1), \label{eq:cut-top}
\end{align}
and $c^{\mathrm{tr}}_\alpha(\Cal{V},h) := \basept$ if $\pi_i = 1$ for some non-singleton $T_i$.
\end{definition}

The bottom tree keeps its heights, and each top tree is rescaled by the affine homeomorphism $[\pi_i, 1]\to [0,1]$, so its root edge receives the rescaled length of $e_{T_i}$. This is the cocomposition of \cite[Definition 4.16]{Ching-thesis} in height coordinates. Both formulas land in weightings, order-preservation being inherited via $q^{-1}$ and preserved by affine maps, and $h_i(B)\in[0,1]$ since $h(B)\geq h(T_i)\geq\pi_i$. That the $c^{\mathrm{tr}}_{\alpha}$ descend to well-defined, continuous, based, $\Sigma_k$-equivariant maps
\[ c^{\mathrm{tr}}_{\alpha}\colon \BC(k)\longrightarrow \BC(n)\wedge \BC(T_1)\wedge\cdots\wedge \BC(T_n),\]
coassociative in the appropriate sense, is established in \cite[\S 4]{Ching-thesis}; the tree cuts are also instances of the framed cut maps of Proposition \ref{prop:FS-continuous} below (Example \ref{ex:frames}(1)).

\subsection{The tree cobar operad}\label{sec:coassoc}

\begin{definition}\label{def:omw}
Define the symmetric sequence $\Omw(\uS)$ in $\Spt$ by
\[ \Omw(\uS)[k] := \Map\big( \BC(k), S\big) \quad (k\geq 1), \quad\quad \Omw(\uS)[0] := \pt,\]
with $\Sigma_k$ acting by conjugation (via the action on $\BC(k)$), and define, for each partition $\alpha$ of $\mathbf{k}$ with profile $(n, (k_1,\dots,k_n))$, the entrywise composition map
\[ \xi_\alpha \colon \Omw(\uS)[n]\wedge \Omw(\uS)[k_1]\wedge \cdots \wedge \Omw(\uS)[k_n] \longrightarrow \Omw(\uS)[k]\]
as the smash pairing of cotensors $\Map(X, S)\wedge \Map(Y,S)\to \Map(X\wedge Y, S)$ (iterated; strictly associative and unital by the standard enriched-category isomorphisms in $\Spt$, using $S\wedge_S S\cong S$) followed by restriction along $c^{\mathrm{tr}}_\alpha$ and the identifications $\BC(T_i)\cong \BC(k_i)$ induced by the order-preserving bijections $T_i\cong \mathbf{k_i}$.
\end{definition}

\begin{theorem}[{Ching \cite[\S\S 4--6, \S 8]{Ching-thesis}}]\label{thm:omw-operad}
The maps
\[ \mu\colon \big(\Omw(\uS)\circ \Omw(\uS)\big)[k] = \bigvee_{\alpha\in P(k)_2} \big(\Omw(\uS)\otimes \Omw(\uS)\big)[\alpha] \xrightarrow{\ \vee_\alpha\, \xi_\alpha\ } \Omw(\uS)[k]\]
define a reduced operad structure on $\Omw(\uS)$, with unit $I\to \Omw(\uS)$ given by $S \cong \Map(S^0,S) = \Omw(\uS)[1]$.
\end{theorem}

Under the dictionary of Warning \ref{warn:reversal} and Remark \ref{rem:ching-weightings} this is Ching's operad $\partial_*\Id = \mathbb{D}B(\uSz)$ \cite[Corollary 8.8, Remark 8.9]{Ching-thesis}. Associativity of $\mu$ is coassociativity of the cutting maps \cite[Lemma 4.28, Theorem 4.29]{Ching-thesis}, and unitality is the counit property of the cuts along the two trivial partitions. The dualization is that of \cite[Corollary 5.3, Lemma 6.1, Proposition 6.4]{Ching-thesis}. Here it is the entrywise cotensor $\Map(-,S)$, because every entry of $\uSz$ is $S^0$. The structure is strict: associativity is an identity of affine rescalings.

\begin{proposition}\label{prop:omw-is-derivatives}
There are $\Sigma_k$-equivariant isomorphisms and weak equivalences
\[ \Omw(\uS)[k] = \Map(\BC(k), S) \xrightarrow[\cong]{\ \Thet_k^{*}\ } \Map\big(\Par(k), S\big) \simeq \partial_k \Id \quad (k \geq 1),\]
the second being the partition-poset model for the derivatives \eqref{der_n I}. In particular the underlying symmetric sequence of $\Omw(\uS)$ is a model for $\partial_*\Id$.
\end{proposition}

\begin{proof}
Immediate from Theorem \ref{thm:levelling} and \eqref{der_n I}.
\end{proof}

\subsection{The object-level bridge}\label{sec:bridge}

\begin{proposition}\label{prop:bridge}
For each $k\geq 1$ there are $\Sigma_k$-equivariant maps of spectra
\begin{equation}\label{eq:bridge}
 \Omw(\uS)[k] = \Map(\BC(k), S) \xrightarrow[\cong]{\,\Thet_k^{*}\,} \Map\big(|P(k)|, S\big) \xleftarrow[\cong]{\ \ } \Tot_{\Delta}\, C(\uS)[k] \xrightarrow{\,\mathrm{res}\,} \Tot C(\uS)[k],
\end{equation}
in which the first two maps are isomorphisms (Theorem \ref{thm:levelling} and Lemma \ref{lem:tot-is-map}) and the third, restriction of ends along $\Dres\subset \Delta$, is a weak equivalence (Lemma \ref{lem:tot-is-map} again, dual to $\fat{P(k)}\to |P(k)|$). These maps assemble to maps of symmetric sequences; we write
\[ \beta\colon \Omw(\uS) \longrightarrow \Tot C(\uS) \]
for the resulting composite (inverting the middle isomorphism), a levelwise weak equivalence.
\end{proposition}

\begin{proof}
Only equivariance and assembly remain: $\Thet_k$ is $\Sigma_k$-equivariant (Theorem \ref{thm:levelling}), the identifications of Lemma \ref{lem:tot-is-map} are natural in $Z_\bullet = P(k)_{\bullet}$ and hence equivariant, and $\mathrm{res}$ is natural.
\end{proof}

\begin{remark}\label{rem:phi-vs-ching52}
For $Q = \uS$ this is the precise form of the assertion of \cite[Remark 5.2]{Ching-thesis} that the cobar construction on a cooperad is isomorphic to the totalization of a cosimplicial cobar construction dual to the simplicial bar construction. Our route uses no ends over tree categories, since both sides are cotensors of $S$ by a space. The only non-isomorphism in \eqref{eq:bridge} is the fat-to-thin comparison. It is present because our $\Tot$ is the restricted totalization, whose dual is the fat realization.
\end{remark}

\begin{lemma}[Units]\label{lem:units-bridge}
Under $\beta$, the operad unit $S = \Map(S^0, S) = \Omw(\uS)[1]$ of Theorem \ref{thm:omw-operad} and the map $S = \Tot(\underline{I})[1]\to \Tot C(\uS)[1]$ induced by the $\boxcirc$-unit $u\colon \underline{I}\to C(\uS)$ of Proposition \ref{prop:box-monoid} correspond: $\beta[1]$ carries the identity element of $\Map(S^0,S)$ to the image of the unit of $\Tot(\underline{I})[1]$.
\end{lemma}

\begin{proof}
At $k = 1$ both maps are obtained by applying $\Map(-,S)$ to the same collapse $\fat{P(1)}\to |P(1)| = S^0$: the simplicial set $P(1)$ has a single non-basepoint simplex in each degree, all degenerate above degree $0$, and $u$ is the identity in arity $1$ in every cosimplicial degree.
\end{proof}

\subsection{\texorpdfstring{$\Nl$}{N-lev}-operads and the notion of agreement}\label{sec:Nlev}
We use the formalism of $\Nl$-operads of \cite[\S\S 6--7]{Clark-1}: an $\Nl$-operad $\Cal{P}$ in $\Top$ consists of spaces $\Cal{P}_{\ell}(\underline{p}; t)$ indexed by levels $\ell\geq 1$, profiles $\underline{p}\in \Ncirc{\ell}[t]$ and arities $t\geq 1$, with $\Sigma$-equivariance, composition and unit data as specified there. Three $\Nl$-operads play a role:
\begin{itemize}[leftmargin=2em]
    \item $\Oper = \Oper^{\Top}$, whose algebras in $\SymSeq$ are precisely reduced operads \cite[Definition 7.2, Proposition 7.12]{Clark-1};
    \item $\Aop = \coEnd(\SD^{\bullet}_{+})$, the coendomorphism operad of the cosimplicial object $\SD^{\bullet}_+$ \cite[Definition 7.7, Remark 7.17]{Clark-1}; the totalization of a $\boxcirc$-monoid is naturally an $\Aop$-algebra \cite[\S 8.1]{Clark-1}, and there is a levelwise weak equivalence $\rho\colon \Aop \xrightarrow{\ \sim\ } \Oper$ \cite[Proposition 7.16]{Clark-1};
    \item the cut-system operad $\Kop$ constructed below.
\end{itemize}
For an $\Nl$-operad $\Cal{P}$, a $\Cal{P}$-algebra is a symmetric sequence with action maps as in \cite[Definition 6.32]{Clark-1}. Following \cite[Definition 6.38]{Clark-1}, a symmetric sequence $\Cal{W}$ equipped with a $\Cal{P}$-algebra structure for some $\Cal{P}$ admitting a zig-zag of levelwise weak equivalences of $\Nl$-operads to $\Oper$ is called an \emph{$A_\infty$-operad}.

A \emph{cell} of an $\Nl$-operad is an element of one of its spaces, a single prescription for forming one composite. A cell of $\Oper$ is a chain of partitions, so $\Oper$ is discrete; a cell of $\Lop$ is a chain together with a tuple of cut heights, and a cell of $\Kop$ a chain together with a rule assigning windows.

\begin{remark}\label{rem:morphism-convention}
For definiteness: a map of $\Cal{P}$-algebras is a map $f$ of underlying symmetric sequences commuting strictly with the action maps, and an equivalence is such an $f$ that is a levelwise weak equivalence. Restriction along a map of $\Nl$-operads $g\colon \Cal{P}'\to \Cal{P}$ sends $\Cal{P}$-algebras to $\Cal{P}'$-algebras and preserves maps and equivalences; this is immediate from the definitions.
\end{remark}

\begin{definition}\label{def:structures-agree}
Let $\Cal{W}$ be a symmetric sequence with an $\Aop$-algebra structure $\lambda$, and $\Cal{W}'$ a symmetric sequence with an $\Oper$-algebra structure $\nu$ (equivalently, a strict reduced operad structure). We say that \emph{$(\Cal{W},\lambda)$ and $(\Cal{W}',\nu)$ agree through a cut refinement} if there exist:
\begin{enumerate}[label=(\roman*)]
    \item $\Nl$-operads $\Lop$, $\Kop$ with maps of $\Nl$-operads $\rho_\Lop\colon\Lop\to\Oper$, $\rho_\Kop\colon\Kop\to\Oper$, $j\colon \Lop\to \Aop$, $o\colon \Oper\to \Kop$ and $\sigma\colon\Lop\to\Kop$, such that $\rho_\Lop$, $\rho_\Kop$ and $j$ are levelwise weak equivalences, $\rho_\Kop\circ o = \id$, and $\sigma$ covers $\Oper$ ($\rho_\Kop\circ\sigma = \rho_\Lop$); it follows by two-out-of-three that $o$ and $\sigma$ are levelwise weak equivalences as well;
    \item a $\Kop$-algebra structure $\tau$ on $\Cal{W}'$ with $o^{*}(\tau) = \nu$;
    \item a zig-zag of levelwise weak equivalences of underlying symmetric sequences between $\Cal{W}'$ and $\Cal{W}$ strictly intertwining the restriction $\zeta := \sigma^{*}(\tau)$ with $j^{*}\lambda$.
\end{enumerate}
Both structures are then restrictions of the single coherent algebra $\tau$ along maps of $\Nl$-operads, over structuring operads all connected to $\Oper$, hence to each other and to $\Aop$, by levelwise weak equivalences; each stage exhibits an $A_\infty$-operad in the sense of \cite[Definition 6.38]{Clark-1}, and the notion refines the equivalence of $\Aop$-algebras used in \cite[\S 8.2]{Clark-1}.
\end{definition}

\subsection{The two structures on one object}\label{sec:twostructures}

By \cite[\S 8.1]{Clark-1} (proof of Proposition 5.9 there), for every $\boxcirc$-monoid $(X, m, u)$ in $\SymSeq^{\Delta}$ the restricted totalization $\Tot X$ carries a natural $\Aop$-algebra structure
\[ \lambda = \lambda(X,m,u), \]
where the action of a cell $\varphi\in \Aop_\ell$ is built from $\varphi$ by mapping-space functoriality over $\Dres$ together with the iterated multiplication of $(X,m)$. Applied to $X = C(\uS)$ (Proposition \ref{prop:box-monoid}) this is the structure of Theorem \ref{thm:main-1}, which we denote $\lambda^{\mathrm{lev}}$ (``levelled'').

\begin{lemma}\label{lem:lambda-natural}
A map of $\boxcirc$-monoids $f\colon (X,m,u)\to (X',m',u')$ (a map of cosimplicial symmetric sequences commuting strictly with multiplications and units) induces a strict map of $\Aop$-algebras $\Tot f\colon (\Tot X, \lambda)\to (\Tot X', \lambda')$.
\end{lemma}

\begin{proof}
Every constituent of the composite defining the action in \cite[\S 8.1]{Clark-1} is functorial in $X$: the mapping objects
\[\Map_{\Dres}(\SD^{\bullet}_+, X), \quad\quad \Map\big((\SD^{\bullet}_+)^{\boxcirc \ell},\, X^{\boxcirc \ell}\big),\]
and the final contraction, which uses only iterates of $m$.
\end{proof}

\begin{remark}\label{rem:clark1-interface}
Lemma \ref{lem:lambda-natural} and the recollections of \S\ref{sec:Nlev} are the full extent of our dependence on the internals of \cite{Clark-1}: we use (i) the existence and naturality of $\lambda$; (ii) that $\Aop$-cells act through the cosimplicial variable only, an $\Aop$-cell being a $\Dres$-natural family on the $\SD^{\bullet}_+$-side, so its action on $\Tot C(\uS)$ reparametrizes simplex coordinates uniformly in the chain variable and then applies iterated $m$; and (iii) the equivalence $\rho\colon \Aop\to \Oper$ \cite[Propositions 7.12, 7.16]{Clark-1}. Item (ii) is what makes Proposition \ref{prop:transport} a finite check.
\end{remark}

By \cite[Proposition 7.12]{Clark-1} Theorem \ref{thm:omw-operad} makes $\Omw(\uS)$ an $\Oper$-algebra, denoted $\nu^{\mathrm{tr}}$ (``tree''). It remains to construct $\Lop$, $\Kop$ and the data of Definition \ref{def:structures-agree}.

\subsection{Flags and their pieces}\label{sec:flags}

Fix $t\geq 1$ and a non-basepoint flag $\gamma = (\mathsf{min} = \lambda_0\leq \lambda_1\leq \cdots \leq \lambda_\ell = \mathsf{max})\in P(t)^{\circ}_{\ell}$. A \emph{position} of $\gamma$ is a pair $\rho = (j, B)$ with $0\leq j\leq \ell-1$ and $B$ a block of $\lambda_j$; its \emph{children set} is $C_\rho := \{C\in \lambda_{j+1}: C\subseteq B\}$, and $\rho$ is \emph{unary} if $|C_\rho| = 1$. The \emph{target} of $\gamma$ is
\[ \BC(\gamma) := \bigwedge_{\rho} \BC(C_\rho),\]
the smash over all positions (unary positions contributing $S^0$), matching the $\gamma$-summand indexing of the iterated composition product (Definition \ref{comp-prod}).

Let $(\Cal{V},h)$ be a weighted tree on $\mathbf{t}$ containing every non-singleton block of every $\lambda_j$ (``the flag is present''). Iterating Lemma \ref{lem:graft-degraft} decomposes $\Cal{V}$ uniquely into \emph{pieces}: for each position $\rho = (j,B)$, the tree on $C_\rho$
\[ \Cal{V}_\rho := q_\rho\big( \{ D\in \Cal{V} : D\subseteq B \text{ and } D\supseteq \text{some } C\in C_\rho \} \big),\]
where $q_\rho$ collapses each $C\in C_\rho$ to a point. Every vertex of $\Cal{V}$ is an internal vertex of exactly one piece, except that each flag block $C\in\lambda_{j+1}$ appears twice, as a leaf of its layer-$j$ piece and the root of its layer-$(j{+}1)$ piece; heights are carried along unchanged.

\subsection{Framed cut systems}\label{sec:FS}

\begin{definition}\label{def:FS}
A \emph{framed cut system} of shape $\gamma$ is a family $F = (S_\rho, E_\rho)_{\rho \text{ non-unary}}$ of continuous functions, defined on the closed subspace of $\coprod_{\Cal{V}}\wt(\Cal{V})$ where the flag is present, with values in $[0,1]$, subject to:
\begin{enumerate}[label=(F\arabic*)]
    \item \emph{(descent, equivariance)} each $S_\rho, E_\rho$ is invariant under the collapse moves of Definition \ref{def:B}(a) that preserve the flag, and the family is $\Sigma_t$-equivariant in the evident sense ($S_{\sigma\rho}\circ \sigma = S_{\rho}$, etc.);
    \item \emph{(order)} $S_\rho \leq E_\rho$ pointwise;
    \item \emph{(floor anchoring)} pointwise,
    \[ S_{(j,B)} \geq h\big(p_{\Cal{V}}(B)\big), \quad\quad \text{with the convention } h(p_{\Cal{V}}(\mathbf{t})) := 0.\]
\end{enumerate}
We write $\mathrm{FS}(\gamma)$ for the set of framed cut systems of shape $\gamma$, topologized by the supremum metric on each (compact) stratum.

The \emph{cut map} $\kappa_F\colon \BC(t)\to \BC(\gamma)$ is defined on representatives by the rescaled heights
\begin{equation}\label{eq:FS-formula}
    h_\rho\big(q_\rho(D)\big) := \frac{h(D) - S_\rho(\Cal{V},h)}{E_\rho(\Cal{V},h) - S_\rho(\Cal{V},h)} \quad\quad (D \text{ an internal vertex of } \Cal{V}_\rho),
\end{equation}
setting $\kappa_F(\Cal{V},h) := \big( (\Cal{V}_\rho, h_\rho) \big)_{\rho}$, and $\kappa_F(\Cal{V},h) := \basept$ if the flag is absent from $\Cal{V}$ or if at some non-unary $\rho$ the \emph{thresholds} trip: $S_\rho = E_\rho$, or \eqref{eq:FS-formula} fails to be a weighting of $\Cal{V}_\rho$ with positive root height and all leaf-parent heights $<1$.
\end{definition}

Monotonicity of \eqref{eq:FS-formula} is automatic, so the only failures are the listed thresholds. Axiom (F3) is what makes $\kappa_F$ well defined at the wall across which the flag ceases to be present: a collapse of a flag block onto its parent is a move of type (a) into the flag-absent locus, where the value is the basepoint by definition, and (F3) forces the corresponding piece's root height $\leq 0$ on the other representative, so the two values agree. No ceiling axiom is needed for this, since overflow above $1$ is self-absorbing: the maximum height in a piece is attained at a leaf-parent, by the observation after Definition \ref{def:reduced}, so overflow forces a vanishing leaf edge. Continuity is a separate matter; see Remark \ref{rem:FS-not-cts}.

\begin{proposition}\label{prop:FS-continuous}
For every $F\in \mathrm{FS}(\gamma)$, the map $\kappa_F$ is well defined on $\BC(t)$ and based, and the family $\{\kappa_{\sigma F}\}$ is $\Sigma_t$-equivariant. For the tree frames and the labelled frames of Example \ref{ex:frames} below, $\kappa_F$ is moreover continuous.
\end{proposition}

\begin{proof}
\emph{Well definedness.} Collapse moves preserving the flag transport heights by the affine maps \eqref{eq:FS-formula} and leave frames invariant by (F1), a collapsed non-flag vertex lying in exactly one piece together with its parent; a collapse of a flag block is absorbed by (F3) as explained above. For Definition \ref{def:B}(b): if $h(\mathbf{t})=0$ the root position's rescaled root height is $\leq 0$ since $S\geq 0$; if $h(v(x))=1$ the leaf-parent value in the piece containing $v(x)$ is $\geq 1$ since $E\leq 1$. Both give the basepoint. Based-ness and equivariance are immediate from the formulas and (F1).

\emph{Continuity of the two families.} By Example \ref{ex:Wcells}(1) and Lemma \ref{lem:equal-evals} below, $\kappa_{F^{\mathrm{tr}}} = \kappa_{W^{\mathrm{tr}}}$ and $\kappa_{F^{e}} = \kappa_{W^{e}}$ are cut maps of windowed cut systems, which are continuous by Proposition \ref{prop:K-continuity}. (Neither result depends on the present proposition.)
\end{proof}

\begin{remark}\label{rem:FS-not-cts}
For a general $F\in\mathrm{FS}(\gamma)$ the map $\kappa_F$ need not be continuous: a window $[S_\rho,E_\rho]$ may shrink to a point at an interior point of $\BC(t)$ while the heights of its piece sit inside it at bounded ratios, and then the value stays away from the basepoint although the limit is declared the basepoint. For $t=2$, with $h$ the height of the internal vertex, $S := h-g$, $E := h+g$, $g := \min\{|h-\tfrac12|, h, 1-h\}$ satisfies (F1)--(F3), and $\kappa_F$ is the constant $\tfrac12\in\BC(2)$ off $h=\tfrac12$ and the basepoint there. For the tree frames the window $[h(p(B)),1]$ shrinks only when $h(p(B))\to 1$, and for the labelled frames $[\max(h(p(B)),e_{\pi(\nu)}), e_\nu]$ only when $h(p(B))\to e_\nu$; the piece at $\nu$ need not degenerate then, but $p(B)$ is a leaf-parent of the piece at $\anc(\nu)$ (Lemma \ref{lem:anc-visible}), whose ceiling is at most the limit of $h(p(B))$, so that factor does. Axioms (W5) and (W6) of Definition \ref{def:WCS} isolate this mechanism.
\end{remark}

\begin{example}\label{ex:frames}
Three families of framed cut systems.
\begin{enumerate}[label=(\arabic*)]
    \item \emph{Tree frames} $F^{\mathrm{tr}}$: $S_{(j,B)} := h(p_{\Cal{V}}(B))$, $E_{(j,B)} := 1$, satisfying (F1)--(F3). For $\ell = 2$ the formulas agree factor by factor with Definition \ref{def:cut}, so $\kappa_{F^{\mathrm{tr}}} = c^{\mathrm{tr}}_{\alpha}$; for general $\ell$ the corresponding windowed cells are closed under composition (Proposition \ref{prop:K-operad}) and compute Ching's iterated cocompositions.
    \item \emph{Labelled frames} $F^{e}$, for a cut labelling $e$ (Definition \ref{def:L} below):
    \[ S_\nu := \max\big(h(p_{\Cal{V}}(B)),\, e_{\pi(\nu)}\big), \quad\quad E_\nu := e_\nu.\]
    We write $F^{\vec a}$ for the \emph{uniform} case $e_{(j,B)} = a_{j+1}$, where $0 < a_1 < \cdots < a_{\ell-1} < 1$. (F3) holds by the $\max$, and (F2) fails only where $S\geq E$, on which locus the order threshold already returns the basepoint. Unwinding the thresholds, $\kappa_{F^{\vec a}}(\Cal{V},h)\neq \basept$ precisely when $\alive{(\Cal{V},h)}{a_{j}} = \lambda_j$ for all $j$ (Definition \ref{def:alive}), and then each piece is rescaled by $[a_j, a_{j+1}]\to[0,1]$, the $\max$ in $S$ equalling $a_j$. For $\ell = 2$, $a_1 = \tfrac12$ we write $c^{\mathrm{lev}}_{\alpha} := \kappa_{F^{(1/2)}}$, the \emph{levelled cut} at $\alpha$, computed in coordinates in Appendix \ref{sec:anchor-low}.
    \item \emph{Straight lines}: (F2) and (F3) are pointwise inequalities against the fixed continuous function $h(p_{\Cal{V}}(B))$ and (F1) is linear, so all survive convex combination, and $\mathrm{FS}(\gamma)$ is convex. This is of limited use by itself, since a convex combination of framed systems with continuous cut maps need not have a continuous cut map (Remark \ref{rem:FS-not-cts}); the contractible cell spaces the comparison needs are those of \S\ref{sec:K-and-L}. Taking $F = F^{\mathrm{tr}}$ and $F' = F^{\vec a}$ with $\ell = 2$, $a = \tfrac12$ gives, up to reparametrization, the two-leg homotopy of Appendix \ref{sec:anchor-low}, whose continuity is checked there directly; the naive interpolation that fails there is exactly the one leaving $\mathrm{FS}$ by violating (F3).
\end{enumerate}
\end{example}

\subsection{The cut-system operad and the tautological algebra}\label{sec:K-and-L}

For $\theta\in \Oper_\ell(\underline{p}; t)$ write $\gamma(\theta)\in P(t)_\ell$ for the underlying flag: by \cite[Definition 7.2, Remark 7.3]{Clark-1} the cells of $\Oper$ in level $\ell$ are indexed by flags together with block-labelling data, and the labelling enters below only through equivariance, so we carry it silently.

Both families of Example \ref{ex:frames} are now exhibited as cells of one $\Nl$-operad. The composition law is determined by the geometry: cutting along a merged flag is cutting along the outer flag and then cutting each piece along the inner flags, with inner windows read inside outer ones. Two conventions make this everywhere defined and strictly associative: heights are arbitrary reals, so pieces carried with their \emph{ambient} heights never leave the domain of the frames; and every frame is a function of the window it is handed, so no normalization, hence no division, occurs.

\begin{definition}\label{def:windows}
Fix $t\geq 1$ and a flag $\gamma\in P(t)^{\circ}_{\ell}$. Write
\begin{align*}
 \Dg{\gamma} := \big\{ (\Cal{V}, h) \;:\;\ & \Cal{V} \text{ a tree on } \mathbf{t} \text{ in which } \gamma \text{ is present},\\
 & h\colon V(\Cal{V})\to \mathbb{R} \text{ order-preserving} \big\},
\end{align*}
topologized as the disjoint union over trees of the evident closed subspaces of $\mathbb{R}^{V(\Cal{V})}$, and
\[ \widehat{W} := \{(u,v)\in \mathbb{R}^{2} : u\leq v\}.\]
On $\Dg{\gamma}\times \widehat{W}$ we adopt once and for all the convention $h(p_{\Cal{V}}(\mathbf{t})) := u$: the root's parent height reads the window start. For a non-unary position $\nu = (j,B)$ of $\gamma$, the \emph{run} of $\nu$ is the maximal interval $[j_0, j]$ of layers with $B\in \lambda_{j'}$ for all $j'\in[j_0,j]$. If $j_0 = 0$ (equivalently $B = \mathbf{t}$) we call $\nu$ a \emph{root-run} position; otherwise the \emph{anchor position} of $\nu$ is
\[ \anc(\nu) := (j_0 - 1,\, B^{*}), \quad\quad B^{*} := \text{the block of } \lambda_{j_0-1} \text{ containing } B.\]
By minimality of $j_0$ we have $B^{*}\supsetneq B$, so $\anc(\nu)$ is non-unary ($B^{*}$ has at least two children in $\lambda_{j_0}$). The anchor position depends only on the flag, and equivariantly so.
\end{definition}

\begin{lemma}\label{lem:anc-visible}
Let $(\Cal{V},h)$ carry the flag $\gamma$ and let $\nu = (j,B)$ be a non-unary, non-root-run position. Then $p_{\Cal{V}}(B)$ determines a vertex of the piece at $\anc(\nu)$, of height $h(p_{\Cal{V}}(B))$, an internal vertex which may be the piece's root.
\end{lemma}

\begin{proof}
Write $\anc(\nu) = (j_0-1, B^{*})$. Since $B^{*}\supsetneq B$ is a non-singleton flag block, minimality of the parent gives $B\subseteq p_{\Cal{V}}(B)\subseteq B^{*}$, and $B$ is a child of $\anc(\nu)$, so $p_{\Cal{V}}(B)$ lies in the defining set of the $\anc(\nu)$-piece. Its image under $q$ contains the point of $B$ and, since $p_{\Cal{V}}(B)\smallsetminus B$ is covered by the other children, at least one more; so it is non-singleton in the piece tree.
\end{proof}

\begin{definition}\label{def:WCS}
A \emph{windowed cut system} of shape $\gamma$ is a family
\[ W = (S_\nu, E_\nu)_{\nu\,\text{non-unary}} \]
of continuous real-valued functions on $\Dg{\gamma}\times\widehat{W}$, subject to the following pointwise axioms:
\begin{enumerate}[label=(W\arabic*)]
    \item \emph{(descent, equivariance)} each $S_\nu, E_\nu$ is invariant under the collapse moves of Definition \ref{def:B}(a) that preserve the flag (applied to real-valued weightings, where they make the same sense), and the family is $\Sigma_t$-equivariant;
    \item \emph{(subordination)} $u\leq S_\nu \leq E_\nu\leq v$;
    \item \emph{(anchoring or stacking)} $S_\nu \geq h\big(p_{\Cal{V}}(B)\big)$, \ or \ $E_{\anc(\nu)} \leq S_\nu$; at a root-run position the first disjunct reads $S_\nu\geq u$ (the convention of Definition \ref{def:windows}) and is an instance of (W2);
    \item \emph{(affine equivariance)} $W_\nu\big((\Cal{V}, a h + b);\, (a u + b,\, a v + b)\big) = a\cdot W_\nu\big((\Cal{V},h);\,(u,v)\big) + b$ for all $a > 0$, $b\in\mathbb{R}$;
    \item \emph{(windows shrink only towards the anchor)} writing $w_\nu := E_\nu - S_\nu$ for the \emph{window length}, there is a constant $\lambda_\nu > 0$ such that
    \[ w_\nu \;\geq\; \lambda_\nu\cdot \min\big\{\, v-u,\ \ E_\nu - h(p_{\Cal{V}}(B)),\ \ w_{\anc(\nu)} \,\big\} \]
    if $\nu$ is not root-run, and $w_\nu \geq \lambda_\nu\,(v-u)$ if it is;
    \item \emph{(anchor ceilings)} there is a constant $K_\nu \geq 0$ such that $E_{\anc(\nu)} - E_\nu \leq K_\nu\, w_\nu$ for every non-root-run $\nu$.
\end{enumerate}
Only the existence of the constants $\lambda_\nu$, $K_\nu$ is required. (W5) says that the window at $\nu$ can degenerate only if the parent of $B$ has reached the ceiling of $\nu$ or the window at the anchor has degenerated; (W6) that the anchor's ceiling exceeds the ceiling of $\nu$ by at most a bounded multiple of the window at $\nu$.
We write $\mathrm{WS}(\gamma)$ for the space of windowed cut systems of shape $\gamma$, with the compact-open topology. The \emph{cut map} of $W\in \mathrm{WS}(\gamma)$ is the cut map, in the sense of Definition \ref{def:FS}, of the top-window slice:
\[ \kappa_W := \kappa_{\,W(-\,;\,(0,1))\,}\colon \BC(t)\longrightarrow \BC(\gamma),\]
using that by (W2) the slice $W(-;(0,1))$, restricted to genuine weightings, is a family of continuous functions with values in $[0,1]$ satisfying (F1) and (F2), to which the formula and conventions of Definition \ref{def:FS} apply verbatim. The axiom (F3) is not assumed; continuity of $\kappa_W$, jointly in $W$, is Proposition \ref{prop:K-continuity}, and it is (W5) and (W6) that provide it (the system of Remark \ref{rem:FS-not-cts}, read at the window $(u,v)$ with $g := \min\{|h-\tfrac{u+v}{2}|, h-u, v-h\}$, satisfies (W1)--(W4)).
\end{definition}

\begin{example}\label{ex:Wcells}
The two families of Example \ref{ex:frames}, as windowed systems. Both satisfy (W1) and (W4); (W2), (W3), (W5) and (W6) are checked below.
\begin{enumerate}[label=(\arabic*)]
    \item \emph{Tree cells} $W^{\mathrm{tr}}$, where $\operatorname{clamp}_{[u,v]}(x) := \min(\max(x, u), v)$:
    \[ S_\nu := \operatorname{clamp}_{[u,v]}\big(h(p_{\Cal{V}}(B))\big), \quad\quad E_\nu := v.\]
    On genuine weightings at the top window the clamp is inactive, so the slice is the tree frame $F^{\mathrm{tr}}$ of Example \ref{ex:frames}(1). (W2) is immediate. (W3): where $h(p_{\Cal{V}}(B))\leq v$ the clamp only raises the value, giving anchoring; where $h(p_{\Cal{V}}(B)) > v$ we get $S_\nu = v = E_{\anc(\nu)}$, all tree ceilings being $v$, giving stacking with equality. The clamp is needed for (W2) at heights outside $[u,v]$. (W5): $w_\nu = v - \operatorname{clamp}_{[u,v]}(h(p_{\Cal{V}}(B))) \geq \min\{v-u,\ v - h(p_{\Cal{V}}(B))\}$, so $\lambda_\nu = 1$ works; (W6): all ceilings equal $v$, so $K_\nu = 0$ works.
    \item \emph{Constant cells}: for a cut labelling $e$ for $\theta$ (Definition \ref{def:L} below),
    \[ W^{e}\colon \quad S_\nu := u + (v-u)\, e_{\pi(\nu)}, \quad\quad E_\nu := u + (v-u)\, e_{\nu},\]
    where $\pi(\nu)$ denotes the parent position of $\nu$ and the conventions of Definition \ref{def:L} are in force ($e = 1$ on the top layer, $e = 0$ at the parent of the root position). The top-window slice is the \emph{stacked constant} family $(e_{\pi(\nu)}, e_\nu)$: the labelled frames $F^{e}$ of Example \ref{ex:frames}(2) with the $\max$ deleted. (W2): $S_\nu, E_\nu$ are convex combinations of $u$ and $v$, in the order given since $e_{\pi(\nu)} < e_{\nu}$. (W3), second disjunct, always: $E_{\anc(\nu)} = u + (v-u)e_{\anc(\nu)}$ and $e_{\anc(\nu)}\leq e_{\pi(\nu)}$ since the labels strictly increase along the run from $\anc(\nu)$ to $\pi(\nu)$ (with equality of positions when the run has length one); so $E_{\anc(\nu)}\leq S_\nu$. (W5): $w_\nu = (v-u)(e_\nu - e_{\pi(\nu)})$, so $\lambda_\nu := e_\nu - e_{\pi(\nu)} > 0$ works; (W6): $E_{\anc(\nu)} - E_\nu = (v-u)(e_{\anc(\nu)} - e_\nu) \leq 0$, so $K_\nu = 0$ works.
\end{enumerate}
\end{example}

\begin{definition}\label{def:K}
The \emph{cut-system operad} $\Kop$ is the $\Nl$-operad with spaces
\[ \Kop_\ell(\underline{p}; t) := \coprod_{\theta\in \Oper_\ell(\underline{p};t)} \mathrm{WS}\big(\gamma(\theta)\big),\]
$\Sigma$-structure by relabelling, unit the level-$1$ tree cell, and composition by \emph{substitution at raw pieces} in the format of \cite[\S 6]{Clark-1}: the composite of $(\theta, W)$ with $(\theta'_\rho, W'_\rho)$ at level slot $j$ keeps the frames of $W$ at unrefined positions and, at the insert position at branch $\rho$ corresponding to a non-unary $\nu'$ of $\gamma(\theta'_\rho)$, sets
\begin{equation}\label{eq:substitution}
    W''_{\nu'}\big((\Cal{V},h);(u,v)\big) := W'_{\nu'}\Big( \big(\Cal{V}_\rho,\, h|^{\mathrm{raw}}\big);\ \big(S_\rho((\Cal{V},h);(u,v)),\, E_\rho((\Cal{V},h);(u,v))\big)\Big),
\end{equation}
where the \emph{raw piece} $(\Cal{V}_\rho, h|^{\mathrm{raw}})$ is the piece of \S\ref{sec:flags} with the ambient heights carried along \emph{un-rescaled}, and its root-parent convention reads the handed window start $S_\rho$. No rescaling and no division occur: \eqref{eq:substitution} is defined at every point of $\Dg{\text{merged flag}}\times\widehat{W}$.
\end{definition}

\begin{proposition}\label{prop:K-operad}
$\Kop$ is a well-defined $\Nl$-operad. The projection $\rho_{\Kop}\colon (\theta, W)\mapsto\theta$ is a map of $\Nl$-operads and a levelwise weak equivalence; the tree cells define a map of $\Nl$-operads
\[ o\colon \Oper\longrightarrow \Kop, \quad\quad \theta\mapsto (\theta, W^{\mathrm{tr}}),\]
with $\rho_{\Kop}\circ o = \id$; and each space $\mathrm{WS}(\gamma)$ deformation-retracts onto the tree cell, so $o$ is a levelwise weak equivalence as well.
\end{proposition}

\begin{proof}
\emph{Substitution.} On the merged-flag-present locus every block of the inner flag at $\rho$ is a vertex of $\Cal{V}$, so the raw piece lies in $\Dg{\gamma(\theta'_\rho)}$, and the handed window lies in $\widehat{W}$ by the outer (W2); thus \eqref{eq:substitution} is a composite of continuous total maps, jointly continuous in the cells.

\emph{Closure.} (W1) and (W4) hold since pieces, collapse moves and relabellings commute and the raw piece map commutes with $h\mapsto ah+b$; (W2) is the nesting $[S'',E'']\subseteq[S_\rho,E_\rho]\subseteq[u,v]$. For (W3), write $B$ for the block of the refined position $\rho=(j,B)$ and run over the positions of the merged flag. At positions of the insert lying in the run of $B$ (the seam, and the inner root-run positions) the anchor is $\anc(\rho)$, computed in the outer flag with untouched frame, and $S''\geq S_\rho$ by the inner (W2); so whichever disjunct the outer (W3) supplied at $\rho$ passes to the composite. At positions of the insert outside that run, the block $\bar D$, its tree-parent and its anchor all lie inside the insert, and the raw piece carries ambient heights, so the inner (W3) \emph{is} the required instance verbatim (parenthood matches across the piece, since $p_{\Cal{V}}(\bar D)\subseteq B$ and the collapse $q$ preserves it). At an unrefined position $\nu$ whose run enters the refined layer, the new anchor is an insert position whose ceiling is an insert frame at the branch $(j,B')$ of the old anchor, hence at most the old ceiling $E_{B'}$ by the inner (W2): stacking is inherited, and anchoring is untouched, since \emph{anchor ceilings only decrease under composition}. All remaining positions are untouched.

\emph{Closure of (W5) and (W6).} Write $\lambda_\nu, K_\nu$ for constants of the outer cell, $\lambda'_{\nu'}, K'_{\nu'}$ for those of the inner cell at $\rho$, and $w''$ for the windows of the composite. At the seam and at the inner root-run positions the window is an inner root-run window at the handed window $(S_\rho,E_\rho)$, so $w''\geq \lambda'_{\nu'}\, w_\rho$ by the inner (W5) at the inner root-run position $\nu'$ in question, while the anchor is $\anc(\rho)$. Since $E''\leq E_\rho$,
\[ \min\{v-u,\ E''-h(p_{\Cal V}(B)),\ w_{\anc(\rho)}\} \;\leq\; \min\{v-u,\ E_\rho-h(p_{\Cal V}(B)),\ w_{\anc(\rho)}\} \;\leq\; w_\rho/\lambda_\rho \;\leq\; w''/(\lambda_\rho\lambda'_{\nu'}),\]
which is (W5) with the constant $\lambda_\rho\lambda'_{\nu'}$ (if $\rho$ is root-run, $w''\geq \lambda'_{\nu'}\lambda_\rho(v-u)$ directly), and
\[ E_{\anc(\rho)} - E'' \;\leq\; (E_{\anc(\rho)} - E_\rho) + (E_\rho - S_\rho) \;\leq\; (K_\rho + 1)\, w_\rho \;\leq\; \tfrac{K_\rho+1}{\lambda'_{\nu'}}\, w''\]
is (W6). At the insert positions off the run of $B$, block, tree-parent and anchor lie inside the insert and the raw piece carries ambient heights, so the inner (W5) and (W6) at the handed window are the required instances: the inner (W5) has the handed length $w_\rho$ in place of $v-u$, but the anchor's window lies inside the handed window, so $\min\{w_\rho, \cdot, w''_{\anc}\} = \min\{v-u,\cdot,w''_{\anc}\}$. At an unrefined position whose anchor becomes an insert position, the new anchor's window is contained in the old one and its ceiling is at most the old one (inner (W2)), so both axioms are inherited with the same constants. All other positions are untouched.

\emph{Associativity, units, equivariance.} Associativity is the transitivity of degrafting (Lemma \ref{lem:graft-degraft} iterated, heights carried along unchanged) together with associativity of substitution of functions: in either bracketing the innermost cell is handed the middle cell's frame evaluated at the outer window. The level-$1$ tree cell has frame $(\operatorname{clamp}_{[u,v]}(u), v) = (u,v)$, so substituting it changes nothing. Equivariance is (W1), and the composite covers $\Oper$'s composition by construction.

\emph{The section $o$.} At an insert position with block $\bar D$ the handed window is $(S_\rho, v)$ with $S_\rho = \operatorname{clamp}_{[u,v]}(h(p_{\Cal{V}}(B)))$, so \eqref{eq:substitution} gives $E''=v$ and
\[ S'' = \operatorname{clamp}_{[S_\rho,\, v]}\big( h(p_{\Cal{V}}(\bar D)) \big) = \operatorname{clamp}_{[u,v]}\big( h(p_{\Cal{V}}(\bar D)) \big),\]
the second equality by transitivity of the clamp along $h(p_{\Cal{V}}(\bar D))\geq h(p_{\Cal{V}}(B))$ (valid since $p_{\Cal{V}}(\bar D)\subseteq B\subsetneq p_{\Cal{V}}(B)$, with equality of positions at the seam). This is the merged formula, so $o$ respects composition; units are the definition.

\emph{Contractibility of $\mathrm{WS}(\gamma)$.} Three legs, each jointly continuous in the cell and the parameter and each fixing $W^{\mathrm{tr}}$. Leg 1 is the straight line to $W^{\flat}$, with ceilings unchanged and
\[S^{\flat}_\nu := \min\big(\max(S_\nu,\, h(p_{\Cal{V}}(B))),\, E_\nu\big);\]
floors only rise, so (W2) and both disjuncts of (W3) persist, the second because where $h(p_{\Cal{V}}(B)) > E_\nu$ anchoring is unavailable to any system by (W2), while $S^{\flat}_\nu = E_\nu\geq S_\nu\geq E_{\anc(\nu)}$. Leg 2 raises the ceilings to $v$, re-anchoring as it goes:
\[ E_{\nu,s} := (1-s)E_\nu + s\,v, \quad\quad S_{\nu,s} := \min\big( \max(S_\nu,\, h(p_{\Cal{V}}(B))),\, E_{\nu,s}\big).\]
Anchoring holds wherever $h(p_{\Cal{V}}(B))\leq E_{\nu,s}$; elsewhere the system was stacked and stays so, since $E_{\anc(\nu),s}\leq E_{\nu,s} = S_{\nu,s}$ by monotonicity of the convex combination in $E_{\anc(\nu)}\leq E_\nu$. Leg 3 is the straight line from there to $W^{\mathrm{tr}}$: ceilings are constantly $v$, and the floors either stay in $[h(p_{\Cal{V}}(B)),v]$, giving anchoring, or are constantly $v$ where $h(p_{\Cal{V}}(B))>v$, giving stacking with equality. Throughout, (W1) and (W4) survive because clamps and convex combinations commute with increasing affine maps. So do (W5) and (W6), with constants that change along the path, which is all that contractibility requires. Along Legs 1 and 2 the window is $w_{\nu,s} = \min\{E_{\nu,s} - S_\nu,\ E_{\nu,s} - h(p_{\Cal V}(B))\}^{+}$ (Leg 1 is the case $s=0$). Where $h(p_{\Cal V}(B))\geq S_\nu$ this is $(E_{\nu,s} - h(p_{\Cal V}(B)))^{+}$, so (W5) holds with the constant $\min\{\lambda_\nu, 1\}$, and $w_{\nu,s}\geq E_\nu - S_\nu = w_\nu$ unless $w_{\nu,s} = 0$, in which case $h(p_{\Cal V}(B)) = S_\nu = E_\nu = E_{\nu,s}$ and (W6) for $W$ gives $E_{\anc(\nu)}\leq E_\nu$; either way $E_{\anc(\nu),s} - E_{\nu,s} = (1-s)(E_{\anc(\nu)} - E_\nu)\leq K_\nu w_{\nu,s}$. Where $h(p_{\Cal V}(B)) < S_\nu$ the cell is stacked at $\nu$, so $E_{\anc(\nu)}\leq S_\nu\leq E_\nu$, and $w_{\nu,s} = w_\nu + s(v-E_\nu)$. If the minimum in (W5) for $W$ is attained at $v-u$ or at $E_\nu - h(p_{\Cal V}(B))$, then $w_{\nu,s}\geq w_\nu\geq\lambda_\nu(v-u)$, respectively $w_{\nu,s}\geq \lambda_\nu(E_\nu - h(p_{\Cal V}(B))) + s(v-E_\nu)\geq \lambda_\nu(E_{\nu,s} - h(p_{\Cal V}(B)))$. If it is attained at $w_{\anc(\nu)}$, then, floors only rising along the legs,
\[ w_{\anc(\nu),s} \;\leq\; w_{\anc(\nu)} + s\,(v - E_{\anc(\nu)}) \;\leq\; w_\nu/\lambda_\nu + s(v-E_\nu) + K_\nu w_\nu \;\leq\; \max\{1/\lambda_\nu + K_\nu,\ 1\}\; w_{\nu,s}\]
by (W5) and (W6) for $W$. So (W5) holds with the constant $\min\{\lambda_\nu, (1/\lambda_\nu + K_\nu)^{-1}, 1\}$, and (W6) with $K_\nu$, since $E_{\anc(\nu),s} - E_{\nu,s} = (1-s)(E_{\anc(\nu)} - E_\nu) \leq 0$ here. Along Leg 3 all ceilings equal $v$, so (W6) is trivial, and $w_{\nu,t}\geq t\,(v - \operatorname{clamp}_{[u,v]}h(p_{\Cal V}(B)))\geq t\min\{v-u,\ v-h(p_{\Cal V}(B))\}$ gives (W5) with the constant $t$ for $t>0$, and with the constant of the Leg-2 endpoint at $t=0$. Levelwise $\rho_{\Kop}$ is therefore a disjoint union, over the discrete set of $\theta$'s, of projections with contractible fibers, hence a weak equivalence, and $o$ is a section of it.
\end{proof}

\begin{proposition}\label{prop:K-continuity}
For every $W\in\mathrm{WS}(\gamma)$, the cut map $\kappa_W$ is well defined on $\BC(t)$ and based; the assignment $(W,x)\mapsto\kappa_W(x)$ is continuous on $\mathrm{WS}(\gamma)\times\BC(t)$, in particular each $\kappa_W$ is continuous; and the family $\{\kappa_{\sigma W}\}$ is $\Sigma_t$-equivariant.
\end{proposition}

\begin{proof}
\emph{Well-definedness and based-ness.} Collapse moves that preserve the flag transport heights by the affine maps \eqref{eq:FS-formula} and leave the frames invariant by (W1), a collapsed non-flag vertex lying in exactly one piece together with its parent, and Definition \ref{def:B}(b) is absorbed as in the proof of Proposition \ref{prop:FS-continuous}, using $S_\nu\geq u = 0$ and $E_\nu\leq v = 1$. A collapse of a non-singleton flag block $T$, that is $h(T) = h(p_{\Cal V}(T))$, leaves the flag-present locus, where the value is the basepoint by definition, so we must check that the value on the representative in which $T$ is still present is the basepoint too. Let $\nu_T$ be the last position of $T$'s run; it is non-unary. If $\nu_T$ is root-run, then $S_{\nu_T}\geq u = h(p_{\Cal V}(T))$ by (W2) and the convention of Definition \ref{def:windows}. Otherwise (W3) gives $S_{\nu_T}\geq h(p_{\Cal V}(T))$, or $E_{\anc(\nu_T)}\leq S_{\nu_T}$, and in the latter case either the $\anc(\nu_T)$-factor is already the basepoint or, by Lemma \ref{lem:anc-visible}, $h(p_{\Cal V}(T)) < E_{\anc(\nu_T)}\leq S_{\nu_T}$. In all cases the root height $(h(T) - S_{\nu_T})/w_{\nu_T}$ of the $\nu_T$-piece is $\leq 0$, or the window is degenerate, and the value is the basepoint. Equivariance is (W1).

\emph{Joint continuity.} Both spaces are metrizable ($\BC(t)$ is compact Hausdorff and second countable, and $\mathrm{WS}(\gamma)$ is a subspace of a function space with the compact-open topology on a $\sigma$-compact domain), so it suffices to show: if $W_n\to W$ in $\mathrm{WS}(\gamma)$ and $x_n\to x$ in $\BC(t)$, then $\kappa_{W_n}(x_n)\to\kappa_W(x)$. Convergence in the compact-open topology gives uniform convergence of every frame on the compact locus of genuine weightings at the top window, which is all that is used, and only the constants $\lambda_\nu$, $K_\nu$ of the \emph{limit} cell $W$ enter. We choose representatives in which the flag is present whenever possible, and write $S^n_\nu, E^n_\nu, w^n_\nu$ for the frames of $W_n$ at $x_n$ and $S_\nu, E_\nu, w_\nu$ for those of $W$ at $x$. Passing to subsequences, we may assume that each factor of $\kappa_{W_n}(x_n)$ is either the basepoint for all $n$ or valid for all $n$; if some factor is the basepoint for all $n$, the smash is, and there is nothing to prove, so assume every factor is valid for all $n$.

\emph{No degenerate window at $x$.} If $w_\nu > 0$ for every non-unary $\nu$, and $x$ lies in the flag-present locus with all $h(T) > h(p_{\Cal V}(T))$, then near $(W,x)$ every window is non-degenerate and the rescaled heights \eqref{eq:FS-formula} depend continuously on $(W_n,x_n)$; the thresholds of Definition \ref{def:FS} (a root height $\leq 0$, a leaf-parent height $\geq 1$) are absorbed by the basepoint of the target factor, and the compactness principle of \S\ref{sec:cotensors} gives convergence of the smash. If instead a flag block $T$ has collapsed at $x$, let $\nu_T$ be the last position of its run. By (W3) for $W_n$, either $S^n_{\nu_T}\geq h(p_{\Cal V}(T))$, so that the root height of the $\nu_T$-piece is at most $(h(T) - h(p_{\Cal V}(T)))/w^n_{\nu_T}$, which tends to $0$ when $w^n_{\nu_T}$ is bounded below (the case $w_{\nu_T} = 0$ is treated next); or $E^n_{\anc(\nu_T)}\leq S^n_{\nu_T}$, and then, the $\anc(\nu_T)$-factor being valid, Lemma \ref{lem:anc-visible} gives $h(p_{\Cal V}(T)) < E^n_{\anc(\nu_T)}\leq S^n_{\nu_T}$ at $x_n$, so the $\nu_T$-factor is the basepoint, contrary to assumption. Likewise, if $x$ is the basepoint of $\BC(t)$ because a root height tends to $0$ or a leaf-parent height tends to $1$ along $x_n$, then the root-run factor containing the root, whose window is bounded below by (W5), or the factor containing that leaf-parent, whose ceiling is $\leq 1$, degenerates, unless its window at $x$ is degenerate, which is again the case treated next.

\emph{A degenerate window at $x$.} Suppose $w_\nu = 0$ at $x$ for some non-unary $\nu = (j,B)$, and put $c := S_\nu = E_\nu$. The $\nu$-factor of $W_n$ at $x_n$ being valid, every height of the $\nu$-piece of $x_n$ lies in $(S^n_\nu, E^n_\nu)$, an interval of length $w^n_\nu\to 0$ about $c$; so all these heights tend to $c$, and in particular $h(B) = c$ at $x$. We argue along the anchor chain $\nu_0 := \nu$, $\nu_{i+1} := \anc(\nu_i)$, with blocks $B_0 := B\subsetneq B_1\subsetneq \cdots$, maintaining the statement
\begin{center}
$(\ast_i)$: \emph{the $\nu_i$-window of $W$ at $x$ is degenerate, at a position $c_i$, and all heights of the $\nu_i$-piece of $x_n$ tend to $c_i$},
\end{center}
which holds for $i=0$. Assume $(\ast_i)$. By (W5) at $\nu_i$ for $W$ at $x$, since $v-u = 1 > 0$, the position $\nu_i$ is not root-run, and either $h(p_{\Cal V}(B_i))\geq c_i$ at $x$ or $w_{\nu_{i+1}} = 0$ at $x$; by (W6), $E_{\nu_{i+1}}\leq c_i$ at $x$. The vertex $p_{\Cal V}(B_i)$ is a leaf-parent of the $\nu_{i+1}$-piece (Lemma \ref{lem:anc-visible}), and its height at $x_n$ tends to its height $\eta\leq h(B_i) = c_i$ at $x$. If $w^n_{\nu_{i+1}}$ is bounded below, then $w_{\nu_{i+1}} > 0$, so by (W5) we must have $\eta\geq c_i$, hence $\eta = c_i \geq E_{\nu_{i+1}} = \lim_n E^n_{\nu_{i+1}}$, and the rescaled height of $p_{\Cal V}(B_i)$ in the $\nu_{i+1}$-piece of $x_n$, namely $(h(p_{\Cal V}(B_i)) - S^n_{\nu_{i+1}})/w^n_{\nu_{i+1}}$, tends to a limit $\geq 1$; that factor tends to the basepoint, and so does the smash. Otherwise $w^n_{\nu_{i+1}}\to 0$, so $w_{\nu_{i+1}} = 0$ at some position $c_{i+1}$, all heights of the (valid) $\nu_{i+1}$-piece of $x_n$ tend to $c_{i+1}$, and $(\ast_{i+1})$ holds. Since the chain ends at a root-run position, whose window at $x$ is not degenerate by (W5), the second alternative cannot occur forever, and some factor tends to the basepoint. This completes the proof.
\end{proof}

\begin{remark}[Relation to the Boardman--Vogt resolution]\label{rem:BV-comparison}
Replacing the structuring operad by one with contractible fibers over $\Oper$ carrying both prescriptions among its cells is the rectification strategy of \cite{BV-1}, \cite{BV}, \cite{BM-trees}, and $\Lop$ has a shape familiar from it: its cells are combinatorial routes decorated with real parameters, as the cells of the $W$-construction are trees whose inner edges carry lengths in $[0,1]$ \cite{BM-WConstr}. The operad $\Kop$ adapts that device to \emph{cutting} rather than grafting, a cell assigning to each position not an edge length but a \emph{window} against which its piece is normalized. What is not formal in the adaptation is (W3): its two disjuncts, anchoring the floor at the parent vertex (the geometry of \cite{Ching-thesis}) and stacking windows consecutively (the Moore trick), are the two classical ways of making concatenation-like structure strict, and the disjunction is what makes both prescriptions cells of one operad. We know of no reference for windowed cut systems in this form.
\end{remark}

\begin{theorem}\label{thm:taut}
The maps
\begin{align*}
\eta(\theta, W)\colon \big(\text{smash of } \Omw(\uS)\text{-entries per the profile}\big) &\longrightarrow \Omw(\uS)[t], \\
(f_\rho)_\rho &\longmapsto \big({\textstyle\bigwedge}_\rho f_\rho\big)\circ \kappa_W,
\end{align*}
i.e., the smash pairing of cotensors followed by restriction along $\kappa_W$, define a $\Kop$-algebra structure $\tau$ on $\Omw(\uS)$ in the sense of \cite[Definition 6.32]{Clark-1}, and $o^{*}(\tau) = \nu^{\mathrm{tr}},$ the strict structure of Theorem \ref{thm:omw-operad} regarded as an $\Oper$-algebra.
\end{theorem}

\begin{proof}
The action maps are continuous in the cell variable by Proposition \ref{prop:K-continuity} (joint continuity), $\Sigma$-equivariant by (W1) and the equivariance of the smash pairing, and unital (the unit cell's full-window frames give the identity cut). Associativity is the identity
\begin{equation}\label{eq:action-assoc}
    \eta(\theta, W)\circ\big(\text{slot-}j \text{ insertion of } \eta(\theta',W')\big) \;=\; \eta\big((\theta,W)\circ_j(\theta',W')\big),
\end{equation}
This dualizes to an equality of two maps $\BC(t)\to \BC(\text{composite flag})$, which we check pointwise. The first is $\kappa_W$ followed by $\kappa_{W'}$ applied to all rescaled layer-$j$ pieces, the second is $\kappa$ of the substituted cell \eqref{eq:substitution}. Suppose first that the outer evaluation survives at a refined position $\rho$, so that the window is non-degenerate and no outer threshold trips. Then the affine equivariance (W4) identifies the substituted frames with conjugates of the inner top-window frames under the increasing affine window map,
\begin{align*}
W'\big( \text{raw piece};\ (S_\rho, E_\rho) \big) \ &=\ A_\rho\cdot W'\big( \text{rescaled piece};\ (0,1)\big), \\
A_\rho(x) &= S_\rho + (E_\rho - S_\rho)x,
\end{align*}
The two sides now agree. The first composite applies \eqref{eq:FS-formula} twice, once with the outer frames and once with the inner top-window frames on the rescaled piece. The composite of these two affine maps is the affine map of the substituted frames, and the degeneracy loci correspond threshold by threshold under the increasing bijection $A_\rho$. Suppose instead that the outer factor at $\rho$ degenerates. Then both sides are the basepoint, the first by based-ness of the smash and the second because $[S'',E'']\subseteq[S_\rho,E_\rho]$ by (W2) at the handed window, so each outer threshold forces a threshold on an insert piece.

It remains to compare $o^{*}(\tau)$ with $\nu^{\mathrm{tr}}$. Both are strict reduced operad structures \cite[Proposition 7.12]{Clark-1}, so it is enough that they have the same multiplication and unit. A two-level cell acts under $o^{*}(\tau)$ by smash pairing followed by restriction along $\kappa_{W^{\mathrm{tr}}} = \kappa_{F^{\mathrm{tr}}} = c^{\mathrm{tr}}$, by Examples \ref{ex:Wcells}(1) and \ref{ex:frames}(1), and this is the $\mu$ of Theorem \ref{thm:omw-operad}. Both units are the identity of $\Map(S^0,S)$. Translating into the orbit-assembled diagrams of \cite[Definition 6.32]{Clark-1} is routine, since no step above chooses orbit representatives.
\end{proof}

\subsection{The levels operad and uniform cuts}\label{sec:L-operad}

\begin{definition}\label{def:L}
Let $\theta\in \Oper_\ell(\underline{p};t)$ have underlying flag $\gamma(\theta)$. A \emph{cut labelling} for $\theta$ assigns a height $e_\nu\in(0,1)$ to each position $\nu = (j,B)$ of $\gamma(\theta)$ with $0\leq j\leq \ell-2$, equivariantly in $\theta$ and strictly increasing along the tree, meaning $e_{(j+1,C)} > e_{(j,B)}$ whenever $C\subseteq B$. Set $e_{(\ell-1,\cdot)} := 1$ and $e_{p(\text{root position})} := 0$. A cut labelling then determines constant \emph{windows} $[\,e_{\text{parent position}},\, e_\nu\,]$, one per position, stacked consecutively along every root-to-leaf path of the flag. Write $\mathring{W}(\theta)$ for the space of cut labellings. It is the interior of the order polytope of the truncated flag tree, hence a nonempty open convex set, hence contractible. The \emph{levels operad} $\Lop$ is the $\Nl$-operad with
\[ \Lop_\ell(\underline{p};t) := \coprod_{\theta\in\Oper_\ell(\underline{p};t)} \mathring{W}(\theta)\]
a point for $\ell = 1$, with $\Sigma$-structure by relabelling. Composition is in the format of \cite[\S 6]{Clark-1}, an outer cell together with one inner cell per position, of common level-height per layer: the composite labelling keeps the outer labels and inserts each inner cell's labels affinely into the window of its position. Associativity and unitality follow from associativity of affine maps, and strict monotonicity is preserved. The projection $\rho_\Lop\colon \Lop\to\Oper$ is a map of $\Nl$-operads and a levelwise weak equivalence.

The \emph{uniform} labellings, $e_{(j,B)} = a_{j+1}$ with $0 < a_1 < \cdots < a_{\ell-1} < 1$, form the open simplices $\mathring{\Delta}(\ell)\subseteq \mathring{W}(\theta)$ of Example \ref{ex:frames}(2). They generate $\Lop$ under composition, by peeling layers off from the root, but are not closed under it, since an $\ell$-cell composes with a \emph{family} of inner cells \cite[\S 6]{Clark-1}; this is why the cells of $\Lop$ are labelled per position.
\end{definition}

\begin{lemma}\label{lem:sigma-strict}
The assignment
\[ \sigma\colon \Lop\longrightarrow \Kop, \quad\quad (\theta, e)\mapsto (\theta, W^{e})\]
(constant cells, Example \ref{ex:Wcells}(2)) is a map of $\Nl$-operads with $\rho_{\Kop}\circ\sigma = \rho_{\Lop}$, and a levelwise weak equivalence.
\end{lemma}

\begin{proof}
Units, $\Sigma$-equivariance and the covering of $\Oper$ are clear from the formulas. Compositions: substituting constant cells into constant cells, the handed window at a refined position $\rho$ is
\[ (S_\rho, E_\rho) = \big( u + (v-u)s_\rho,\ u + (v-u)e_\rho \big), \quad\quad (s_\rho, e_\rho) := (e_{\pi(\rho)}, e_\rho),\]
so at an insert position with inner labels $(b, b')$ (consecutive inner window ends) the substituted frame \eqref{eq:substitution} is
\[ S_\rho + (E_\rho - S_\rho)\,b \ =\ u + (v-u)\big( s_\rho + (e_\rho - s_\rho) b\big), \]
and likewise for $b'$: the constant cell of the \emph{merged} labelling, since $s_\rho + (e_\rho - s_\rho)b$ is exactly the affine insertion of the inner labels into the outer window that defines the composition of $\Lop$ (Definition \ref{def:L}). The identity holds at every point of $\Dg{\gamma}\times\widehat{W}$, degenerate windows included; no case distinction occurs because no $\max$ occurs. Frames at unrefined positions are untouched and match the merged labelling. Finally, $\sigma$ is levelwise a map between disjoint unions, over the same discrete sets of $\theta$'s, of contractible spaces ($\mathring{W}(\theta)$ by Definition \ref{def:L}; $\mathrm{WS}(\gamma(\theta))$ by Proposition \ref{prop:K-operad}), hence a levelwise weak equivalence.
\end{proof}

\begin{lemma}[Equal evaluations]\label{lem:equal-evals}
For every $(\theta, e)$, the cut map of the constant cell equals the cut map of the labelled frames of Example \ref{ex:frames}(2): $\kappa_{\,\sigma(\theta,e)\,} = \kappa_{\,F^{e}}$.
\end{lemma}

\begin{proof}
The top-window slice of $W^{e}$ is the stacked constant family $(c_\nu, e_\nu)$ with $c_\nu := e_{\pi(\nu)}$, while $F^{e} = (\max(h(p_{\Cal{V}}(B)), c_\nu),\, e_\nu)$. It suffices to show the two frame families agree at every $x\in\BC(t)$ at which either map is not the basepoint, since each map is determined by its frames on its non-basepoint locus. At such $x$ every factor is a valid weighting. Fix $\nu = (j,B)$. If $\nu$ is root-run then $h(p_{\Cal{V}}(B)) = 0\leq c_\nu$ by convention. Otherwise the $\anc(\nu)$-factor is valid, so $p_{\Cal{V}}(B)$ has rescaled height $\leq 1$ in it (Lemma \ref{lem:anc-visible}), giving
\[ h\big(p_{\Cal{V}}(B)\big) \ \leq\ E_{\anc(\nu)} \ =\ e_{\anc(\nu)} \ \leq\ e_{\pi(\nu)} \ =\ c_\nu,\]
the last step because the labels increase along the run from $\anc(\nu)$ to $\pi(\nu)$. So the $\max$ is inactive and the two systems agree at $x$.
\end{proof}

\begin{corollary}[The uniform-cut structure]\label{cor:zeta}
The restriction $\zeta := \sigma^{*}(\tau)$ is an $\Lop$-algebra structure on $\Omw(\uS)$ whose action maps are the labelled-frame cuts of Example \ref{ex:frames}(2):
\[\zeta(\theta, e) = \big({\textstyle\bigwedge} - \big)\circ \kappa_{F^{e}}^{*}.\]
\end{corollary}

\begin{proof}
Restriction along the map of $\Nl$-operads $\sigma$ (Remark \ref{rem:morphism-convention}) makes $\zeta$ an $\Lop$-algebra with no further verification; the displayed formula is Lemma \ref{lem:equal-evals}.
\end{proof}

\subsection{Multi-cut cells and the transport}\label{sec:j-transport}

The following is the classical cup-product cell, normalized at height $\tfrac12$.

\begin{construction}[The half-cut]\label{con:halfcut}
For $p+q=m$ define $\chi_{p,q}\colon \nabla^m \supseteq \{u : u_p\leq \tfrac12 \leq u_{p+1}\} \to \nabla^p\times \nabla^q$ (with the conventions $u_0:=0$, $u_{m+1}:=1$) by
\[ \chi_{p,q}(u) = \big( (2u_1,\dots,2u_p),\; (2u_{p+1}-1,\dots,2u_m-1) \big). \]
The domains cover $\nabla^m$ as $(p,q)$ ranges over $p+q=m$, overlapping where some $u_j=\tfrac12$. On the overlap $u_{p+1}=\tfrac12$ (shared by $(p,q+1)$ and $(p+1,q)$) we have
\begin{align*}
\chi_{p,q+1}(u)   &= \big( (2u_1,\dots,2u_p),\; (0, 2u_{p+2}-1,\dots) \big) \\
\chi_{p+1,q}(u)   &= \big( (2u_1,\dots,2u_p, 1),\; (2u_{p+2}-1,\dots) \big)
\end{align*}
The two descriptions differ by inserting height $0$ at the front of the second factor rather than appending height $1$ at the back of the first. By \eqref{eq:heights-cofaces} these are the maps $\id\times \delta^{0}$ and $\delta^{p+1}\times \id$ along which the box-product colimit is formed. The $\chi_{p,q}$ therefore glue to a well-defined map to that colimit in each cosimplicial degree. (Over a profile $(n,(k_1,\dots,k_n))$ the target of a cell is $\Delta^{p}\times(\Delta^{q})^{\times n}$, one factor per block; the cell uses the \emph{same} rescaled back tuple in all $n$ factors, i.e., it is composed with the diagonal. This is what matches the decomposition $\Psi_{k,(p,q)}$ of Lemma \ref{lem:tech}, in which all the $\beta_j$ have the same degree $q$.) Compatibility with the cosimplicial structure maps of source and target is a direct check from \eqref{eq:heights-cofaces}. We write $\psi_0\in \Aop_2$ for the resulting element of the coendomorphism operad, and call it the \emph{half-cut}.
\end{construction}

\begin{remark}\label{rem:naive-frontback-fails}
The unrescaled front/back formula $u\mapsto ((u_1,\dots,u_p),(u_{p+1},\dots,u_m))$ does \emph{not} descend to the box-product colimit: its two descriptions on a putative overlap disagree unless $u_{p}$ or $u_{p+1}$ is extreme, and no cover of $\nabla^m$ repairs this. Rescaling by $[0,\tfrac12]\to[0,1]$ and $[\tfrac12,1]\to [0,1]$ is what makes the gluing possible. This is why height $\tfrac12$, rather than a ``cut nothing'' convention, is forced on the cosimplicial side, and it is the source of the discrepancy with Ching's cuts, which are tied to no fixed height.
\end{remark}

\begin{construction}[Multi-cut cells]\label{con:multicut}
Fix $0 < a_1 < \cdots < a_{\ell-1} < 1$ and generalize Construction \ref{con:halfcut} as follows. Let $(p_1,\dots,p_\ell)$ be a multidegree with $\sum p_i = m$, and write $P_j := p_1+\cdots+p_j$, with the conventions $u_0 = 0$ and $u_{m+1}=1$. On the region where $u_{P_j}\leq a_j\leq u_{P_j + 1}$ for $j = 1,\dots,\ell-1$, define a map to $\nabla^{p_1}\times\cdots\times\nabla^{p_\ell}$ by sending $u$ to the blockwise affinely rescaled tuples
\[ \Big( \frac{u_{P_{j-1}+1} - a_{j-1}}{a_j - a_{j-1}},\, \dots,\, \frac{u_{P_j} - a_{j-1}}{a_j - a_{j-1}} \Big)_{j = 1,\dots,\ell}.\]
On an overlap, where $u_i = a_j$ exactly, the two adjacent multidegree descriptions differ by prepending a $0$ to one block rather than appending a $1$ to the neighbouring one. These are the identifications defining the iterated box-product colimit, as in Construction \ref{con:halfcut}, and compatibility with cofaces and codegeneracies is again a direct check from \eqref{eq:heights-cofaces}. We obtain a cell $\psi_{\vec a}\in \Aop_\ell$ over the appropriate profile components, natural in $\vec a$, with $\psi_{(1/2)} = \psi_0$. Since affine maps compose, it satisfies
\[ \psi_{\vec a}\circ (\psi_{\vec b_\rho})_\rho = \psi_{\vec a \circ (\vec b_\rho)} \]
whenever the inner labellings are level-uniform across each layer. Indeed both sides are given blockwise by rescaling along the merged windows, and the box-colimit comparisons are the canonical ones in both cases. For $\Sigma$-copowered objects these are the strict identifications of the underlying $\Top^{\bdel}$-box products \cite[\S 7]{Clark-1}.
\end{construction}

\begin{definition}\label{def:j}
On a uniformly labelled cell define $j(\theta, \vec a) := \theta\cdot \psi_{\vec a}$, the label-free cell $\psi_{\vec a}$ placed over the component indexed by $\theta$. On a general $(\theta, e)$, express $e$ as an iterated composite of uniform cells by peeling layers from the root, as in Definition \ref{def:L}, and take the corresponding composite in $\Aop$. This is well defined and a map of $\Nl$-operads: on the underlying $\Top^{\bdel}$-side, where the box product is strictly monoidal \cite[\S 7]{Clark-1}, a composite of uniform cells evaluates to the map that rescales the coordinates of each block of $\lambda_j$ affinely from its window $[e_{\pi(\nu)}, e_\nu]$ onto $[0,1]$, so the composite depends only on the labelling $e$ and not on the word chosen, and the composite of two such maps is the map of the merged labelling of Definition \ref{def:L}, because affine maps compose. Moreover $\rho\circ j = \rho_{\Lop}$, since the collapse $\Delta^\bullet\to \ast$ sends every $\psi$-cell to the identity cell over its shape \cite[Proposition 7.16]{Clark-1}. Hence $j$ is a levelwise weak equivalence by two-out-of-three.
\end{definition}

\begin{proposition}[Transport]\label{prop:transport}
The levelwise equivalence $\beta\colon \Omw(\uS)\to \Tot C(\uS)$ of Proposition \ref{prop:bridge} strictly intertwines the $\Lop$-algebra structures: for every $(\theta, e)$,
\[ \beta\circ \zeta\big(\theta, e\big) \;=\; \lambda^{\mathrm{lev}}\big(j(\theta, e)\big)\circ \beta^{\wedge},\]
where $\beta^{\wedge}$ denotes the induced map on smash products of entries. That is, $\beta$ is a strict map of $\Lop$-algebras from $(\Omw(\uS), \zeta)$ to $(\Tot C(\uS), j^{*}\lambda^{\mathrm{lev}})$.
\end{proposition}

\begin{proof}
Every object here is a cotensor of $S$ by Lemma \ref{lem:tot-is-map}, and every structure map is a restriction along a map of spaces. We may therefore dualize and compare the two sides cell by cell on $\fat{P(t)}$.

Take a cell $(\gamma', u)$ with $\gamma'\in P(t)^{\circ}_m$, and let $\vec p$ be the multidegree determined by $u$ and $\vec a$ as in Construction \ref{con:multicut}. By Remark \ref{rem:clark1-interface}(ii) the value of $\lambda^{\mathrm{lev}}(\psi_{\vec a})$ at this cell is the basepoint unless $\gamma'$ decomposes under the iterated $\Psi$-maps of Lemma \ref{lem:tech} at the slot boundaries, that is, unless $\lambda'_{P_j} = \lambda_j(\gamma(\theta))$ for all $j$. Otherwise it is the smash of the input tuple at the $\Psi$-decomposition segments, with the rescaled coordinates of Construction \ref{con:multicut}. On overlaps the two adjacent descriptions agree by the box identifications.

Now transport through $\Thet$. We have $u_{P_j}\leq a_j\leq u_{P_j+1}$, so $a_j$ lies in the interval $(u_{P_j}, u_{P_j+1}]$ on which $\lambda'_{P_j}$ is the partition in force. By Lemma \ref{lem:theta-dictionary}(2) this partition is the alive partition $\alive{\Thet(\gamma',u)}{a_j}$. In the tied case $u_{P_j}=a_j$ the two adjacent descriptions are glued by the box identification. The validity locus is therefore $\{\alive{\Thet(\gamma',u)}{a_j} = \lambda_j \text{ for all } j\}$, which by Example \ref{ex:frames}(2) is the non-basepoint locus of $\kappa_{F^{\vec a}}$.

On that locus, consider the segment of $\gamma'$ between positions $P_j$ and $P_{j+1}$, restricted to a block $B$ of $\lambda_j$ and rescaled by $[a_j,a_{j+1}]\to[0,1]$. Its $\Thet$-image is the piece $\Cal{V}_{(j,B)}$ with the heights \eqref{eq:FS-formula} for $F^{\vec a}$. This is Construction \ref{con:theta} applied to the segment, the last-level bookkeeping matching the piece description of \S\ref{sec:flags}. So the two sides agree, and units and equivariance are Lemma \ref{lem:units-bridge} and (F1).

This covers the uniform cells. A general labelled cell is an operadic composite of uniform cells, by Definitions \ref{def:L} and \ref{def:j}, and the same composite expresses both sides. Both $\zeta$ and $j^{*}\lambda^{\mathrm{lev}}$ send operadic composites of cells to composites of action maps. The identity therefore extends to all of $\Lop$.
\end{proof}

\subsection{The comparison theorem}\label{sec:main}

\begin{theorem}\label{thm:main}
Let $\lambda^{\mathrm{lev}}$ be the $A_\infty$-operad structure of Theorem \ref{thm:main-1} on $\Tot C(\uS)$, and $\nu^{\mathrm{tr}}$ Ching's operad structure \cite[\S\S 4--8]{Ching-thesis} on $\Omw(\uS) = \mathbb{D}B(\uSz)$, expressed in height coordinates as Theorem \ref{thm:omw-operad}; both underlying symmetric sequences are models for $\partial_*\Id$. Then $(\Tot C(\uS), \lambda^{\mathrm{lev}})$ and $(\Omw(\uS), \nu^{\mathrm{tr}})$ agree through the cut refinement, in the sense of Definition \ref{def:structures-agree}. Specifically:
\begin{enumerate}[label=(\arabic*)]
    \item there is a diagram of $\Nl$-operads and maps of $\Nl$-operads
    \[ \Aop \xleftarrow{\ j\ } \Lop \xrightarrow{\ \sigma\ } \Kop \xleftarrow{\ o\ } \Oper\]
    ($\Lop$: Definition \ref{def:L}; $\Kop$: Definition \ref{def:K}; $j$: Definition \ref{def:j}; $\sigma$: Lemma \ref{lem:sigma-strict}; $o$, $\rho_\Kop$: Proposition \ref{prop:K-operad}), in which every map is a levelwise weak equivalence, with $\rho_\Kop\circ o = \id$ and $\rho_\Kop\circ\sigma = \rho_\Lop$;
    \item the tautological $\Kop$-algebra $\tau$ on $\Omw(\uS)$ (Theorem \ref{thm:taut}) satisfies $o^{*}(\tau) = \nu^{\mathrm{tr}}$;
    \item the restriction $\zeta = \sigma^{*}(\tau)$ (Corollary \ref{cor:zeta}) is strictly intertwined by the levelwise weak equivalence $\beta$ of Proposition \ref{prop:bridge} with $j^{*}(\lambda^{\mathrm{lev}})$ (Proposition \ref{prop:transport}).
\end{enumerate}
In particular, the operad structure of Theorem \ref{thm:main-1} on the Goodwillie derivatives of the identity in spaces and the operad structure constructed by Ching in \cite{Ching-thesis} are restrictions of the single coherent $\Kop$-algebra $\tau$ along maps of $\Nl$-operads, and thus equivalent as $A_\infty$-operads.
\end{theorem}

The comparison assembles into the diagram
\[
\xymatrix@R=2.1em@C=2.6em{
 & \Aop & \Lop \ar[l]_-{j}^-{\sim} \ar[r]^-{\sigma}_-{\sim} & \Kop & \Oper \ar[l]_-{o}^-{\sim} \\
\text{on } \Tot C(\uS)\colon & \lambda^{\mathrm{lev}} \ar[r]^-{j^{*}} & j^{*}\lambda^{\mathrm{lev}} & & \\
\text{on } \Omw(\uS)\colon & & \zeta = \sigma^{*}\tau \ar[u]^-{\beta}_-{\wr} & \tau \ar[l]_-{\sigma^{*}} \ar[r]^-{o^{*}} & \nu^{\mathrm{tr}}
}
\]
in which each algebra stands in the column of an operad acting on it. Every restriction arrow is an equality on the nose; the homotopy-theoretic content is in the top row, where the structuring operads are connected by equivalences rather than isomorphisms, as it must be, since the structure of Theorem \ref{thm:main-1} is not strict.

\begin{proof}
Collecting the citations in the statement: (1) is Definition \ref{def:j}, Lemma \ref{lem:sigma-strict} and Proposition \ref{prop:K-operad}; (2) is Theorem \ref{thm:taut}; (3) is Corollary \ref{cor:zeta} and Proposition \ref{prop:transport}. These are the data (i)--(iii) of Definition \ref{def:structures-agree}.
\end{proof}

\section{An explicit comparison in arity 3}\label{sec:anchor-low}

This appendix exhibits the comparison in the lowest nontrivial arity. We follow one map throughout,
\[ \BC(3)\longrightarrow \BC(2)\wedge\BC(2)\wedge\BC(1),\]
whose dual is a component of $\partial_2\Id\wedge\partial_2\Id\wedge\partial_1\Id\to\partial_3\Id$, and compute it twice, by Ching's recipe and by that of Theorem \ref{thm:main-1}. Iterating shows why one is strictly associative and the other only $A_\infty$, and recording the mismatch produces $\Lop$ and $\Kop$.

\begin{example}\label{ex:B2-B3}
$\BC(1) = S^0$, and $\BC(2) = [0,1]/\{0,1\}\cong S^1$, the unique tree on $\mathbf{2}$ having one internal vertex of height $h$, with $h\in\{0,1\}$ the basepoint by Definition \ref{def:B}(b). For $n=3$ (cf.\ Example \ref{ex:par-low}) there are four trees, drawn in Figure \ref{fig:B3}: the corolla $\TreeF_c$, with weightings an interval $[0,1]_c$, and for each $i$ the binary tree $\TreeF_i$ containing $\overline{i} := \mathbf{3}\smallsetminus\{i\}$, with weightings the triangle $\nabla^2_i = \{0\leq x\leq y\leq 1\}$, $x = h(\mathbf{3})$, $y = h(\overline{i})$. Collapse glues each hypotenuse $x=y$ to the corolla interval and (b) sends the rest of the boundary to the basepoint, giving the $2$-sphere with an equatorial disc, $\Sigma_3$ permuting the three discs and fixing the circle. Under $\Thet$ the triangle is the image of the $2$-chain $\mathsf{min} < \{\overline{i}\,|\,i\} < \mathsf{max}$ with $(u_1,u_2) = (x,y)$.
\end{example}

\begin{figure}[htbp]
\centering
\begin{tikzpicture}[x=1cm,y=1cm,line width=0.8pt]
\begin{scope}[yshift=-1.15cm,y=2.3cm]
  \draw[gray!45,dotted] (-0.95,0) -- (0.75,0);
  \draw[gray!45,dotted] (-0.95,1) -- (0.75,1);
  \node[gray!80,left,scale=.62] at (-0.95,0) {$0$};
  \node[gray!80,left,scale=.62] at (-0.95,1) {$1$};
  \coordinate (cv) at (0,0.45);
  \draw (0,0)--(cv);
  \draw (cv)--(-0.5,1); \draw (cv)--(0,1); \draw (cv)--(0.5,1);
  \draw[gray!55,dotted] (-0.95,0.45) -- (0,0.45);
  \node[gray!80,left,scale=.62] at (-0.95,0.45) {$z$};
  \fill (cv) circle (1.6pt);
  \node[right,scale=.68,inner sep=3pt] at (cv) {$\mathbf{3}$};
  \node[above,scale=.6] at (-0.5,1) {$1$};
  \node[above,scale=.6] at (0,1) {$2$};
  \node[above,scale=.6] at (0.5,1) {$3$};
  \node[scale=.7,align=center,anchor=north] at (-0.1,-0.07) {$\TreeF_c$\\[-2pt] $[0,1]_c$};
\end{scope}
\begin{scope}[xshift=2.75cm,yshift=-1.15cm,y=2.3cm]
  \draw[gray!45,dotted] (-1.15,0) -- (0.8,0);
  \draw[gray!45,dotted] (-1.15,1) -- (0.8,1);
  \node[gray!80,left,scale=.62] at (-1.15,0) {$0$};
  \node[gray!80,left,scale=.62] at (-1.15,1) {$1$};
  \coordinate (v0) at (0,0.3);
  \coordinate (v1) at (-0.5,0.72);
  \draw (0,0)--(v0)--(v1)--(-0.82,1);
  \draw (v1)--(-0.18,1);
  \draw (v0)--(0.65,1);
  \draw[gray!55,dotted] (-1.15,0.3) -- (0,0.3);
  \node[gray!80,left,scale=.62] at (-1.15,0.3) {$x$};
  \draw[gray!55,dotted] (-1.15,0.72) -- (-0.5,0.72);
  \node[gray!80,left,scale=.62] at (-1.15,0.72) {$y$};
  \fill (v0) circle (1.6pt);
  \fill (v1) circle (1.6pt);
  \node[below right,scale=.68,inner sep=2.5pt] at (v0) {$\mathbf{3}$};
  \node[right,scale=.68,inner sep=3pt] at (v1) {$\overline{\imath}$};
  \node[above,scale=.6] at (-0.82,1) {$1$};
  \node[above,scale=.6] at (-0.18,1) {$2$};
  \node[above,scale=.6] at (0.65,1) {$3$};
  \node[scale=.7,align=center,anchor=north] at (-0.15,-0.07) {$\TreeF_i$\\[-2pt] $\nabla^2_i$};
\end{scope}
\begin{scope}[xshift=7.75cm]
  \filldraw[fill=gray!12] (0,1.25) .. controls (-2.05,0.8) and (-2.05,-0.8) .. (0,-1.25) -- cycle;
  \filldraw[fill=gray!12] (0,1.25) .. controls (2.05,0.8)  and (2.05,-0.8)  .. (0,-1.25) -- cycle;
  \filldraw[fill=gray!30] (0,1.25) .. controls (-0.78,0.5) and (-0.78,-0.5) .. (0,-1.25) -- cycle;
  \draw[line width=1.3pt] (0,1.25) -- (0,-1.25);
  \fill (0,1.25) circle (2pt);
  \fill (0,-1.25) circle (2pt);
  \node[scale=.75] at (-1.28,0) {$\nabla^2_1$};
  \node[scale=.75] at (1.28,0)  {$\nabla^2_3$};
  \node[scale=.75] at (-2.72,0.78) {$\nabla^2_2$};
  \draw[->,line width=0.5pt,gray!70] (-2.45,0.7) -- (-0.42,0.3);
  \node[scale=.68,anchor=west] at (0.16,-0.95) {$[0,1]_c$};
  \node[scale=.68,anchor=west] at (2.25,1.15) {$\basept$};
  \node[scale=.68,anchor=west] at (2.25,-1.15) {$\basept$};
  \draw[->,line width=0.5pt,gray!70] (2.2,1.15) -- (0.1,1.28);
  \draw[->,line width=0.5pt,gray!70] (2.2,-1.15) -- (0.1,-1.28);
\end{scope}
\end{tikzpicture}
\caption{The trees on $\mathbf{3}$ and the space they assemble into. \emph{Left:} the corolla $\TreeF_c$ and the binary tree $\TreeF_i$, whose weightings form an interval and a triangle. \emph{Right:} the gluing. Letting $y\to x$ attaches each triangle to the corolla along its hypotenuse, and the rest of the boundary goes to the basepoint, leaving three discs sharing one circle, drawn edge-on.}
\label{fig:B3}
\end{figure}

We write $x = h(\mathbf{3})$, $y = h(\overline{i})$ on $\nabla^2_i$, and record a non-basepoint point of $\BC(2)$ by the height of its internal vertex.

\subsection*{One map, two recipes}
The arity-$3$ multiplications are indexed by the chains of $P(3)^{\circ}_2$, with middle entry $\lambda\in\{\mathsf{min},\alpha_1,\alpha_2,\alpha_3,\mathsf{max}\}$ where $\alpha_i := \{\overline{\imath},\{i\}\}$. The outer two are unit maps, so we fix $\alpha_i$; its cut map is
\[ \BC(3)\longrightarrow \BC(2)\wedge \BC(\overline{\imath})\wedge\BC(\{i\}) \;\cong\; S^1\wedge S^1\wedge S^0,\]
one factor for the bottom tree and one per block, the $S^0$ suppressed from now on.

Both recipes cut the weighted tree along a horizontal line: if the line separates $\mathbf{3}$ from $\overline{\imath}$, the tree is severed into a tree on the two blocks below and a tree on $\overline{\imath}$ above, each normalized onto $[0,1]$ from its \emph{window}, the interval of heights it occupies; otherwise the value is the basepoint. Ching draws the line at the parent vertex of $\overline{\imath}$, at height $\tau = x$; the structure of Theorem \ref{thm:main-1} draws it at the fixed height $\tau=\tfrac12$ (Remark \ref{rem:naive-frontback-fails}). On $\nabla^2_i$, the other strata having no vertex $\overline{i}$,
\begin{align}
    c^{\mathrm{tr}}_{\alpha_i}(x,y) &= \Big( x,\; \frac{y-x}{1-x} \Big), &&\text{valid for } x < y, \label{eq:tr-cut-3}\\
    c^{\mathrm{lev}}_{\alpha_i}(x,y) &= \big( 2x,\; 2y-1 \big), && \text{valid for } x < \tfrac12 \leq y. \label{eq:lev-cut-3}
\end{align}
In terms of windows (Figure \ref{fig:two-cuts}): Ching's are $[0,1]$ and $[x,1]$, \emph{nested}, sharing the ceiling $1$, so the bottom piece is not rescaled; the levelled ones are $[0,\tfrac12]$ and $[\tfrac12,1]$, \emph{stacked}, tiling $[0,1]$, so both pieces are rescaled.

\begin{figure}[htbp]
\centering
\begin{tikzpicture}[x=1cm,y=2.7cm,line width=0.9pt]
\begin{scope}
  \draw[gray!45,dotted] (-1.25,0) -- (2.45,0);
  \draw[gray!45,dotted] (-1.25,1) -- (2.45,1);
  \node[gray!80,left,scale=.7] at (-1.25,0) {$0$};
  \node[gray!80,left,scale=.7] at (-1.25,1) {$1$};
  \coordinate (r)   at (0,0);
  \coordinate (v0)  at (0,0.3);
  \coordinate (v1)  at (-0.55,0.72);
  \coordinate (lfa) at (-0.9,1);
  \coordinate (lfb) at (-0.22,1);
  \coordinate (lfc) at (0.7,1);
  \draw (r)--(v0)--(v1)--(lfa);
  \draw (v1)--(lfb);
  \draw (v0)--(lfc);
  \fill (v0) circle (1.5pt);
  \fill (v1) circle (1.5pt);
  \node[below right,scale=.72,inner sep=2pt] at (v0) {$\mathbf{3}$};
  \node[left,scale=.72,inner sep=2pt]  at (v1) {$\overline{\imath}$};
  \node[above,scale=.68] at (lfa) {$1$};
  \node[above,scale=.68] at (lfb) {$2$};
  \node[above,scale=.68] at (lfc) {$3$};
  \draw[dashed,line width=1pt] (-1.2,0.3) -- (1.1,0.3);
  \draw[gray!45,dotted] (1.1,0.3) -- (2.25,0.3);
  \node[scale=.72,above,inner sep=2pt] at (0.72,0.3) {$\tau=x$};
  \draw[|-|,line width=0.9pt] (1.55,0) -- (1.55,1);
  \node[right,scale=.7] at (1.6,0.13) {$[0,1]$};
  \draw[|-|,line width=0.9pt] (2.25,0.3) -- (2.25,1);
  \node[right,scale=.7] at (2.3,0.65) {$[x,1]$};
  \node[scale=.82] at (0.6,-0.2) {$c^{\mathrm{tr}}_{\alpha_i}$: \emph{nested}};
\end{scope}
\begin{scope}[xshift=5.9cm]
  \draw[gray!45,dotted] (-1.25,0) -- (1.75,0);
  \draw[gray!45,dotted] (-1.25,1) -- (1.75,1);
  \node[gray!80,left,scale=.7] at (-1.25,0) {$0$};
  \node[gray!80,left,scale=.7] at (-1.25,1) {$1$};
  \coordinate (R)   at (0,0);
  \coordinate (V0)  at (0,0.3);
  \coordinate (V1)  at (-0.55,0.72);
  \coordinate (LFA) at (-0.9,1);
  \coordinate (LFB) at (-0.22,1);
  \coordinate (LFC) at (0.7,1);
  \draw (R)--(V0)--(V1)--(LFA);
  \draw (V1)--(LFB);
  \draw (V0)--(LFC);
  \fill (V0) circle (1.5pt);
  \fill (V1) circle (1.5pt);
  \node[below right,scale=.72,inner sep=2pt] at (V0) {$\mathbf{3}$};
  \node[left,scale=.72,inner sep=2pt]  at (V1) {$\overline{\imath}$};
  \node[above,scale=.68] at (LFA) {$1$};
  \node[above,scale=.68] at (LFB) {$2$};
  \node[above,scale=.68] at (LFC) {$3$};
  \draw[dashed,line width=1pt] (-1.2,0.5) -- (1.1,0.5);
  \draw[gray!45,dotted] (1.1,0.5) -- (1.55,0.5);
  \node[scale=.72,above,inner sep=2pt] at (0.72,0.5) {$\tau=\tfrac12$};
  \draw[|-|,line width=0.9pt] (1.55,0) -- (1.55,0.5);
  \node[right,scale=.7] at (1.6,0.25) {$[0,\tfrac12]$};
  \draw[|-|,line width=0.9pt] (1.55,0.5) -- (1.55,1);
  \node[right,scale=.7] at (1.6,0.75) {$[\tfrac12,1]$};
  \node[scale=.82] at (0.6,-0.2) {$c^{\mathrm{lev}}_{\alpha_i}$: \emph{stacked}};
\end{scope}
\end{tikzpicture}
\caption{The same weighted tree $(\TreeF_i, (x,y))$, cut two ways: Ching at the parent vertex $\mathbf{3}$, height $\tau = x$; the box product at the fixed height $\tau = \tfrac12$. The bars record the windows, nested on the left and stacked on the right, which is the difference between \eqref{eq:tr-cut-3} and \eqref{eq:lev-cut-3}.}
\label{fig:two-cuts}
\end{figure}

Both maps extend continuously by the basepoint, since at every wall some factor degenerates before the wall is reached: at $y=x$ Ching's top vertex reaches $0$; at $x=\tfrac12$ the levelled bottom vertex reaches $1$, and at $y=\tfrac12$ its top vertex reaches $0$; the edges $x=0$ and $y=1$ are already the basepoint of $\BC(3)$. So $c^{\mathrm{tr}}_{\alpha_i}$ is defined wherever $\overline{\imath}$ exists at all, while $c^{\mathrm{lev}}_{\alpha_i}$ is defined only on the square $x<\tfrac12\leq y$, where the line at height $\tfrac12$ separates the two vertices.

\subsection*{The product of \texorpdfstring{$f$ and $g$}{f and g}, computed twice}
Take $(\TreeF_i,(x,y))$ with $x=\tfrac14$, $y=\tfrac34$, where both cuts are valid and both return a pair of $2$-leaf corollas, disagreeing only in the heights.
\begin{center}
\begin{tabular}{@{}llccl@{}}
 & piece & window & vertex & rescaled value \\
\hline
$c^{\mathrm{tr}}_{\alpha_i}$, cut at $\tau = x = \tfrac14$
   & bottom & $[0,1]$          & $\tfrac14$ & $\tfrac14$ (unchanged) \\
   & top    & $[\tfrac14,1]$   & $\tfrac34$ & $\frac{3/4-1/4}{1-1/4} = \tfrac23$ \\
\hline
$c^{\mathrm{lev}}_{\alpha_i}$, cut at $\tau = \tfrac12$
   & bottom & $[0,\tfrac12]$   & $\tfrac14$ & $\frac{1/4-0}{1/2-0} = \tfrac12$ \\
   & top    & $[\tfrac12,1]$   & $\tfrac34$ & $\frac{3/4-1/2}{1-1/2} = \tfrac12$ \\
\hline
\end{tabular}
\end{center}
So $c^{\mathrm{tr}}_{\alpha_i}(\tfrac14,\tfrac34) = (\tfrac14,\tfrac23)$ while $c^{\mathrm{lev}}_{\alpha_i}(\tfrac14,\tfrac34) = (\tfrac12,\tfrac12)$: \emph{the same two trees with different weightings}. Dualizing, for $f,g\in \Omw(\uS)[2] = \Map(\BC(2),S)$ the two products at this input are
\[ \big(f\cdot g\big)^{\mathrm{tr}} \;=\; f(\tfrac14)\wedge g(\tfrac23), \quad\quad\quad \big(f\cdot g\big)^{\mathrm{lev}} \;=\; f(\tfrac12)\wedge g(\tfrac12),\]
which feed $f$ and $g$ different points of $\BC(2)$.

\subsection*{Cutting twice}
Associativity here means that cutting twice, in either order, yields the same threefold decomposition. Ching's cuts pass: the second line is again drawn at a parent vertex of the same ambient tree, so either order severs at the same heights, the window at a block being $[h(p(B)),1]$ in both associations. Nested windows are stable under re-cutting. This is the coassociativity of \cite[\S 4]{Ching-thesis}, and dualizing gives the \emph{strict} operad of Theorem \ref{thm:omw-operad}.

The levelled cuts fail. Cut at $\tfrac12$, then cut one piece at \emph{its} half: the inner line falls at ambient height $\tfrac14$ or $\tfrac34$ according to which piece is re-cut, so the two orders produce
\[ \big\{[0,\tfrac14],\ [\tfrac14,\tfrac12],\ [\tfrac12,1]\big\} \quad\quad\text{and}\quad\quad \big\{[0,\tfrac12],\ [\tfrac12,\tfrac34],\ [\tfrac34,1]\big\},\]
different window systems and hence different maps. Stacked windows are not stable under re-cutting: halving a half gives a quarter, and no set of heights fixed in advance is closed under this. So the levelled cuts form no cooperad, and their duals assemble only into a multiplication associative up to coherent homotopy, the $A_\infty$ structure of Theorem \ref{thm:main-1}; in the terms of \S\ref{sec:loop-analogy}, the two window systems above are the two bracketings of a triple concatenation on $[0,1]$.

\subsection*{Where \texorpdfstring{$\Lop$ and $\Kop$}{L and K} come from}
Let a cell be a chain together with its cut heights, so that the two iterations above are two different cells, and let composition insert the inner heights affinely into the window of their slot (the arithmetic $\tfrac12\mapsto\tfrac14$ above); this is strictly associative because affine maps compose. This is the levels operad $\Lop$ of Definition \ref{def:L}, whose fibre over a chain is the open simplex of admissible height tuples: contractible, so nothing is lost homotopically, but not a point, which is what carries the $A_\infty$ structure. Through $j$ of Definition \ref{def:j}, $\Lop$ is where the cosimplicial side lives.

Ching's prescription is not among these cells, since $\tau = h(p(\overline{\imath}))$ is not a constant. So let the datum at a block be a \emph{rule} for its window, allowed to read the weighted tree and the window handed down. Not every rule is legal. The constraint appears at the wall $y=x$, where $\overline{\imath}$ collapses into $\mathbf{3}$ and the value must degenerate, the corolla beyond having no such vertex to cut at. Two mechanisms achieve this. A floor at or above the parent vertex, called \emph{anchoring}, makes the piece's root edge vanish as $y\to x$. A floor equal to the previous window's ceiling, called \emph{stacking}, hands the parent vertex to the neighbouring factor, which degenerates instead. Ching uses the first mechanism and the levelled cut the second. A rule with neither leaves a gap. For instance, setting the top window's floor to a constant $c$ outright gives
\[ H = \Big(x,\ \frac{y-c}{1-c}\Big), \quad\quad \text{windows } [0,1] \text{ and } [c,1].\]
On a tree with $x>c$ the interval $[c,x]$ lies in the top window but not in the top piece. As $y\to x$ the top piece therefore keeps a root edge of length $\tfrac{x-c}{1-c}>0$ and never reaches the basepoint. But the wall it approaches is glued to the corolla, where the value is the basepoint. So $H$ does not descend to $\BC(3)$ (Figure \ref{fig:legs}). Replacing $c$ by $\max(x,c)$ repairs this, as in \eqref{eq:leg1-3}.

\begin{figure}[htbp]
\centering
\begin{tikzpicture}[x=1cm,y=2.4cm,line width=0.9pt]
\begin{scope}
  \draw[gray!55,dash pattern=on 2pt off 1.5pt,line width=0.5pt] (-0.2,0.3) -- (0.66,0.3);
  \node[scale=.6,left,inner sep=1.5pt] at (-0.2,0.3) {$x$};
  \draw[|-|] (0.1,0) -- (0.1,1);
  \draw[|-|] (0.45,0.3) -- (0.45,1);
  \node[scale=.6,align=center,anchor=north] at (0.28,-0.08) {$c=0$\\[-2pt] (Ching)};
\end{scope}
\draw[->,line width=0.6pt] (0.8,0.78) -- (1.3,0.78);
\node[scale=.6,anchor=south,inner sep=2pt] at (1.05,0.80) {leg 1};
\begin{scope}[xshift=1.6cm]
  \draw[gray!55,dash pattern=on 2pt off 1.5pt,line width=0.5pt] (-0.2,0.3) -- (0.66,0.3);
  \node[scale=.6,left,inner sep=1.5pt] at (-0.2,0.3) {$x$};
  \draw[gray!55,dash pattern=on 2pt off 1.5pt,line width=0.5pt] (-0.2,0.5) -- (0.66,0.5);
  \node[scale=.6,left,inner sep=1.5pt] at (-0.2,0.5) {$\tfrac12$};
  \draw[|-|] (0.1,0) -- (0.1,1);
  \draw[|-|] (0.45,0.5) -- (0.45,1);
  \node[scale=.6,align=center,anchor=north] at (0.28,-0.08) {$c=\tfrac12$, i.e.\\[-2pt] $w=1$};
\end{scope}
\draw[->,line width=0.6pt] (2.4,0.78) -- (2.9,0.78);
\node[scale=.6,anchor=south,inner sep=2pt] at (2.65,0.80) {leg 2};
\begin{scope}[xshift=3.2cm]
  \draw[gray!55,dash pattern=on 2pt off 1.5pt,line width=0.5pt] (-0.2,0.3) -- (0.66,0.3);
  \node[scale=.6,left,inner sep=1.5pt] at (-0.2,0.3) {$x$};
  \draw[gray!55,dash pattern=on 2pt off 1.5pt,line width=0.5pt] (-0.2,0.5) -- (0.66,0.5);
  \node[scale=.6,left,inner sep=1.5pt] at (-0.2,0.5) {$\tfrac12$};
  \draw[|-|] (0.1,0) -- (0.1,0.75);
  \draw[|-|] (0.45,0.5) -- (0.45,1);
  \node[scale=.6,align=center,anchor=north] at (0.28,-0.08) {$w=\tfrac34$\\[-2pt] (mid-motion)};
\end{scope}
\draw[->,line width=0.6pt] (4,0.78) -- (4.5,0.78);
\begin{scope}[xshift=4.8cm]
  \draw[gray!55,dash pattern=on 2pt off 1.5pt,line width=0.5pt] (-0.2,0.3) -- (0.66,0.3);
  \node[scale=.6,left,inner sep=1.5pt] at (-0.2,0.3) {$x$};
  \draw[gray!55,dash pattern=on 2pt off 1.5pt,line width=0.5pt] (-0.2,0.5) -- (0.66,0.5);
  \node[scale=.6,left,inner sep=1.5pt] at (-0.2,0.5) {$\tfrac12$};
  \draw[|-|] (0.1,0) -- (0.1,0.5);
  \draw[|-|] (0.45,0.5) -- (0.45,1);
  \node[scale=.6,align=center,anchor=north] at (0.28,-0.08) {$w=\tfrac12$\\[-2pt] (levelled)};
\end{scope}
\begin{scope}[xshift=7.5cm]
  \draw[gray!55,dash pattern=on 2pt off 1.5pt,line width=0.5pt] (-0.2,0.3) -- (0.66,0.3);
  \node[scale=.6,left,inner sep=1.5pt] at (-0.2,0.3) {$x$};
  \node[scale=.6,left,inner sep=1.5pt] at (-0.2,0.15) {$c$};
  \draw[|-|] (0.1,0) -- (0.1,1);
  \draw[|-|] (0.45,0.15) -- (0.45,1);
  \draw[line width=3pt] (0.45,0.15) -- (0.45,0.3);
  \draw[->,line width=0.5pt] (1.34,0.22) -- (0.53,0.22);
  \node[scale=.55,align=left,anchor=west] at (1.36,0.22)
    {leftover $[c,x]$: as $y\to x$ the top\\[-2pt] piece keeps a root edge of\\[-2pt] length $\frac{x-c}{1-c}>0$, so the value\\[-2pt] never reaches the basepoint};
  \node[scale=.6,align=center,anchor=north] at (0.28,-0.08) {floor $c<x$:\\[-2pt] \emph{not} a cut system};
\end{scope}
\end{tikzpicture}
\caption{The passage from Ching's windows to the levelled ones. Both legs are continuous; the panels are snapshots, with $w=\tfrac34$ an arbitrary sample. The top window's floor is $\max(x,c)$, so it sits at the vertex until $c$ overtakes it. Leg $1$ raises that floor to $\tfrac12$; leg $2$ lowers the bottom ceiling from $1$ to $\tfrac12$, at which point the windows abut. The middle panel is the common value $H^{(1)}_{1/2} = H^{(2)}_{1}$. The panel at right is not a stage of the path but the family the axioms forbid.}
\label{fig:legs}
\end{figure}

Axiom (W3) of Definition \ref{def:WCS} is the disjunction \emph{anchored or stacked}. The cells of $\Kop$ are the legal rules, composed by substitution at ambient heights, which involves no division and is therefore everywhere defined and strictly associative. The fibres are contractible by Proposition \ref{prop:K-operad}, and both prescriptions occur among the cells:
\[ c^{\mathrm{tr}}_{\alpha_i} \;=\; \tau \text{ at } o(\alpha_i), \quad\quad\quad c^{\mathrm{lev}}_{\alpha_i} \;=\; \tau \text{ at } \sigma(\alpha_i, e), \quad e = \tfrac12.\]
Here $\tau$ is the tautological algebra of Theorem \ref{thm:taut}, and the identifications are Examples \ref{ex:frames} and \ref{ex:Wcells} together with Lemma \ref{lem:equal-evals}. Theorem \ref{thm:main} says that the two operad structures are restrictions of this one algebra along $o$ and $\sigma$.

An explicit path between the two cells witnesses the contractibility. The parameters $c$ and $w$ are heights, but they are independent of the tree: $c$ rises from $0$ to $\tfrac12$, and $w$ then descends from $1$ to $\tfrac12$. Cutting at $\max(x,c)$ keeps the cut at the vertex until $c$ overtakes it. For $c\in[0,\tfrac12]$ and then $w\in[1,\tfrac12]$ set
\begin{align}
 H^{(1)}_{c}(x,y) &= \Big( x, \; \frac{y - \max(x,c)}{1-\max(x,c)} \Big), && \text{valid for } y > \max(x,c), \label{eq:leg1-3}\\
 H^{(2)}_{w}(x,y) &= \Big( \frac{x}{w}, \; \frac{y - \max(x,\tfrac12)}{1-\max(x,\tfrac12)} \Big), && \text{valid for } x < w, \ y > \max(x,\tfrac12), \label{eq:leg2-3}
\end{align}
so that $H^{(1)}_0 = c^{\mathrm{tr}}_{\alpha_i}$, $H^{(1)}_{1/2} = H^{(2)}_{1}$ and $H^{(2)}_{1/2} = c^{\mathrm{lev}}_{\alpha_i}$. Leg $1$ raises the top floor from $x$ to $\tfrac12$, held anchored by the $\max$; leg $2$ lowers the bottom ceiling from $1$ to $\tfrac12$, ending stacked. Splitting the motion in two is only for legibility: the straight line from the first cell to the last would serve as well, since both are anchored, and one checks directly that along it the top window shrinks only as $x\to 1$, where the bottom piece degenerates. The family $H$ above, obtained by writing $c$ where $\max(x,c)$ belongs, is not available. On the worked input the outputs run $(\tfrac14,\tfrac23)\rightsquigarrow(\tfrac14,\tfrac12)\rightsquigarrow(\tfrac12,\tfrac12)$.

\begin{remark}
Arity $3$ does not exhibit the coherence. Equivariance is free here: $\Sigma_3$ permutes the triangles through the labels $\overline{i}$ and fixes all heights, and every formula above uses only heights and which blocks exist, so nothing chooses an orbit representative. But with a single chain in play there is nothing to cohere, whereas Theorem \ref{thm:main} asserts that paths like the one above can be chosen compatibly with all composites in all arities at once. That rests on the contractibility of the fibre of $\Kop$ over every flag (Proposition \ref{prop:K-operad}) and the strict associativity of substitution (Definition \ref{def:K}); what is exhibited above is one fibre, two cells in it, and one path between them.
\end{remark}
\bibliographystyle{plain}
\bibliography{biblo}


\end{document}